\documentclass[11pt] {article} 
\usepackage[backref]{hyperref}
  \usepackage{amsmath}
    \usepackage{amssymb}
    
    \usepackage{pdfsync}

\newtheorem{example}{Example}[section]
\newtheorem{theorem}{Theorem}[section]
\newtheorem{lemma}{Lemma}[section]
\newtheorem{corollary}{Corollary}[section]

\newtheorem{remark}{Remark}[section]

\newcommand{\eqnsection}{
   \renewcommand{\theequation}{\thesection.\arabic{equation}}
   \makeatletter
   \csname @addtoreset\endcsname{equation}{section} 
   \makeatother}
   
\def \ov{\overline}

\def \be{\begin{equation}}
\def \ee{\end{equation}}
\def \bt{\begin{theorem}} 
\def \et{\end{theorem}} 
\def \bl{\begin{lemma}} 
\def \el{\end{lemma}}
\def \bea{\begin{eqnarray}}
\def \eea{\end{eqnarray}}
\def \bas{\begin{eqnarray*}}
\def \eas{\end{eqnarray*}}

\def \al{\alpha}  
\def \bb{\beta}
\def \ga{\gamma} 
 
\def \de{\delta} 
\def \De{\Delta} 
\def \ep{\epsilon}

\def \la{\lambda}  

\def \ka{\kappa}
\def \om{\omega}
\def \Om{\Omega}
\def \va{\varrho}

\def \si{\sigma}

\def \th{\theta}
\def \vth{\vartheta}

\def \ze{\zeta}

\def \ff{\infty}
\def \wh{\widehat}
\def \wt{\widetilde}

\def \rar{\rightarrow}

\def \cd{\,\cdot\,}

\def \AA{{\cal A}}
\def \BB{{\cal B}}

\def \DD{{\cal D}}
\def \EE{{\cal E}}
\def \FF{{\cal F}}

\def \HH{{\cal H}}
\def \II{{\cal I}}

\def \LL{{\cal L}}

\def \PP{{\cal P}}

\def \SS{{\cal S}}

\def \VV{{\cal V}}

\def \({\left(}
\def \){\right)}

\def \lc{\left\{}
\def \rc{\right\}}

\def \nn{\nonumber}
\def \Proof{\noindent{\bf Proof $\,$ }}

\def \bc{\begin{center} }
\def \ec{\end{center} }
\def \bs{\begin{slide} }
\def \es{\end{slide} }

\def\square{{\vcenter{\vbox{\hrule height.3pt
        \hbox{\vrule width.3pt height5pt \kern5pt
           \vrule width.3pt}
        \hrule height.3pt}}}}
\def\qed{{\hfill $\square$ \bigskip}}

 \def \Rev({\mbox{\rm Rev}(}

\eqnsection

\begin{document}
\title{Intersection local times, loop soups  and permanental Wick powers}

 \author{Yves Le Jan\,\, \,\, Michael B. Marcus\,\, \,\, Jay Rosen \thanks{Research of   M. B. Marcus and J. Rosen 
  was partially supported by  grants from the National Science Foundation and
PSC CUNY.  In addition, research of J. Rosen   was partially supported by the Simons Foundation.}}
\maketitle

\begin{abstract}   
  Several  stochastic processes related to transient L\'evy processes with potential densities $u(x,y)=u(y-x)$, that need not be symmetric nor bounded on the diagonal, are defined and studied.  
 They are real valued processes on a   space of measures $\VV$  endowed with a metric  $d$.  Sufficient  conditions   are obtained for the continuity of these processes on $(\VV,d)$.   The processes  include $n$-fold   self-intersection local times of transient L\'evy processes and  permanental chaoses, which are    `loop soup  $n$-fold   self-intersection local times' constructed from the loop soup of the L\'evy process.   Loop soups are also used  to define permanental Wick powers, which generalizes standard Wick powers,  a class of  $n$-th  order Gaussian chaoses.   Dynkin type isomorphism theorems are obtained that relate  the various processes. 
 
   Poisson chaos processes are defined and permanental Wick powers are shown to have a Poisson chaos decomposition. Additional properties of  Poisson chaos processes are studied and a martingale extension is obtained for many of the processes described above.

 \end{abstract}

\bibliographystyle{amsplain}

  \section{Introduction} 
   We define and study several  stochastic processes related to transient L\'evy processes with potential densities $u(x,y)=u(y-x)$ that need not be symmetric nor bounded on the diagonal. We are particularly interested in the case when $u(x,y)$ is not symmetric since some of the processes we consider have already  been studied in the case when it is  symmetric. Significantly the results we obtain give the known results  when $u(x,y)$ is  symmetric.
 
 The processes we consider are real valued processes on a   space of measures $\VV$ endowed with a metric $d$.  We obtain sufficient  conditions for the continuity of these processes on 
 $(\VV, d)$.  Specifically, we study $n$-fold   self-intersection local times of transient L\'evy processes and    permanental chaoses, which are  `loop soup   $n$-fold   self-intersection  local times' constructed from the loop soup of the L\'evy process.  We also use loop soups to define permanental Wick powers, which,    are generalizations of  standard Wick powers,  a class of  $n$-th  order Gaussian chaoses.   We develop the concept of Poisson chaos decompositions and describe the Poisson chaos decomposition of  permanental Wick powers. This illuminates  the relationship between  permanental Wick powers and the  permanental chaos processes constructed from   self-intersection local times.   Lastly we define and study the exponential  Poisson  chaos and   show that the processes described above  have a natural extension as martingales.

 \medskip	 Let   
$Y=\{Y_{t},t\in R^{+}\}$ be a  L\'evy process in $R^d$,   $d=1,2$, with characteristic exponent $\bar \ka$, i.e.,
\begin{equation}
E\(e^{i\xi Y_{t}}\)=e^{-t \bar \ka (\xi)}.\label{sl.00}
\end{equation}  
   We assume that for some $\ga\geq 0$,
 \be
  \Big|{1 \over \ga+  \bar\ka(\xi)}\Big|\le  \frac{1 }{ \va_{\al}(|\xi|) }  \label{s7.8q},
\ee
for some  function $ \va_{\al}(|\xi|)$  that is regularly varying at infinity with index 
 \be  \(1-\frac{1}{2n}\)d<\al\le d,\label{8.34} 
 \ee 
 for some $n\ge 2$. (The precise value of $n$ depends on the results we prove.)
 
    We assume that 
    \begin{equation}
   \int_{K}  \(\va_{\al}(|\xi|)\)^{-1} \,d\xi  <\ff  \qquad \mbox{for all compact sets $K\in R^{d} $.}\label{1.4w}
   \end{equation}
   and that
   \begin{equation}
      \int_{R^{d}}  \(\va_{\al}(|\xi|)\)^{-1} \,d\xi=\ff.\label{1.5}
   \end{equation} 
 We  also assume    that $ \bar\ka $ satisfies the sectorial condition,  
\begin{equation}
 |\mbox{Im } \bar\ka(\xi)| \le C (\ga+ \mbox{ Re }\bar\ka(\xi)) \qquad\forall \xi\in R^{d},\label   {8.21qwa}
  \ee 
  for some $C<\ff$.
  
   Let  
$X=\{X_{t},t\in R^{+}\}$    be the process obtained by killing $Y$ at an independent exponential time with mean $1/\ga$,  $\ga>0$. When $\ga=0$,    we take $X=Y$. In the rest of this paper we simply refer to the transient L\'evy process $X$ without specifying whether or not it is an exponentially killed L\'evy process.

   We  first
 consider the  $n$-fold   self-intersections of    $X $ in $R^{d}$, $d=1,2$. This entails studying  functionals of the form
\be
 \al_{n,\ep}(\nu,t)
\stackrel{def}{=}\int \!\!\int_{[0,t]^n} 
\prod_{j=1}^n f_{\ep}(X(t_j)- y)\,dt_1\cdots\,dt_n
\,d\nu(y),\label{3}
\ee  where $f_{\ep}$ is an approximate $\delta$-function at zero and $\nu$ is a finite  measure on $R^{d}$.
 Ideally we would like to take the limit   of $ \al_{n,\ep}(\nu,t)$ as $\ep$ goes to 0, but  if the potential density of $X$ is unbounded at the origin, which is always the case in dimension $d\geq 2$,  the limit is
infinite   for all  $n\ge 2$. To  deal with this we use a technique  called renormalization, which consists of
forming a
linear combination of  the  
$\{\al_{k,\ep}(\nu,t)\}_{k=1}^n$ which has a finite limit, $ L_{n}(\nu,t)$, as $\ep\to 0$.   We study  the behavior of  $L_{n}(\nu):= L_{n}(\nu,\ff)$.

We set 
$L_{1,\ep}(\nu)=\al_{1,\ep}(\nu,\ff)$ and define   recursively  
\begin{equation}
L_{n,\ep}(\nu)= \al_{n,\ep}(\nu,\ff)-\sum_{j=1}^{n-1}c_{n,j,\ep}L_{j,\ep}(\nu)\label{zzx}
\end{equation}
where the $c_{n,j,\ep}$ are constants which diverge as $\ep\rar 0$; see (\ref{rilt.1})   and (\ref{2.47}).  

 We introduce  a $\si$-finite measure $\mu$ on the path space of $X$, called the loop measure, and show that for a certain class of positive measures $\nu$
\begin{equation}
L_{n}(\nu):=\lim_{\ep \rar 0}L_{n,\ep}(\nu)\quad\mbox{    exists in  $L^{p}(\mu)$, for all $p\ge 2$.}\label{rilt.9intro}
\end{equation}
We refer to $L_{n}(\nu)$ as the $n$-fold renormalized self-intersection
local time   of ${ X}$, with respect to   $\nu $. 
\medskip

 We show in the beginning of Section \ref{sec-iltac} that $X$ has potential densities $u(x,y)=u(y-x)$. For $\ga\ge 0$ we define 
   \begin{equation}
\hat u(\xi)  ={1 \over \ga+  \bar\ka(\xi)}.\label{sl.12e}
\end{equation}
By (\ref{s7.8q})
 \be
  |\hat u(\xi)|\le  \frac{1 }{ \va_{\al}(|\xi|) }  \label{7.8q}.
 \ee    
  Let  
\begin{equation}
\tau_{n}(\xi):= \overset{\mbox{$n$-times}}
 {\overbrace {|\hat u|\ast  |\hat u|\ast \cdots  \ast  |\hat u|}}  (\xi)\label{l2.1}
\end{equation}
denote the $n$-fold convolutions of   $|\hat u|$.   We  define
 \be
\|\nu \|_{2,\tau_{2n}}:=\(\int \tau_{2n}(\la )|\hat{\nu}(\la)|^{2}\,d\la\)^{1/2}<\ff.\label{12.30mm}
 \ee
 
  Let $f$ and $g $ be functions on $R_{+}^{1}$. By the rotation invariance of Lebesgue measure on $R^{d}$, $\int_{R^{d}}f(|\eta-\xi|)g(|\eta|)\,d\eta$ depends only on $|\xi|$. We let $f\underset{d}{\ast}g$ denote the function on $R_{+}^{1}$ which satisfies
\begin{equation}
f\underset{d}{\ast}g(|\xi|)=\int_{R^{d}}f(|\eta-\xi|)g(|\eta|)\,d\eta.\label{dconv}
\end{equation}
We refer to $f\underset{d}{\ast}g$ as the d-convolution of $f$ and $g$.   (The letter $d$ refers to the dimension of the space we are integrating on.)

Let $\wt\tau_{n}$ denote the   $n$-fold d-convolution  of   $(\va_{\al} )^{-1}$. (We use the notation the notation $(\va_{\al} )^{-1}$ for $1/\va_{\al})$. That is,  
\begin{equation}\qquad 
\wt\tau_{n}(|\xi|):= \overset{\mbox{$n$-times}}
 {\overbrace { (\va_{\al} )^{-1} \underset{d}{\ast}  (\va_{\al} )^{-1}\underset{d}{\ast}\cdots  \underset{d}{\ast} (\va_{\al} )^{-1} }}  (|\xi|).\label{l2.1f}
\end{equation} 
It follows from (\ref{8.34}) and Lemma \ref{lem-9.2} that    $\wt\tau_{2n}(|\xi|)<\ff$ and $\lim_{|\xi|\to\ff}\wt\tau_{2n}(|\xi|)=0$.

Let $  \mathcal {B}_{2n}(R^{d})$ denote the set of finite   signed measures 
 $\nu$ on $R^{d}$  such that  
 \be
\|\nu \|_{2,\wt\tau_{2n}}:=\(\int \wt\tau_{2n}(|\la |)|\hat{\nu}(\la)|^{2}\,d\la\)^{1/2}<\ff.\label{12.30m}
 \ee

 \bt\label{theo-multiriltintro-m} Let
$X=\{X(t),t\in R^+\}$ be a L\'evy process in $R^d$, $d=1,2$,   as described above  
   and let $\nu\in \mathcal {B}_{2n}(R^{d})$,     $n\ge 1$. Then (\ref{rilt.9intro}) holds.
   \et

    We are   also concerned with the continuity of  $\{L_{n}(\nu),\nu\in \VV\}$, where $\VV$ is some metric space. Here is a particularly straight forward example of our results. 
    
\medskip	      For any finite positive measure $\nu$ on
$R^d$, let $\nu_{x}(A)=\nu(A-x)$.   

\bt\label{theo-cont14.1b}  Under the hypotheses of Theorem  \ref{theo-multiriltintro-m} 
\be
\int_{1}^{\ff}  \(\int_{|\xi|\geq x}\tau_{2n}(\xi)|\hat{\nu}(\xi)|^{2}\,d\xi\) ^{1/2} \frac{
 (\log x)^{n-1}}{x}\,dx<\ff,\label{b101}
\ee   
is a sufficient condition for $\{L_{n}(\nu_x),\,x\in R^d\}$ to be   continuous  $P^{y}$ almost surely, for all $y\in R^{d}$. \et (As usual, $P^{y}$ denotes the probability of the L\'evy process $X$ starting at $y\in R^{d}$.)

Several concrete examples are given at the end of Section \ref{sec-iltcont}. As a sample, we note that  when (\ref{7.8q}) holds for $1/\va_{\al}  (|\xi|) =O(|\xi|^{-d })$ as   $|\xi|\to\ff$, $\{L_{n}(\nu_x),\,x\in R^d\}$ exists for all $n\ge 2$, and is continuous almost surely when 
\begin{equation}
  |\hat\nu(\xi)|^{2}= O\( \frac{1}{(\log |\xi|)^{4n   +\de}}\)   \qquad\mbox{as $|\xi|\to\ff$},
   \end{equation} 
for any $\de>0$.  

  When $\va_{\al}  (|\xi|) =|\xi|^{\al}$, for $d(1-\frac{1}{2n})<\al<d$,  $\{L_{n}(\nu_x),\,x\in R^d\}$ exists  and is continuous almost surely when 
\be
 |\hat\nu(\xi)|^{2}=O\(\frac{1}{|\xi|^{2n(d-\al) }  (\log |\xi| )^{2n +1 +\de}} \)    \qquad\mbox{as $|\xi|\to\ff$},
   \end{equation} 
for any $\de>0$.

\medskip	
   Theorems \ref{theo-multiriltintro-m} and \ref{theo-cont14.1b}     follow easily from the  proof of the next theorem.  

 \begin{theorem}   \label{theo-1.3}     Let
$X $ be as in Theorem  \ref{theo-multiriltintro-m} 
and let $n=n_{1}+\cdots+n_{k}$,      $k\ge 2$,  and $ \nu_{i}\in \BB_{2n_{i}}(R^{d})$. Then 
 \bea 
\mu\(  \prod_{i=1}^{k} L_{n_{i}}(\nu_{i} ) \)\label{rilt.15intro}&
=&{\prod_{i=1}^{k}(n_{i}!) \over n} \sum_{\pi\in \mathcal{M}_{a}}  \int  \prod_{j=1}^{n}u (x_{\pi(j)},x_{\pi(j+1)}) \prod_{i=1}^{k}\,d\nu_{i}  (x_{i} )\nn\\ 
&\le& \frac{| \mathcal{M}_{a}|}{n} \prod_{i=1}^{k}n_{i}! \|\nu_{i}\|_{2, \tau_{2n_{i}}},
\eea 
 where $ \pi(n+1) = \pi(1) $ and $\mathcal{M}_{a}$ is the set of maps $\pi:[1,n]\mapsto [1,k] $  with
 $|\pi^{-1}(i)|=n_{i}$ for each $i$  and  such that, if $\pi(j)=i$ then $\pi(j+1)\neq  i$. (The subscript `a'  in $\mathcal{M}_{a}$ stands for alternating).
 \et
 
  (Note that although the hypothesis of Theorem \ref{theo-1.3} requires that $\|\nu_{i}\|_{2, \wt\tau_{n_{i}}}<\ff$, the bound in (\ref{rilt.15intro}) is in terms of  the possibly smaller  norms $\|\nu_{i}\|_{2, \tau_{2n_{i}}}$.)

\medskip

In \cite{MRmem} we study   self-intersection local times of L\'evy processes with symmetric potential densities. We obtain sufficient conditions for continuity, such as Theorem \ref{theo-cont14.1b} above,  by associating the   self-intersection local time $ L_{n}(\nu)$ with a $2n$-th order Gaussian chaos $:G^{2n}\nu:$,    called a $2n$-th Wick power, that is constructed from the Gaussian field   with covariance 
\begin{equation}
E\(   :G^{2n}\nu::G^{2n}\mu:\)=(2n)!\int\!\!\int\(u(x,y)\)^{2n}\,d\nu(x)\,d\mu(y).\label{1.9}
   \end{equation}
(In (\ref{1.9}) $u$ must be symmetric.) The association is by a Dynkin type isomorphism theorem that allows us to infer  results
about  $ \{L_{n}(\nu),\nu\in\VV\}$ from results about  $\{:G^{2n}\nu:,\nu\in\VV\}$.
The advantage here is that the continuity results  we consider are known for  the  Gaussian chaoses,   so once we have the isomorphism theorem it is easy to extend them to the associated intersection local times.

In this paper     $u$ need not be symmetric. We define and obtain continuity results about   self-intersection local times directly, without relating them   to any other stochastic process. Nevertheless, the question remains, is there a Dynkin type isomorphism theorem that relates them  to another process and more specifically to what other process. We give two answers to this question. In the first we relate the  intersection local times to the  `loop soup   $n$-fold   self-intersection local time' $\psi_{n}(\nu)$, which  we also call an $n$-th order permanental chaos. To construct $\psi_{n}(\nu)$ we take a   Poisson  process  $\LL_{\al}$ with intensity measure       $\al\mu$,  $\al>0$, on $\Om_{\De}$, the space of the paths of $X$. This process is called a loop soup. The loop soup   self-intersection local time $\psi_{n}(\nu)$ is the   renormalized sum of the  $n$-fold   self-intersection local times, i.e., the $L_{n}(\nu)$, of the paths in $\LL_{\al}$. In   Theorem  \ref{theo-ljriltintro} we give a Dynkin type isomorphism theorem that relates $\{L_{n}(\nu),\nu\in\VV\}$ and $\{\psi_{n}(\nu),\nu\in\VV\}$.

\medskip	  Analogous to (\ref{rilt.15intro}) we have the following joint moment formula for   $\psi_{n_{1}}(\nu),\ldots \newline \psi_{n_{k}}(\nu)$, which is proved in Section  \ref{sec-loopilt}:

  \begin{theorem} \label{theo-1.4} Let
$X $ be as in Theorem  \ref{theo-multiriltintro-m} and 
let $n=n_{1}+\cdots+n_{k}$ and $c(\pi)$ equal to the number of cycles in the permutation $\pi$. Then   
 \be 
E_{\mathcal{L}_{\al}}\( \prod_{i=1}^{k}\psi_{n_{i}}(\nu_{i}) \)\label{20.15qintro}= \sum_{\pi\in \mathcal{P}_{0}}\al^{c(\pi)}  \int  \prod_{j=1}^{n}u (z_{j},z_{\pi(j)}) \prod_{i=1}^{k}\,d\nu_{i}  (x_{i} ) ,
\ee 
 where $z_{1},\ldots,z_{n_{1}}$ are all equal to $x_{1}$,  the next $n_{2}$ of the $\{z_{j}\}$ are all equal to $x_{2}$, and so on, so that the last   $n_{k}$ of the $\{z_{j}\}$ are all equal to $x_{k}$    
and  $\mathcal{P}_{0}$ is the set of permutations $\pi$ of $[1,n]$ with cycles that alternate the variables $\{x_{i}\}$; (i.e., for all $j$, if $z_{j}=x_{i}$ then $z_{\pi(j)}\neq x_{i}$), and in addition, for each $i=1,\ldots,k$, all the $ \{z_{j}\}$ that are equal to $x_{i}$ appear in the same cycle.  
 \end{theorem}

\medskip	 Using standard results about Poisson processes, in Section  \ref{sec-IT}, we   obtain an isomorphism theorem that relates the   self-intersection local times $L_{n}(\nu)$ of $X$ and the loop soup   self-intersection local times $\psi_{n}(\nu)$ of $X$.   
\bt[Isomorphism Theorem  I]  \label{theo-ljriltintro}  For any   positive measures $\rho,\phi\in  \BB  _{2}(R^{d})$ there exists a random variable $\th^{\rho,\phi}$ such that for any finite measures $\nu_{j} \in \BB_{2n_{j}}(R^{d})$ ,  $j=1,2,\ldots$, and bounded 
measurable functions  $F$ on $R^\ff_+$,  
\begin{equation}
E_{\mathcal{L}_{\al}}\int Q^{x,x}\(L_{1}(\phi)\, F \( \psi_{n_{j}}(\nu_{j})+L_{n_{j}}(\nu_{j})\) \)\,d\rho(x)={1 \over \al}E_{\mathcal{L}_{\al}} \(\th^{\rho,\phi} \,F\(\psi_{n_{j}}(\nu_{j})\)\).\label{6.1intro}
\end{equation} 
  (Here we use the notation $F( f(x_{ i})):=F( f(x_{ 1}),f(x_{ 2}),\ldots 
)$.)
\et

  It is interesting to compare (\ref{6.1intro}) with    \cite[Theorem 1.3]{LMR} in which  all $n_{j}=1$, and with \cite[Theorem 3.2]{EK} which is for local times.

\medskip	   The measure $Q^{x,y}$, which is used to define the loop measure $\mu$, is defined in   (\ref{10.1q}).    The term   $\th^{\rho,\phi} $ is a positive random variable with $E_{\mathcal{L}_{\al}} \((\th^{\rho,\phi}) ^{k}\)<\ff$  for all integers $k\ge 1$. In particular 
 \begin{equation}
E_{\mathcal{L}_{\al}} \(\th^{\rho,\phi} \)=\al \int Q^{x,x}\(L_{1}(\phi) \)\,d\rho(x). \label{6.2}
 \end{equation}
 
It is actually quite simple to obtain Theorems  \ref{theo-1.4} and \ref{theo-ljriltintro}. However they are not really generalizations of the results in \cite{MRmem}. The loop soup   self-intersection local time $\psi_{n}(\nu)$ is not a $2n$-th Wick power when the potential density  of $X$ is symmetric.  

 	 Our second answer to the questions raised in the paragraphs preceding Theorem  \ref{theo-1.4} 
 relates  the   self-intersection local times $L_{n}(\nu)$ of $X$ with a process $\wt\psi_{n}(\nu)$ which we call a   permanental Wick power.  The   rationale for this name is that this process is a $2n$-th Wick power when   $u$ is symmetric and $\al=1/2$. It seems significant to have a generalization of Wick powers that does not require that the kernel that defines them is symmetric. This is done in Section \ref{sec-powers} in which we give analogues of Theorems    \ref{theo-1.4} and \ref{theo-ljriltintro} for $\wt\psi$. (The analogue of Theorem     \ref{theo-1.4}   for $\wt\psi$ is exactly the same as Theorem     \ref{theo-1.4}, except that the final phrase ``and in addition, for each $i=1,\ldots,k$, all the $ \{z_{j}\}$ that are equal to $x_{i}$ appear in the same cycle.'' is omitted.)   
 
  In Section \ref{sec-Poisson} we develop the concept of Poisson chaos decompositions.  In Section \ref{sec-lspower} we obtain a   Poisson chaos decomposition of   the   permanental Wick power $\wt\psi_{n}(\nu)$  and relate it to the   loop soup self-intersection local time $\psi_{n}(\nu)$. 
The process $\psi_{n}(\nu)$ incorporates the  self-intersection local times $L_{n}(\nu)$ of each path in the loop soup. In addition to the  self-intersection local times, the process $\wt\psi_{n}(\nu)$ also incorporates  the mutual  intersection local times between different  paths   in the loop soup.
 In Theorem \ref{theo-pjr} we give an isomorphism theorem relating permanental Wick powers and   self-intersection local times.
 
 In Section \ref{sec-expP} we  define and study  exponential  Poisson  chaoses and   show that many of  the  processes we consider  have a natural extension as martingales.
 
\medskip	  Here is a summary of the processes we study and a reference to the first place they appear. They are all related to a L\'evy process $X$ and are functions of  the potential density of $X$.

 \bigskip
\begin{tabular}{ll}
  $ \al_{n,\ep}(\nu,t) $ & approximate  $n$-fold   self-intersection local time of   $X$,  \\&(\ref{3}).
   \\&\\
  $L_{n}(\nu) $& $n$-fold renormalized self-intersection local time  (of $X$), \\&(\ref{2.47}), (\ref{rilt.9intro}). \\&\\
 $\mu$ &      loop measure on the path space of $X$,  (\ref{ls.3}).\\&\\
  $\LL_{\al}$&   Poisson point process on the path space of $X$ with \\&intensity measure $\al\mu$ called  the loop soup of $X$, \\& Section \ref{sec-loopilt}. \\&\\
$\psi_{n}(\nu)$& loop soup $n$-fold   self-intersection local time (of $X$). \\&Also  called a permanental chaos. $\psi_{1}(\nu)$ is also called   \\
&a permanental field, (\ref{ls.10a}), (\ref{ls.10}).
 \\&\\
$\wt \psi_{n}(\nu)$&  $n$-th order renormalized permanental field (of $X$). Also \\&called a $2n$-th  permanental Wick power, (\ref{11.2j}),  \\
&Theorems \ref{theo-multivar}  and
\ref{theo-multiloop}.
 \\&\\
$I_{n}(g_{1},\ldots,g_{n})$& Poisson Wick product, page  \pageref{p56}.
\\&\\
$\II_{l_{1},\ldots,l_{k}}(\nu)$& $l_{1}+\cdots+l_{k}=n$, $k$-path, $n$-fold    intersection local \\
&time  (of $X$),  (\ref{10.24}).  \\&\\ $ \oplus _{n=0}^{\ff }\,\, H_{n}$ &  Poisson chaos decomposition   of $L^{2}\(P_{\mathcal{L}_{\al}}\)$,  (\ref{11.0r}).\\&\\
$    \sum  
\II_{|D_{1}|,\ldots, |D_{l}|}(\nu)$&Poisson chaos decomposition   of $\wt\psi_{n}(\nu)$, Theorem  
\ref{theo-multiloop}.
\\&\\
  $\EE(g)$& exponential Poisson chaos (\ref{10.10}).
\\&\\
$\II_{n}^{(\al)} $, $\wt \psi_{n}^{(\al)}$ etc. &$ (E_{\LL}, \FF_{\al}) $ martingales,  page \pageref{cor-mart}.\end{tabular}
 
 \bigskip
Several relationships between these processes are given. For example, it follows immediately from  in Theorem  
\ref{theo-multiloop}  and Corollary \ref{cor-7.1}  that 
\be
\wt\psi_{n}(\nu)=  \psi_{n}(\nu)+\sum_{D_{1}\cup\cdots\cup D_{l}=[1,n];l\ne 1}
\II_{|D_{1}|,\ldots, |D_{l}|}(\nu).
\ee

Critical estimates are used in all the proofs that require understanding properties of convolutions of regularly varying functions. These are studied in  Section \ref{sec-8}.

\section{Loop measures and renormalized intersection local times }\label{sec-iltac}

 Let  
$X=\{X_{t},t\in R^{+}\}$ be a  L\'evy process in $R^d$ as described in the Introduction. 
 It follows from  (\ref{8.21qwa}) that 
   \begin{equation}
 {1 \over \ga+  \mbox{ Re } \bar\ka(\xi)} \le {C' \over |\ga+  \bar\ka(\xi)|}, \qquad\forall \xi\in R^{d},  \label{s77}
   \end{equation}
   for some constant  $C'<\ff$. 
  Together with (\ref{s7.8q}) this shows that for each $t>0$
\begin{equation}
e^{-t \bar \ka (\xi)}\in L^{1}(R^{d}).\label{sl.10}
\end{equation}
We define   the continuous function 
\begin{equation}
p_{t}(x)={1 \over (2\pi)^{d}}\int e^{-ix\xi}e^{-t (\ga+  \bar\ka(\xi))}\,d\xi.\label{sl.11}
\end{equation}
 Note that for any $f\in \SS(R^{d})$,   the space of   rapidly decreasing $C^{\ff}$  functions,  
\begin{eqnarray}
\int p_{t}(x)f(x)\,dx&=& \int \wh f(\xi) e^{-t (\ga+  \bar\ka(\xi))}\,d\xi
\label{ft1}\\
&=& \int \wh f(\xi) E\(e^{i\xi X_{t}}\)\,d\xi \nonumber\\
&=&  E\(f(X_{t})\). \nonumber
\end{eqnarray}
This shows that  $p_{t}(x)$ is a (sub)probability density function for $X_{t}$. In particular it is integrable on $R^{d}$. Therefore,  we can invert the transform in (\ref{sl.11}), and  using  the fact that $\bar \ka (\xi)$ is continuous, we obtain
\begin{equation}
e^{-t (\ga+  \bar\ka(\xi))}=\int e^{ix\xi}p_{t}(x)\,dx.\label{sl.11k}
\end{equation}
We have
\be
\wh {p_{t}\ast p_{s}}(\xi)=\wh p_{t}(\xi)    \wh p_{s}(\xi)=e^{-(t+s) (\ga+  \bar\ka(\xi))}.\label{2.10q}
\ee
  We define $p_{t}(x,y)=p_{t}(y-x)$. It follows from (\ref{2.10q}) that   $\{p_{t}(x,y),(x,y,t)\in R^{1}_{+}\times R^{2d} \}$ is  a jointly continuous semigroup of transition densities for  $X$. We    define  
\begin{equation}
u (x,y) =\int_{0}^{\ff} p_{t}(x,y)\,dt. \label{sl.12a}
\end{equation}
Note that  $u(x,y)$ are potential densities for   $X$.

To justify the definition of $\hat u$ in (\ref{sl.12e}) 
note that for any $f\in \SS(R^{d})$,   
\be
\int p_{t}(x)f(x)\,dx= \int \hat f(\xi)e^{-t (\ga+  \bar\ka(\xi))}\,d\xi,\label{2.8}
\ee
by (\ref{sl.11}).  Consequently, by (\ref{sl.12a}),
  \begin{equation}
\int u (x)f(x)\,dx=\int {1 \over \ga+  \bar\ka(\xi)}\,\hat f(\xi) \,d\xi,\label{sl.12d}
\end{equation} 
in which the use of  Fubini's theorem  is justified by (\ref{s77})   and (\ref{s7.8q}). In this context,  $ ( \ga+  \bar\ka(\xi))^{-1}$ is the Fourier transform of $u$, considered  as a distribution in $\SS'$.

 Let  $X$ be a  L\'evy process in  $R^{d}$   with  transition densities $p_{t}(x,y)=p_{t}(y-x)$ and  potential densities     $u(x,y)=u(y-x)$ as described above.  We assume that $u(x)<\ff$ for $x\ne 0$, but since we are interested in L\'evy processes that do not have local times we are primarily concerned with the case when $u(0)=\ff$. We also assume that $\int_{\de}^{\ff}p_{t}(0)\,dt<\ff$ and that $0<p_{t}(x)<\ff$.  Most significantly, we do not require that   $p_{t}(x,y)$,   and therefore that $u(x,y)$, is  symmetric. 

Under   these assumptions  it  follows, as in  \cite{FPY}, that for all $0<t<\ff$ and  $x,y\in R^{d}$,   there exists a finite measure $Q_{t}^{x,y}$ on $\mathcal{F}_{t^{-}}$, of total mass $p_{t}(x,y)$, such that
\begin{equation}
Q_{t}^{x,y}\(1_{\{\ze>s\}}\,F_{s}\)=P^{x}\( F_{s}  \,p_{t-s}(X_{s},y) \),\label{10.1}
\end{equation}
for all $F_{s} \in \mathcal{F}_{s}$,    $s<t$.

We   take  $\Om$ to be  the set of right continuous paths $\om$ in  $R^{d}_{\De}=R^{d}\cup \De$ where  $\De\notin R^{d}$, and is such that  $\om_{t}=\De$ for all \label{page10}$t\geq \ze=\inf \{t>0\,|\om_{t}=\De\}$.   We set  $X_{t}(\om)=\om_{t}$ and define a $\si$-finite measure $\mu$
on $(\Om, \mathcal{F})$ by the following    formula:
\begin{equation}
\int F\,d\mu  =\int_{0}^{\ff}{1\over t}\int\,Q_{t}^{x,x}\(F\circ k_{t} \)\,dm (x)\,dt, \label{ls.3}
\end{equation}
for all $\mathcal{F}$-measurable functions $F$ on $\Om$.
Here $k_{t}$ is the killing operator defined by $k_{t}\om(s)=\om(s)$ if $s<t$ and $k_{t}\om(s)=\De$ if $s\geq t$, so that $k_{t}^{-1}\mathcal{F}\subset \mathcal{F}_{t^{-}}$. (We  often write $\mu (F)$ for $\int F\,d\mu$.)

 \medskip  The  next lemma is  \cite[Lemma 2.1]{LMR}, (with   $\nu_{j}(dx)=g_{j}(x)\,dx$, $j=1,\ldots,k$).

\begin{lemma}\label{lem-ls} Let $k\geq 2$ and  $g_{j}$, $ j=1,\ldots,k$ be bounded   integrable  functions on $R^{d}$. Then the loop measure $\mu$ defined in (\ref{ls.3}) satisfies  
\bea
\lefteqn{
\mu\(   \prod_{j=1}^{k}\int_{0}^{\ff}g_{j}(X_{t})\,dt\) \label{ls.4}}\\
&& = \sum_{\pi\in \mathcal{P}^{\odot}_{k }}\int u(y_{\pi(1)},y_{\pi(2)})\cdots   u(y_{\pi(k-1)},y_{\pi(k)})u(y_{\pi(k)},y_{\pi(1)})  \prod_{j=1}^{k} g_{j}(y_{j})\,dy_{j}\nn
\eea
where $\mathcal{P}_{k}^{\odot}$ denotes the set of permutations of $[1,k]$ on the circle.   (For example,  $(1,2,3)$,  $(3,1,2)$ and $(2,3,1)$ are considered to be one permutation   $\pi\in \mathcal{P}^{\odot}_{3 }$.)
\end{lemma}

 	\noindent 
  	We show in (\ref{ls.2mm})   that  when  $u(0)=\ff$,
\be
\mu\(\int_{0}^{\ff}g_{j}(X_{t})\,dt\)=\infty.\label{2.4}
\ee

\subsection{Renormalized intersection local times }\label{subsec-rilt}

  Let $f(y)$ be a positive smooth function supported in the unit ball of $R^{d}$ with $\int f(x)\,dx=1$. Set  $f_{r}(y)=r^{-d}f(y/r)$,  and $f_{ x,r}(y)=f_{r}(y-x)$. 
Let 
\begin{equation}
L(x,r):=\int_{0}^{\ff}f_{ x,r}(X_{t})\,dt. \label{rilt.1q}
\end{equation} 
  $L(x,r)$ can be thought of as the approximate total  local time of $X$ at the point $x\in R^{d}$.  When $u(0)=\ff$, the local time of $X$ does not exist and we can not take the limit $\lim_{r\rar 0}L(x,r)$.    Nevertheless, it is often the case that renormalized intersection local times exist.    We proceed to define renormalized intersection local times.
  
 We begin with the definition of the chain functions 
  \begin{equation}
  \mbox{ch}_{k}(r) =  \int   u (ry_{1},ry_{2}) \cdots u(ry_{k},ry_{k+1})   \prod_{j=1}^{k+1} f(y_{j}) \,dy_{j},\hspace{.2 in} k\geq 1.\label{1L.1}
  \end{equation} 
  Note that $\mbox{ch}_{k}(r)$ involves $k$ factors of the potential density $u$,   but $k+1$ variables of integration.
    For any $  \si=(k_{1}, k_{2},\ldots )$ let   
      \begin{equation}
 |\si|=\sum_{i=1}^{\ff}  ik_{i}\qquad\mbox{and}\qquad  |\si|_{+}=\sum_{i=1}^{\ff}   (i+1)k_{i}.\label{rilt.2a}
  \end{equation}
  We   set   $L_{1}(x,r)=L(x,r)$ and define  recursively 
  \begin{eqnarray}
&&
L_{n}(x,r)  = L^{n} (x,r)  -\sum_{\{ \si\,|\, 1\leq | \si| < | \si|_{+} \leq n\} } J_{n}( \si,r) \label{rilt.1},
  \end{eqnarray}
  where
\begin{equation}
J_{n}( \si,r)={n! \over \prod_{i=1}^{\ff}k_{i}! (n-| \si|_{+})!}\prod_{i=1}^{\ff} \(\mbox{ch}_{i}( r)\)^{k_{i} } 
    L_{n-|\si|} (x,r) .\label{rilt.2b}
\end{equation} 
(Note that $n-|\si|_{+} \ge 0$.)

To help in understanding (\ref{rilt.1}) we note that 
\begin{equation}
L_{2}(x,r)  = L^{2} (x,r)  -2 \,  \mbox{ch}_{1}(r)L(x,r)\label{rilt.2ex}
\end{equation}  
and 
\bea
L_{3}(x,r) & =& L^{3} (x,r)  -6\,  \mbox{ch}_{1}(r)L_{2}(x,r)-6\mbox{ch}_{2}(r)L(x,r)\label{rilt.3ex}\\
& =&  L^{3} (x,r)  -6\,   \mbox{ch}_{1}(r)L^{2}(x,r)+(12\,  \mbox{ch}^{2}_{1}(r)-6\, \mbox{ch}_{2}(r))L(x,r).\nn
\eea

 We show in Remark \ref{rem-2.2} that $L_{n}(x,r)$ can also be defined by a generating function.
     
\medskip	\noindent  {\bf  Proof  of  Theorems \ref{theo-multiriltintro-m} and \ref{theo-1.3} }   Set  
\be L_{n,\ep}(\nu)=  \int   L_{n}(x,\ep) \,d\nu(x).\label{2.47}
\ee We show that
 for $\nu\in \BB_{2n}(R^{d})$, $d=1,2$ and  $n\ge 1$,
\begin{equation}
L_{n}(\nu) =\lim_{\ep \rar 0}L_{n,\ep}(\nu)\quad\mbox{    exists in  $L^{p}(\mu)$, for all $p\ge 2$.}\label{rilt.9introq}
\end{equation}
  The techniques  used in the proof  (\ref{rilt.9introq}) allow us to
show that for   $n=n_{1}+\cdots+n_{k}$,      $k\ge 2$,  and $ \nu_{i}\in \BB_{2n_{i}}(R^{d})$,
 \bea 
\mu\(  \prod_{i=1}^{k} L_{n_{i}}(\nu_{i} ) \)\label{rilt.15introww}&
=&{\prod_{i=1}^{k}(n_{i}!) \over n} \sum_{\pi\in \mathcal{M}_{a}}  \int  \prod_{j=1}^{n}u (x_{\pi(j)},x_{\pi(j+1)}) \prod_{i=1}^{k}\,d\nu_{i}  (x_{i} )\nn\\ 
&\le& \frac{| \mathcal{M}_{a}|}{n} \prod_{i=1}^{k}n_{i}! \|\nu_{i}\|_{2, \tau_{2n_{i}}},
\eea 
 where $ \pi(n+1) = \pi(1) $ and $\mathcal{M}_{a}$ is the set of maps $\pi:[1,n]\mapsto [1,k] $  with
 $|\pi^{-1}(i)|=n_{i}$ for each $i$  and  such that, if $\pi(j)=i$ then $\pi(j+1)\neq  i$.

We begin by showing that for    approximate identities $f_{r,x}$  and    $f_{i}=f_{r_{i}, y_{i}},i= 1,\ldots,m$,     
  \bea
\lefteqn{
\mu\(L_{n }(x,r)   \prod_{i=1}^{m}\int_{0}^{\ff}f_{i}(X_{t})\,dt\)\label{20.15s}}\\
&=& \sum_{\pi\in \mathcal{P}^{\odot}_{m,n }} \int  \prod_{j=1}^{m+n}u (z_{\pi(j)},z_{\pi(j+1)})     \prod_{i=1}^{m} f_{i}  (z_{i} )\,dz_{i} \prod_{i=m+1}^{m+n}f_{r,x}  (z_{i} )\,dz_{i} \nn\\
&&\hspace{1 in} +  \int E_{r}(x,  {\bf z })\prod_{i=1}^{m} f_{i}  (z_{i} )\,dz_{i} \prod_{i=m+1}^{m+n}f_{r,x}  (z_{i} )\,dz_{i} \nn,
\eea
  where 
  $ \pi(m+n+1) = \pi(1) $ and  $\mathcal{P}^{\odot}_{m,n }$ is the subset of permutations $\pi\in \mathcal{P}^{\odot}_{m+n }$  with the property   that for all $j$, when $\pi(j)\in [m+1,m+n]$, then $\pi(j+1)\in [1,m]$.    That is, under the permutation $\pi$, no two elements of $[m+1,m+n]$ are adjacent,   mod $m+n$.   When $m<n$, $\mathcal{P}^{\odot}_{m,n }$ is empty.

   The last term in (\ref{20.15s}) is an error term.  It is actually the sum of many terms, some of which may depend on some of the $z_{1},\ldots,z_{n+m}  $.  We use ${\bf z }$ to designate $z_{1},\ldots,z_{n+m}  $.    Since the functions    $f_{\cd}$  are  probability 
density functions we write  last term in (\ref{20.15s})  as an expectation,
 \be 
\EE_{{\bf f }}(E_{r}(x,  {\bf z })):=\int E_{r}(x,  {\bf z })\prod_{i=1}^{m} f_{i}  (z_{i} )\,dz_{i} \prod_{i=m+1}^{m+n}f_{r,x}  (z_{i} )\,dz_{i}.\label{2.13}
\ee
 We show  in Section \ref{sec-iltan} that   for   $\nu\in\BB_{2n}(R^{d})$,  
  \begin{equation}
 \lim_{r\rar 0} \sup_{\forall |z_{i}|\le M}\int E_{r}(x, {\bf z }) \,d\nu(x)=0,\label{err0}
  \end{equation}
 for some finite number M, which implies that 
  \begin{equation}
    \lim_{r\rar 0}\int  \EE_{{\bf f }}(E_{r}(x,  {\bf z }))\,d\nu(x)=0.
   \end{equation}
 We are using the fact that  since $f$ is supported on the unit ball in $R^{d}$, $f_{r,x}$  has compact support, so that  the  variables $|z_{i}|$ are uniformly bounded.) We deal with all the additional  error terms that are introduced similarly.

The integral on the right-hand side of (\ref{20.15s}) is the same for permutations $\pi$ and $\pi'$ that  differ by a permutation of $\{m+1,\ldots, m+n\}$. We call these internal permutations. (For example, suppose $m=3$ and $ n=2$ and $\pi=(1,4,3,2,5)$ and $\pi'=(1,5,3,2,4)$.  The permutations $\pi$ and $\pi'$differ by a permutation of $\{4,5\}$).  There are $n!$ internal permutations   in $\mathcal{P}^{\odot}_{m,n }$.   This accounts for the  factor $n!$ in   (\ref{rilt.2b}).

\medskip	For a function $g(x), x\in R^{d}$ we define $\De_{h}g(x)=g(x+h)-g(x)$.
  Note that since $u(x,y)=u(y-x)$
  \begin{equation}
   \De_{h}u(x,y)=u(x,y+h)-u(x,y)=u(x-h,y)-u(x,y).
   \end{equation}
By (\ref{ls.4}) and   (\ref{rilt.1q})
   \bea
&&
\mu\(L^{2 }(x,r)   \prod_{i=1}^{m}\int_{0}^{\ff}f_{i}(X_{t})\,dt\)\label{20.15s2}\\
&&\qquad= \sum_{\pi\in \mathcal{P}^{\odot}_{m+2}}  \int  \prod_{j=1}^{m+2}u (z_{\pi(j)},z_{\pi(j+1)})     \prod_{i=1}^{m}  f_{i}  (z_{i} )\,dz_{i}\prod_{i=m+1}^{m+2}f_{r,x} (z_{i } ) \,dz_{i} \nn,
\eea
where $ \pi(m+3) = \pi(1) $. Considering  (\ref{rilt.2ex}) it is clear that   to obtain (\ref{20.15s}) for $n=2$ we need only show that  
\begin{equation}
2 \,  \mbox{ch}_{1}(r) \mu\(L(x,r) \prod_{i=1}^{m}\int_{0}^{\ff}f_{i}(X_{t})\,dt\)
   \end{equation}
 is equal to   the second line in  (\ref{20.15s2}), except that the sum is taken over permutations $\pi\in \mathcal{P}^{\odot}_{m+2 }-\mathcal{P}^{\odot}_{m,2 }$, plus 
  an error term that satisfies (\ref{err0}). 
 Note that by Lemma \ref{lem-ls}
    \bea 
   && \mu\(L(x,r) \prod_{i=1}^{m}\int_{0}^{\ff}f_{i}(X_{t})\,dt\) \label{2.17m}\\
   &&\qquad =  \sum_{\pi\in \mathcal{P}^{\odot}_{m+1 }}  \int  \prod_{j=1}^{m+1}u (z_{\pi(j)},z_{\pi(j+1)})      \prod_{i=1}^{m}  f_{i}  (z_{i} )\,dz_{i}\,f_{r,x} (z_{m+1 } ) \,dz_{m+1} \nn.
 \eea
    
%
%
  
   Let $\ov\pi\in \mathcal{P}^{\odot}_{m+2 }-\mathcal{P}^{\odot}_{m,2 }$, i.e.   for some $j, j+1$ mod $m+2$, we have  $\ov\pi(j)=m+1,\ov\pi(j+1)= m+ 2 $,    and consider  the term on the right-hand side of  (\ref{20.15s2}) for $\ov\pi$,
    \begin{equation}
     \int  \prod_{j=1}^{m+2}u (z_{\ov\pi(j)},z_{\ov\pi(j+1)})     \prod_{i=1}^{m}  f_{i}  (z_{i} )\,dz_{i}\prod_{i=m+1}^{m+2}f_{r,x} (z_{i } ) \,dz_{i}.\label{2.17a}
   \end{equation}
   Note that there is a sequence  in (\ref{2.17a}) of the form
\be u (z_{a},z_{m+1})u (z_{m+1},z_{m+2 })u (z_{m+2 },z_{b})\ee
where $a,b\neq m+1$ or $m+2$.

 	    Consider  a portion of  (\ref{2.17a}) of the form
\begin{eqnarray}
\lefteqn{ \int u (z_{a},z_{m+1})u (z_{m+1},z_{m+2 })u (z_{m+2 },z_{b}) \prod_{i=m+1}^{m+2}f_{r,x} (z_{i } ) \,dz_{i}
\label{an.1a} }\\
&& =\int u (z_{a},x+rz_{m+1})u (rz_{m+1},rz_{m+2 })u (x+rz_{m+2 },z_{b})\nonumber\\
&& \hspace{2 in}f (z_{m+1} )f (z_{m+2} )\,dz_{m+1}\,dz_{m+2}. \nonumber 
\end{eqnarray}
Note that   
\bea 
 &&  u (z_{a},x+rz_{m+1})u (x+rz_{m+2 },z_{b})\label{2.20qq}\\
 &&\qquad= \nn  u (z_{a},x)u (x,z_{b})+(\De_{rz_{m+1}}  u (z_{a},x))u (x,z_{b})\\
 && \qquad\quad+ u (z_{a},x)(\De_{-rz_{m+2}}  u (x,z_{b}) )+ (\De_{rz_{m+1}}  u (z_{a},x)  )(\De_{-rz_{2}} ) u (x,z_{b}) \nn.
 \eea 
We abbreviate this by
\bea
&&
 u (z_{a},x+rz_{m+1})u (x+rz_{m+2 },z_{b})\label{2.20q}\\
 &&\qquad=  u (z_{a},x)u (x,z_{b})  +C_{\De}(z_{a},z_{b},x,r,z_{m+1},z_{m+2})\nn.
   \eea
In this notation we can write (\ref{an.1a}) as
\bea
\lefteqn{  
		  \mbox{ch}_{1}( r)\,u (z_{a},x)u (x,z_{b})+\int  C_{\De}(z_{a},z_{b},x,r,z_{m+1},z_{m+2})\label{2.22} }\\
 &&\hspace{1.5 in}
u (rz_{m+1},rz_{m+2 })  f (z_{m+1} )f (z_{m+2} )\,dz_{m+1}\,dz_{m+2} . \nn\eea 
Using (\ref{2.20q}) again with $z_{m+1}=z_{m+2}=z$ and the fact that the integral of $f$ is equal to 1, we see that   (\ref{2.22}) is equal to 
\begin{eqnarray}
\lefteqn{ \mbox{ ch}_{1}( r) \( \int  u (z_{a},x+rz)u (x+rz,z_{b})  \,f(z)\,dz\)\label{2.23q}}\\
&& \hspace{1in} -\mbox{ch}_{1}( r)\( \int  C_{\De}(z_{a},z_{b},x,r,z ,z ) \,f(z)\,dz\) \nonumber \\
&& \hspace{-.3 in}+\int  C_{\De}(z_{a},z_{b},x,r,z_{m+1},z_{m+2})
u (rz_{m+1},rz_{m+2 })  f (z_{m+1} )f (z_{m+2} )\,dz_{m+1}\,dz_{m+2}. \nn
\eea

  Using (\ref{2.23q}) and the identity
\bea 
    \int  u (z_{a},x+rz)u (x+rz,z_{b})  \,f(z)\,dz= \int  u (z_{a},z)u ( z,z_{b})  \,f_{r,x}(z)\,dz
 \eea we can write the integral  in   (\ref{2.17a}) as  
 \be
\mbox {ch}_{1}( r)  \int  \prod_{j=1}^{m+1}u (z_{\pi'(j)},z_{\pi'(j+1)})  f_{r,x} (z_{m+1 } ) \,dz_{m+1}     \prod_{i=1}^{m}  f_{i}  (z_{i} )\,dz_{i} +\EE_{{\bf f }} \(H'_{r}(x,{\bf z })\),\label{2.19} 
   \ee
   where $\pi'$ is the permutation  in $ \mathcal{P}^{\odot}_{m+1 }$ obtained from $\bar\pi$ by removing $m+2$ from the sequence $(\bar\pi(1),\ldots, \bar\pi(m+2))$ and $H'_{r}(x,{\bf z })$ contains the error terms which are given in the last two lines of (\ref{2.23q}). (All of them have at least one factor of the form $\De_{\cdot}u(\cd)$.)   Moreover we can repeat the argument in the last three paragraphs when $\ov\pi(j)=m+2  $ and $\ov\pi(j+1)= m+ 1 $, so that there are a total of 2  terms that can rewritten as the integral in (\ref{2.19}). 
Using this and (\ref{2.17m}) we establish (\ref{20.15s}) for $n=2$.

\medskip	 Assume that   (\ref{20.15s}) is proved for $L_{n'}(x,r)$, $n'<n$. For any $  \si=(k_{1}, k_{2},\ldots )$ let $\mathcal{P}^{\odot}_{m+n }(\si)$ denote the set of permutations $\bar\pi\in \mathcal{P}^{\odot}_{m+n }$  
that contain    $k_{i}$
chains of order $i=1,2,\ldots$ in $[m+1,m+n]$.   (A chain of order $i\geq 1$ is  a  sequence  $\bar\pi(j),\pi(j+1),\ldots, \bar\pi(j+i)$ in $[m+1,m+n]$ which is maximal in the sense that $\bar\pi(j-1),$ and $\bar\pi(j+i+1)$ are not in $[m+1,m+n]$. In such a case we refer to $j,j+1,\ldots, j+i$ as chain integers.)  Note that  $\mathcal{P}^{\odot}_{m+n }-\mathcal{P}^{\odot}_{m,n }=\cup_{|\si|\geq 1}\mathcal{P}^{\odot}_{m+n }(\si).$

  In the same way we obtained (\ref{2.23q})   and (\ref{2.19}), we see that  the term for any  $\bar\pi\in \mathcal{P}^{\odot}_{m+n }(\si)$ in the evaluation of 
\be
 \mu\(L^{n}(x,r)   \prod_{i=1}^{m}\int_{0}^{\ff}f_{i}(X_{t})\,dt\),\label{2.20}
 \ee 
  the generalization of (\ref{20.15s2}),
is   the same as  the term in (\ref{20.15s}) for a particular permutation $\pi'\in \mathcal{P}^{\odot}_{m,n-|\si| }$     in the evaluation of 
\be
\mu\( \prod_{i=1 }^{\ff} \(\mbox{ch}_{i}( r)\)^{k_{i}} 
L_{n-|\si|}(x,r) \prod_{i=1}^{m}\int_{0}^{\ff}f_{i}(X_{t})\,dt\),  \label{2.21}
   \ee 
   up to error terms  $H_{r}(x,{\bf z })$.    \label{page11}   To see this we note that (\ref{2.21}) can be written as 
  \be
\prod_{i=1 }^{\ff} \(\mbox{ch}_{i}( r)\)^{k_{i}} \mu\( 
L_{n-|\si|}(x,r) \prod_{i=1}^{m}\int_{0}^{\ff}f_{i}(X_{t})\,dt\) \label{2.21q}.
   \ee 
 We use Lemma \ref{lem-ls} to write out (\ref{2.20}),  as in (\ref{2.17m}), and (\ref{20.15s})  and the induction hypothesis  to write out (\ref{2.21}) so they can be easily compared.
   The permutation $\pi' $ is obtained from  $\bar\pi$ by a method we   call `remove and relabel',  which is used in the much simpler case considered in (\ref{2.19}).   We illustrate this with an example. Consider the case in which  $m=10 $, $n=8$ 
and
\begin{equation}
\bar\pi=(6, 7, 11,13, 8, 9,10,1,14, 12,16, 2, 3, 4, 15, 5, 17, 18).\label{am1}
\end{equation}
There are three chains   in this sequence:  
\be (11,13)\qquad(14,12,16)\qquad(17, 18)\label{2.30q}
\ee so that  $\si=  (2,1,0,\ldots)$ and $\bar \pi\in \mathcal{P}^{\odot}_{18 }(2,1,0,\ldots)$.
We first remove all but the first element in each chain   in (\ref{am1}) to obtain
\begin{equation}
 (6, 7, 11, 8, 9,10,1,14,  2, 3, 4, 15, 5,17).\label{am2}
\end{equation}
  The permutation  $\pi'\in \mathcal{P}^{\odot}_{10,4 }$ is obtained from (\ref{am2}) by relabeling the remaining elements in $(11,\ldots,18)$ in increasing order from left to right, i.e.,
\begin{equation}
 \pi'=(6, 7, 11, 8, 9,10,1,12,  2, 3, 4, 13, 5, 14).\label{am2q}
\end{equation}
 Let   $ | \mathcal{P}^{\odot}_{18 }(2,1,0,\ldots)_{\pi'}| $ denote the number of permutations in $\mathcal{P}^{\odot}_{18 }(2,1,0,\ldots)$ that  give rise to $\pi'$. We compute $| \mathcal{P}^{\odot}_{18 }(2,1,0,\ldots)_{\pi'}|$. 
 Clearly, any of the $8!$ permutations of the elements $(11,12,\ldots,18)$ in $\bar\pi$ give rise to distinct permutations in $ \mathcal{P}^{\odot}_{18 }(2,1,0,\ldots)_{\pi'}$. We call these  internal permutations. Furthermore, we consider  the single integer $15$ in $\bar \pi$ to be a chain of order zero. Adding this chain to the three in (\ref{2.30q}) allows us to consider that $\bar\pi $ contains four chains. Clearly, each of the $4!$ arrangements  of these four chains  correspond to distinct permutations in $ \mathcal{P}^{\odot}_{18 }(2,1,0,\ldots)_{\pi'}$. However, we do not want to count the  interchanges of the two chains of order one, since they are counted in the  internal permutations. Consequently  
 \be
 | \mathcal{P}^{\odot}_{18 }(2,1,0,\ldots)|_{\pi'}={8!4! \over 2!} .
 \ee

\medskip	
For general $\si$ and $\pi'\in \mathcal{P}^{\odot}_{m+n-|\si| }$, in which,  as in the  example above, the integers $m+1,\ldots, m+n-|\si| $ appear in increasing order,
\begin{eqnarray}
&&| \mathcal{P}^{\odot}_{m+n }(\si)_{\pi'}| ={n! \over \prod_{i=1}^{\ff}k_{i}! }  {(n-|\si|)!\over (n-|\si|_{+})!}.\label{11.20q}
\end{eqnarray} 
 
 To see this first note that there are $n!$ internal permutations. Since for any $\bar\pi\in \mathcal{P}^{\odot}_{m+n }(\si)_{\pi'}$  there are $|\si|_{+}$  integers from $\{m+1,\ldots, m+n\}$ in the chains  of order $1,2,\ldots$, there are also $n-|\si|_{+}$ remaining integers in  $\{m+1,\ldots, m+n\}$ which, as above,  we consider to be   chains of order $0$. The total number of these chains, including those of order $0$, is   
$n-|\si|_{+}+\sum_{i=1}^{\ff}k_{i}=n-|\si|$. Thus any of the $(n-|\si|)!$ permutations  of these chains 
in  $\bar\pi$ are in $\mathcal{P}^{\odot}_{m+n }(\si)_{\pi'}$. However, we do not want to  count the $(n-|\si|_{+})!\prod_{i=1}^{\ff}k_{i}! $ interchanges of chains of order   0 and $k_{i}$ among themselves, since this has already been counted in the  internal permutations.
  Putting all this  together gives 
(\ref{11.20q}). 

  Consider  (\ref{2.21}) again and the particular permutation $\pi'\in \mathcal{P}^{\odot}_{m,n-|\si| }$. We have already pointed out that there are $(n-|\si|)!$ different permutations,  the internal permutations,  in $  \mathcal{P}^{\odot}_{m,n-|\si| }$,  for which 
\begin{equation}
    \mu\( \prod_{i=1 }^{\ff} \(\mbox{ch}_{i}( r)\)^{k_{i}} 
L_{n-|\si|}(x,r) \prod_{i=1}^{m}\int_{0}^{\ff}f_{i}(X_{t})\,dt\)
   \end{equation}
is the same as it is for  $\pi'$. Therefore  up to the error terms, the contribution to (\ref{2.20}) from $\mathcal{P}^{\odot}_{m+n }(\si)$ is equal to 
\bea
&&
{n! \over \prod_{i=1}^{\ff}k_{i}! (n-|\si|_{+})!} \mu\( \prod_{i=1 }^{\ff} \(\mbox{ch}_{i}( r)\)^{k_{i}} 
L_{n-|\si|}(x,r) \prod_{i=1}^{m}\int_{0}^{\ff}f_{i}(X_{t})\,dt\)  \nn\\
&&\qquad=\mu\( J_{n}(\si,r)\prod_{i=1}^{m}\int_{0}^{\ff}f_{i}(X_{t})\,dt\).\label{2.21s}
   \eea 
  Considering (\ref{rilt.1}) and the fact that 
$\mathcal{P}^{\odot}_{m+n }-\mathcal{P}^{\odot}_{m,n }=\cup_{|\si|\geq 1}\mathcal{P}^{\odot}_{m+n }(\si)$, we see that   the induction step in the proof of (\ref{20.15s}) is proved.

\medskip 	We iterate the steps  used in the proof of   (\ref{20.15s}), and use the fact that each of the $L_{n_{i}}(x_{i}, r  )$ are sums of multiples of $L(x_{i}, r  )=\int_{0}^{\ff}f_{r,x_{i}}(X_{t})\,dt$ to obtain   
 \bea
&&
\mu\(  \prod_{i=1}^{k} L_{n_{i}}(x_{i}, r_{i}  ) \) \label{20.15aa}\\
&&\qquad= \sum_{\pi\in  \mathcal{\bar P}^{\odot}_{n_{1},\ldots,n_{k}}}   \int  \prod_{j=1}^{n}u (z_{\pi(j)} ,z_{\pi(j+1)} )\prod_{j=1}^{n} f _{r_{g(j)},x_{g(j)}}  (z_{j} )\,dz_{j}\nonumber\\
&&\hspace{2.5 in} +\EE_{{\bf f }}(E_{r_{1},\ldots, r_{k}}(x_{1},\ldots, x_{k},{\bf z })), \nn
\eea
where $ \pi(n+1) = \pi(1) $ and  $  \mathcal{\bar P}^{\odot}_{n_{1},\ldots,n_{k}}$ is the set of permutations $\pi$ of $[1,n]$ on the circle, with the property   that for all $j$, when $\pi(j)\in \big[1+\sum_{p=1}^{i-1}n_{p},\sum_{p=1}^{i }n_{p}\big]:=B_{i}$ then $\pi(j+1)\notin B_{i}$,  for all $i\in [1,k]$,  and     $g(j)=i$ when  $j\in B_{i}$. (In the last term in (\ref{20.15aa}) we use the notation introduced in (\ref{2.13}).)

 Note that  
\bea 
\lefteqn{ \int  \prod_{j=1}^{n}u (z_{\pi(j)} ,z_{\pi(j+1)} )\prod_{j=1}^{n} f _{r_{g(j)},x_{g(j)}}  (z_{j} )\,dz_{j}=\label{2.23} }\\
 && \nn  \int  \prod_{j=1}^{n}u (x_{g(\pi(j))}+r_{g(\pi(j))}z_{\pi(j)} ,x_{g(\pi(j+1))}+r_{g(\pi(j+1))}z_{\pi(j+1)})\prod_{j=1}^{n} f    (z_{j} )\,dz_{j}.
 \eea
For each $j=1,\ldots,n$ we write  
\bea
&&u (x_{g(\pi(j))}+r_{g(\pi(j))}z_{\pi(j)} ,x_{g(\pi(j+1))}+r_{g(\pi(j+1))}z_{\pi(j+1)})\label{2.8qa}  \\
 &&\quad\qquad=u (x_{g(\pi(j))},x_{g(\pi(j+1))})+\De_{h}u (x_{g(\pi(j))},x_{g(\pi(j+1))})\nn
   \eea
   where $h=r_{g(\pi(j+1))}z_{\pi(j+1)}-r_{g(\pi(j))}z_{\pi(j)}$.
     Substituting this into (\ref{2.23}) and putting all the terms with one or more  $\De_{\cdot}u$ into  $E'_{r_{1},\ldots, r_{k}}(x_{1},\ldots, x_{k},{\bf z })$ we see that (\ref{2.23}) is equal to   
   \begin{equation}
   \int  \prod_{j=1}^{n}u (x_{g(\pi(j))},x_{g(\pi(j+1))})\prod_{j=1}^{n} f    (z_{j} )\,dz_{j}+ \EE_{{\bf f }}\(E'_{r_{1},\ldots, r_{k}}(x_{1},\ldots, x_{k},{\bf z })\).\label{2.35}
   \end{equation}
We use (\ref{2.23}) and (\ref{2.35}) in (\ref{20.15aa}) and sum over $\pi\in \mathcal{M}_{a}$, rather than $\pi\in  \mathcal{\bar P}^{\odot}_{n_{1},\ldots,n_{k}}$, to get  
\bea
&&
\mu\(  \prod_{i=1}^{k} L_{n_{i}}(x_{i}, r_{i}  ) \) \label{2.36}\\
&&\qquad= {\prod_{i=1}^{k} n_{i}!\over n}\sum_{\pi\in \mathcal{M}_{a}}     \prod_{j=1}^{n}u (x_{\pi(j)} ,x_{\pi(j+1)} )  +\EE_{{\bf f }}\(E'_{r_{1},\ldots, r_{k}}(x_{1},\ldots, x_{k},{\bf z })\). \nn 
\eea
We use the fact that each term in the sum in (\ref{2.36})  comes  from $\prod_{i=1}^{k} n_{i}! $ different terms in (\ref{20.15aa}).  The factor ${1/n   }$ comes from the fact that $ \mathcal{\bar P}^{\odot}_{n_{1},\ldots,n_{k}}$  contains   permutations of $[1,n]$ on the circle, whereas $\mathcal{M}_{a}$ does not. 

 We integrate both sides of  (\ref{2.36}) with respect to $\prod_{i=1}^{k} \nu _{i}(x_{i})$  to get
 \bea
  && \mu\(  \prod_{i=1}^{k} \int L_{n_{i}}(x_{i}, r_{i}  )\,d\nu_{i}(x_{i}) \)\label{2.27w}\\
  &&\qquad={\prod_{i=1}^{k} n_{i}!\over n}\sum_{\pi\in \mathcal{M}_{a}}    \int \prod_{j=1}^{n}u (x_{\pi(j)} ,x_{\pi(j+1)} )\prod_{i=1}^{k}\,d \nu _{i}(x_{i})\nn\nn\\
  &&\hspace{1.5in}+\int \EE_{{\bf f }}\( E'_{r_{1},\ldots, r_{k}}(x_{1},\ldots, x_{k},{\bf z })\) \prod_{i=1}^{k}\,d \nu _{i}(x_{i})\nn.
   \eea

  We    show  in  Lemma \ref{lem-bddhat} that for fixed $\{n_{i}\}$, and $\{\nu_{i}\}$
 \begin{equation}
\int \prod_{j=1}^{n}u (x_{\pi(j)} ,x_{\pi(j+1)} )\prod_{i=1}^{k}\,d \nu _{i}(x_{i})\le \prod_{i=1}^{k} \|\nu_{i}\|_{2, \tau_{2n_{i}}}.
   \end{equation}
   We show in Section \ref{sec-iltan} that  
 \begin{equation}
\lim_{|r |\to 0}\sup_{|{ z_{i} }|\le M} \int E'_{r_{1},\ldots, r_{k}}(x_{1},\ldots, x_{k},{\bf z })\prod_{i=1}^{k}\,d\nu_{i}(x_{i})=0\label{esr4},
\end{equation}
where $r=(r_{1},\ldots, r_{k})$.  
 
\medskip	We can now prove Theorem \ref{theo-multiriltintro-m}.  
  By  (\ref{2.27w})-(\ref{esr4}), for $\nu\in\BB_{2n}$,  
\begin{equation}
\mu\(\(  \int   L_{n}(x,r) \,d\nu(x)  \)^{p}\)<\ff\label{2.30}
\end{equation}
for all $0<r\le r_{0}$ and all even  integers $p\ge 2$.
By (\ref{2.27w})
 \bea
  &&\mu\(\( L_{n,r}(\nu) \)^{j}\(  L_{n,r'}(\nu) \)^{p-j}\)\label{2.27wq}\\
  &&\qquad={\prod_{i=1}^{k} (n !)^{p}\over np}\sum_{\pi\in \mathcal{M}_{a}}    \int \prod_{j=1}^{np}u (x_{\pi(j)} ,x_{\pi(j+1)} )\prod_{i=1}^{p}\,d \nu _{i}(x_{i})\nn\nn\\
  &&\hspace{2in}+\int \EE_{\bf f}\(E'_{r,r' ,  j}(x_{1},\ldots, x_{p},{\bf z }) \)\prod_{i=1}^{p}\,d \nu _{i}(x_{i})\nn\\
    &&\qquad:= A+\II_{r,r' ,  j}(\bf z)\nn.  
   \eea
Therefore  
\begin{eqnarray}
 \mu\(\( L_{n,r}(\nu) - L_{n,r'}(\nu) \)^{p}\)
\label{cau}  &=&\sum_{j=0}^{p} (-1)^{j} {p \choose j}\(A+\II_{r,r',j}({\bf z})\)\\
 & = &     \sum_{j=0}^{p} (-1)^{j} {p \choose j} \II_{r,r',j}({\bf z})\le    2^{p}\sum_{j=0}^{p}  \II_{r,r',j}({\bf z})\nonumber.
\end{eqnarray}
Consequently, by (\ref{esr4})
\begin{equation}
  \lim_{r,r'\to 0}  \mu\(\( L_{n,r}(\nu) - L_{n,r'}(\nu) \)^{p}\)=0.
   \end{equation}
This gives  (\ref{rilt.9intro}) for all even integers $p\ge 2$. Furthermore, we can interpolate to see that it holds for all $p\ge 2$.   This     completes  the proof of Theorem \ref{theo-multiriltintro-m}.  

  Now that we have Theorem \ref{theo-multiriltintro-m} we can return to (\ref{2.27w}) and take the limit as the $r_{i}\to 0$ to complete the proof of Theorem \ref{theo-1.3}.\qed

\begin{remark}\label{rem-2.1} {\rm   Theorem \ref{theo-1.3} does not give the value of $\mu(L_{n}(\nu ))$ for any $n$. When $n=1$ it follows from (\ref{2.4}) that $   \mu\(L _{1}(x,r)\)=\ff$.  In general,  for $n\ge 2$, $\mu\(L _{n}(x,r)\)=\pm \ff$.  See, for example, (\ref{rilt.2ex}) and (\ref{rilt.3ex}). }\end{remark}

\subsection{Bounds for  the error terms}\label{sec-iltan} 

  As in (\ref{l2.1}) and (\ref{12.30m})   we define  
\bea
\|\nu\|_{2, \tau_{k}}& =&\(\int |\hat \nu (\la_{1}+\cdots+\la_{k})|^{2}\,\,\prod_{j=1}^{k}|\hat u  (\la_{j})|\,d\la_{j}\)^{1/2}\label{12.21}\\
&\,=&\(\int |\hat \nu (\la)|^{2}\,\,\tau_{k}(\la)\,d\la\)^{1/2}\nn,
\eea
  where $u$ is the potential density of a L\'evy process in $R^{d}$ and $\nu$ is a finite   measure on $R^{d}$.

\begin{lemma}\label{lem-bddhat}

For any $n_{i}, i=1,\ldots, k$   
\begin{equation}
\int\prod_{j=1}^{n} u (x_{\pi(j+1)}-x_{\pi(j)}) \prod_{i=1}^{k}\,d\nu_{i}(x_{i})\leq \frac{1}{(2\pi)^{nd} }\prod_{i=1}^{k} \|\nu_{i}\|_{2, \tau_{2n_{i}}},\label{l2.22}
\end{equation}
  where  $\pi(n+1)=\pi(1)$, $n=\sum_{i=1}^{k}n_{i}$ and $\pi: [1,n]\mapsto [1,k]$ is such that  $|\pi^{-1}(i)|=n_{i}$ for each $i=1,\ldots, k$, and has the property that when $\pi(j)=i$,  $\pi(j+1)\neq i$. 

In particular for any $\pi\in\PP_{n}$ and smooth function $f$ with compact support
\begin{equation}
\int\prod_{j=1}^{n} u (x_{\pi(j+1)}-x_{\pi(j)}) \prod_{i=1}^{n}\,f(x_{i})\,dx_{i}\leq\frac{1}{(2\pi)^{nd} }\(\int |\hat f(\la)|^{2}\,\,\tau_{2}(\la)\,d\la\)^{n/2}<\ff.\label{l2.22a}
\end{equation}
\end{lemma}

\Proof  
  Since the  integrand is positive, and the  $\nu_{i}$ are finite measures, we can use Fubini's theorem to see that
\begin{eqnarray}
&&\int\prod_{j=1}^{n} u (x_{\pi(j+1)}-x_{\pi(j)}) \prod_{i=1}^{k}\,d\nu_{i}(x_{i})
\label{sl.13}\\
&& =\int_{R_{+}^{n}}  \(\int\prod_{j=1}^{n} p_{t_{j}} (x_{\pi(j+1)}-x_{\pi(j)}) \prod_{i=1}^{k}\,d\nu_{i}(x_{i})\) \prod_{j=1}^{n}e^{-\ga t_{j}}\,dt_{j} \nonumber.
\end{eqnarray}

Considering (\ref{sl.10}) we can use   Fubini's theorem  again to see that the inner integral immediately above, is equal to  $(2\pi)^{nd}$ times
\begin{eqnarray}
&& \int   \(\prod_{j=1}^{n} \int e^{-i(x_{\pi(j+1)}-x_{\pi(j)})\cdot \la_{j}}e^{-t_{j} \bar \ka (\la_{j})}\,d\la_{j}\) \prod_{i=1}^{k}\,d\nu_{i}(x_{i}) \label{s3.6}\\
&&\qquad=\int    \( \prod_{i=1}^{k}\int e^{ i(\sum_{J_{i}}\pm \la_{j } ) \cd x_{i}}
\,d\nu_{i}(x_{i})\)   \prod_{j=1}^{n}e^{-t_{j} \bar \ka (\la_{j})}\,d\la_{j} \nonumber\\
&&\qquad=\int   \( \prod_{i=1}^{k}\,\hat \nu_{i}\(\mbox{$\sum_{J_{i}}$}\pm\la_{j}\)\)   \prod_{j=1}^{n}e^{-t_{j} \bar \ka (\la_{j})}\,d\la_{j} \nonumber.
\end{eqnarray}
  In this formulation  $J_{i}=\pi^{-1}(i)\cup \{\pi^{-1}(i)+1\}$,   and   the sum is taken over all $\la_{j}\in J_{i}$, half of which are multiplied by $-1$.  (The cardinality  $|J_{i}|=2n_{i}$, $i=1,\ldots,k$ and each $\la_{j}$, $j=1,\ldots,n$  appears twice, once in each of  two distinct $J_{i}$, and is multiplied by $-1$ in one of its appearances.   It is not necessary to be more specific.)

Using (\ref{s3.6}) in  (\ref{sl.13}) we see that   (\ref{sl.13}) is bounded by 
\begin{eqnarray}
&&\int_{R_{+}^{n}} \int   \( \prod_{i=1}^{k}\,|\hat \nu_{i}\(\mbox{$\sum_{J_{i}}$}\pm\la_{j}\)|\)   \prod_{j=1}^{n}|e^{-t_{j} \bar \ka (\la_{j})}|\,d\la_{j}\prod_{j=1}^{n}e^{-\ga t_{j}}\,dt_{j}
\label{sl.15}\\
&&\qquad=\int_{R_{+}^{n}} \int   \( \prod_{i=1}^{k}\,|\hat \nu_{i}\(\mbox{$\sum_{J_{i}}$}\pm\la_{j}\)|\)   \prod_{j=1}^{n}e^{-t_{j} \mbox{ Re }\bar \ka (\la_{j})}\,d\la_{j}\prod_{j=1}^{n}e^{-\ga t_{j}}\,dt_{j}   \nonumber\\
&&\qquad=\int   \( \prod_{i=1}^{k}\,|\hat \nu_{i}\(\mbox{$\sum_{J_{i}}$}\pm\la_{j}\)|\)   \prod_{j=1}^{n} {1\over  \ga +\mbox{ Re }\bar \ka (\la_{j})}\,d\la_{j}\nonumber\\
&&\qquad\leq C'^{n}\int    \prod_{i=1}^{k}\,|\hat \nu_{i}\(\mbox{$\sum_{J_{i}} \pm\la_{j}$}\)  \!|  \prod_{j=1}^{n}|\hat u (\la_{j})|\,d\la_{j}   \nonumber.
\end{eqnarray}
Here the second equality uses Fubini's theorem since the  integrand is positive, and the last inequality follows from 
(\ref{s77}). 

%

%

%

Repeated applications  of the Cauchy-Schwarz inequality to the final line of (\ref{sl.15})   and the eventual recognition that we can change $\{\pm\la_{j}\}$ to $\{ \la_{j}\}$ gives (\ref{l2.22}).  We give some idea of how this goes.
Assume for definiteness that $\la_{1}$ appears in $J_{l}$ and $J_{m}$. Then by the Cauchy-Schwarz inequality
\begin{eqnarray}
&&\int    \prod_{i=1}^{k}\,|\hat \nu_{i}\(\mbox{$\sum_{J_{i}} \pm\la_{j}$}\)  \!|  \prod_{j=1}^{n}|\hat u (\la_{j})|\,d\la_{j}
\label{l2.24}\\
&&\qquad \leq \int  \prod_{i\neq l,m} \,|\hat \nu_{i}\(\mbox{$\sum_{J_{i}} \pm\la_{j}$}\)  \!|  \(\int |\hat \nu_{l}\(\mbox{$\sum_{J_{l}} \pm\la_{j}$}\)  \!|^{2} |\hat u (\la_{1})|\,d\la_{1} \)^{1/2} \nonumber\\
&&  \hspace{1 in} \(\int |\hat \nu_{m}\(\mbox{$\sum_{J_{m}}\pm \la_{j}$}\)  \!|^{2} |\hat u (\la_{1})|\,d\la_{1} \)^{1/2}\,\,   \prod_{j=2}^{n}|\hat u (\la_{j})|\,d\la_{j}. \nonumber
\end{eqnarray}

Next assume now that $\la_{2}$ appears in $J_{l}$ and $J_{m'}$.
Using  the Cauchy-Schwarz inequality again we see that (\ref{l2.24})
\begin{eqnarray}
\lefteqn{\leq \int  \prod_{i\neq l,m, m'} \,|\hat \nu_{i}\(\mbox{$\sum_{J_{i}} \pm\la_{j}$}\)  \!| \(\int |\hat \nu_{l}\(\mbox{$\sum_{J_{l}}\pm \la_{j}$}\)  \!|^{2} |\hat u (\la_{1})|\hat u (\la_{2})|\,d\la_{1} \,d\la_{2} \)^{1/2} 
\nn}\\
&&\hspace{-.2in}   \(\int |\hat \nu_{m}\(\mbox{$\sum_{J_{m}}\pm \la_{j}$}\)  \!| ^{2}   |\hat u (\la_{1})|\,d\la_{1} \)^{1/2}\(\int |\hat \nu_{m'}\(\mbox{$\sum_{J_{m'}}\pm \la_{j}$}\)  \!|^{2} |\hat u (\la_{2})|\,d\la_{2} \)^{1/2}\nn\\
&&\hspace{3in}  \prod_{j=3}^{n}|\hat u (\la_{j})|\,d\la_{j}.\label{l2.25} 
\end{eqnarray}
Let us continue and concentrate on the the integrals of  $ |\hat \nu_{l} (\mbox{$\sum_{J_{l}} \pm\la_{j}$} )   |^{2}$. Since  $|J_{l}|=2n_{i}$, the procedure above will ultimately result in the term    
\begin{equation}
  \(\int |\hat \nu_{l}\(\mbox{$\sum_{J_{l}}\pm  \la_{j}$}\)  \!|^{2} \prod_{i=1}^{2n_{i}}|\hat u (\wt  \la_{i})|  \,d \wt  \la_{i}  \)^{1/2},\label{3.10}   \end{equation}
where $\wt  \la_{i}$ is an ordering of the $\la_{j}\in J_{l}$.     Now  we  can use  the fact that $|\hat u (-\la)|=|\hat u (\la)|$ to replace $ |\hat \nu_{l} (   \sum_{J_{l}} \pm\la_{j}   )  |$  by $ |\hat \nu_{l} (\sum_{J_{l}}  \la_{j}  )  |$. With this change 
 (\ref{3.10}) is equal to   $\|\nu_{l}\|_{2,2n_{i}}$, in which the norm is written as in the first equation in (\ref{12.21}). The other terms in (\ref{l2.22}) follow similarly.

  That (\ref{l2.22a}) is finite follows from (\ref{7.8q}), Lemma \ref{lem-9.2}   and  the fact that  $\hat f$ is bounded and  rapidly decreasing.\qed

The  following  is an immediate corollary of Lemma \ref{lem-bddhat}.  

\begin{corollary} \label{cor-3.1} 
\begin{equation}
 {\prod_{i=1}^{k} n_{i}!  \over n} \sum_{\pi\in \mathcal{M}_{a}}   \int   \,   \prod_{j=1}^{n}u (x_{\pi(j)} ,x_{\pi(j+1)} )\prod_{i=1}^{k}\,d\nu_{i}(x_{i})\le \frac{| \mathcal{M}_{a}|}{n} \prod_{i=1}^{k}n_{i}! \|\nu_{i}\|_{2, \tau_{2n_{i}}}.
   \end{equation}
 \end{corollary}
 
\medskip	
 
 We now deal with the error terms. It should be clear that these come about in many different ways. We begin by considering a relatively simple way that they occur. 
 Consider (\ref{20.15aa}). The error terms represented by $E_{r_{1},\ldots, r_{k}}(x_{1},\ldots, x_{k},{\bf z })$ contain chains or variables, which when the analysis is complete, give rise to chains. More error terms are introduced when we show that  
 \bea
&& \sum_{\pi\in  \mathcal{\bar P}^{\odot}_{n_{1},\ldots,n_{k}}}   \int  \prod_{j=1}^{n}u (z_{\pi(j)} ,z_{\pi(j+1)} )\prod_{j=1}^{n} f _{r_{g(j)},x_{g(j)}}  (z_{j} )\,dz_{j}\\
&&\qquad={\prod_{i=1}^{k} n_{i}!\over n}\sum_{\pi\in \mathcal{M}_{a}}     \prod_{j=1}^{n}u (x_{\pi(j)} ,x_{\pi(j+1)} ) +\EE_{\bf f}\(F'_{r_{1},\ldots, r_{k}}(x_{1},\ldots, x_{k},{\bf z } )\)\nn
   \eea
 in (\ref{20.15aa})--(\ref{2.36}). (Note that $F'$ is not the same as $E'$ in (\ref{2.36})  because  $E'$ also contains the error terms in $E$ in  (\ref{20.15aa}).)
 
 \begin{lemma}  \begin{equation}
\lim_{|r| \to 0}\sup_{\forall |z_{i}|\le M} \int F'_{r_{1},\ldots, r_{k}}(x_{1},\ldots, x_{k},{\bf z })\prod_{i=1}^{k}\,d\nu_{i}(x_{i})=0\label{esr4qq}.
\end{equation}
 \end{lemma}
 
 \Proof  Consider (\ref{2.23})--(\ref{2.35}) and set 
 \begin{equation}
   h(g,\pi,j)=r_{g(\pi(j+1))}z_{\pi(j+1)}-r_{g(\pi(j))}z_{\pi(j)}.
   \end{equation} 
  The expression $  F'_{r_{1},\ldots, r_{k}}(x_{1},\ldots, x_{k},{\bf z })$  consists of  $2^{n}-1$ error terms created in the transition from  (\ref{2.23})--(\ref{2.35}). 
    Each of them is of  the form
\bea 
 \int  \prod_{j=1}^{\ell-1}u (x_{\pi(j)} ,x_{\pi(j+1)} )\prod_{j=\ell }^{n}\De_{h(g,\pi, j)}u (x_{\pi(j)} ,x_{\pi(j+1)} )\prod_{j=1}^{n} f   (z_{j} )\,dz_{j},\label{2.56}
\eea 
for some $\ell\le n$.

 Let
\begin{equation}
   V (x_{1},\ldots,x_{k};\pi):= \prod_{j=1}^{\ell-1}u (x_{\pi(j)} ,x_{\pi(j+1)} )\prod_{j=\ell }^{n}\De_{h(g,\pi, j)}u (x_{\pi(j)} ,x_{\pi(j+1)} )\label{2.57}.
   \end{equation}
   We go through  the steps in   the proof of  Lemma \ref{lem-bddhat} to see that  
  \be 
 \Big|\!\int V(x_{1},\ldots,x_{k};\pi)\prod_{i=1}^{k}  d\nu_{i}(x_{i}) \Big|
\! \le \!\int    \prod_{i=1}^{k}\,\big|\hat \nu_{i}\(\mbox{$\sum_{J_{i}}\pm \la_{j}$}\)\big| \,\,   \prod_{j=1}^{n}|\wh{T_{h(g,\pi, j) }  u} (\la_{j})|\,d\la_{j},\label{3.15}
   \ee 
where $T_{h(g,\pi, j)}$ is either the identity or $\De_{h(g,\pi, j)}$ and at least for one $1\le j\le n$,   $T_{h(g,\pi, j)}$ is of the form $\De_{h(g,\pi, j)}$. 
 Note that in general
\bea   |\wh{\De_{h(g,\pi, j)} u(\la_{j})}|&=&|1-e^{ih(g,\pi,j)\la_{j}}|\,|\hat u(\la_{j})|
  \label{3.17}\\
&\le &\( |1-e^{ir_{g(\pi(j+1))}z_{\pi(j+1)} \la_{j}}| +|1-e^{ir_{g(\pi(j ))}z_{\pi(j )}\la_{j}}| \)|\hat u(\la_{j})|\nn\\
&\le&4 |\hat u(\la_{j})|.\nn
   \eea 
We use the  bound $ |\wh{\De_{h(g,\pi, j)} u(\la_{j})}|\le 4 |\hat u(\la_{j})|$ for all but one of the terms $ |\wh{\De_{h(g,\pi, j)} u(\la_{j})}|$ and use the bound in the second line of (\ref{3.17}) for only one of the terms  $ |\wh{\De_{h(g,\pi, j)} u(\la_{j})}|$.
 
 To simplify the notation let us suppose that for this choice $j=1$ and we have
 \be  
   |\wh{\De_{h(g,\pi, 1)} u(\la_{1})}| \le  \( |1-e^{ir' z' \la_{1}}| +|1-e^{ir''z''\la_{1}}| \)|\hat u(\la_{1})|.
      \ee 
 (The reader will see that it doesn't matter what the specific values of $r',r'', z',z''$ are.)

%
 
  It follows from this that the expression in (\ref{3.15})  is 
     \bea
&&\le  2^{n}\int    \prod_{i=1}^{k}\,\big|\hat \nu_{i}\(\mbox{$\sum_{J_{i}} \pm\la_{j}$}\)\big| \, |1-e^{ir 'z' \la_{1}}|\prod_{j=1}^{n}|\hat   u  (\la_{j})|\,d\la_{j}\label{3.15a}\\
&&\qquad+ \,2^{n}\int    \prod_{i=1}^{k}\,\big|\hat \nu_{i}\(\mbox{$\sum_{J_{i}} \pm\la_{j}$}\)\big| \, |1-e^{ir'' z'' \la_{1}}|\prod_{j=1}^{n}|\hat   u  (\la_{j})|\,d\la_{j}.\nn
   \eea
  We define $J_{\ell}$ and $J_{m}$ as in (\ref{l2.24}) and apply the Cauchy Schwarz inequality as in  (\ref{l2.24}) and get that (\ref{3.15})
   \bea
&&  \leq 2^{n}\int  \prod_{i\neq l,m} \,|\hat \nu_{i}\(\mbox{$\sum_{J_{i}} \pm\la_{j}$}\)  \!|  \(\int |\hat \nu_{l}\(\mbox{$\sum_{J_{l}} \la_{j}$}\)  \!|^{2} |1-e^{ir' z' \la_{1}}|^{2} |\hat u (\la_{1})|\,d\la_{1} \)^{1/2} \nonumber\\
&&  \hspace{.7 in} \(\int |\hat \nu_{m}\(\mbox{$\sum_{J_{m}}\pm \la_{j}$}\)  \!|^{2} |\hat u (\la_{1})|\,d\la_{1} \)^{1/2}\,\,   \prod_{j=2}^{n}|\hat u (\la_{j})|\,d\la_{j}, \label{2.62}
\end{eqnarray}
 plus a similar term but with $r',z'$ replaced by $r'',z''$. Note that in applying the Cauchy Schwarz inequality we take 
\begin{equation}
    |1-e^{ir' z' \la_{1}}|  |\hat u (\la_{1})|=\(( |1-e^{ir' z' \la_{1}}|  |\hat u (\la_{1})|^{1/2}\) |\hat u(\la_{1})|^{1/2},
   \end{equation}
so that we get $ |1-e^{ir' z' \la_{1}}| ^{2}|\hat u(\la_{1})|$ in one of the terms we integrate with respect to $\nu_{1}$.

 We proceed  as in the proof of Lemma \ref{lem-bddhat} to get 
   \begin{equation}
  \Big|  \int V (x_{1},\ldots,x_{k};\pi)\prod_{i=1}^{k}d\nu_{i}(x_{i}) \Big| \le 2^{n} \sup_{|z'|\le 1 }\|\nu_{1}\|_{2, \tau_{2n_{1}},r'z'}\prod_{i=2}^{k} \|\nu_{i}\|_{2, \tau_{2n_{i}}},\label{3.18}
   \end{equation}
    plus a similar term but with $r',z'$ replaced by $r'',z''$. 
  Here
\bea
\|\nu\|_{2, \tau_{k},r'z'}&=&\(\int |\hat \nu (\la_{1}+\cdots+\la_{k})|^{2}\,\,|1-e^{ir'z'  \la_{1}}|^{2}\prod_{j=1}^{k}|\hat u (\la_{j})|\,d\la_{j}\)^{1/2}\label{l2.21}\\
&=&\(\int |\hat \nu (\la)|^{2}\,\,\tau_ {k ,r'z'}(\la)\,d\la\)^{1/2}\nn
\eea
and
\be \tau_{k,r'z'}(\la)= \int |1-e^{ir'z'  \la_{1}}|^{2 }\,  |\hat u  (\la_{1}) |\tau_{ k-1}(\la-\la_{1})\,d\la_{1} .\label{2.65x}
\ee
We also use the fact that $|z'|\le 1$, since $z'$ is in the domain of $f$ which is the unit ball of $R^{d}$.

Since the right-hand side of   (\ref{3.18}) does not depend on $z'$ and the integral over $\prod_{j=1}^{n} f   (z_{j} )\,dz_{j}$ is equal to one, we  
 see that each of the error terms 
 \begin{equation}
\le 4^{n}\sup_{|z'|\le M} \|\nu_{1}\|_{2, \tau_{2n_{1}},r'z'}\prod_{i=2}^{k} \|\nu_{i}\|_{2, \tau_{2n_{i}}}.\label{2.65}
   \end{equation}
 (Here we combine the terms in $r',z'$ and $r'',z''$.)  
   
 Putting this together we see that there exists a $j\in[1,k]$ such that
\be 
   \int  F'_{r_{1},\ldots, r_{k}}(x_{1},\ldots, x_{k},{\bf z })\prod_{j=1}^{k}d\nu(x_{i}) 
 \le 8^{n}\sup_{|z'|\le M} \|\nu_{j}\|_{2, \tau_{2n_{j}},r_{j}z'}\prod_{i=1,i\ne j}^{k} \|\nu_{i}\|_{2, \tau_{2n_{i}}},
   \ee 
   where  $z'$ may be any of the values $z_{1},\ldots,z_{n}$.
  Note that by  (\ref{7.8q}) 
  \begin{equation}
   \tau_{k,r'z'}(\la)\le \vartheta _{k}(r'z',\la),
   \end{equation}
 but with $h$ replaced by $\va_{\al}$.  (See Lemma \ref{lem-hatbound} and the comments in Section \ref{sec-2.1}.) Therefore, by (\ref{B.7k})   and (\ref{8.65}), for any $2\le k\le 2n_{1}$,  
 and $1\leq i\leq k$
  \begin{equation}
   \|\nu_{i}\|_{2, \tau_{k},r'z'}\le  o\( \wt H ( 1/|r'z'|)  \)^{ n_{1}- k/2} ,\quad\mbox{as $r'\to 0$},\label{2.69}
   \end{equation}
 and, in particular,  
  \begin{equation} 
   \sup_{|z'|\le M} \|\nu_{i}\|_{2, \tau_{2n_{1}},r'z'}=o(1),  \quad\mbox{as $r'\to 0$}.
   \end{equation}  
Consequently
\begin{equation}
   \lim_{|r| \to 0}    \int  F'_{r_{1},\ldots, r_{k}}(x_{1},\ldots, x_{k},{\bf z })\prod_{j=1}^{k}d\nu(x_{i}) =0.
   \end{equation}
 \qed

 We now deal with the error terms created when we form the chains. Let us first note that these are similar regardless of the length of the chain. Suppose we have a chain of length $\ell$. Following the analysis in (\ref{an.1a})--(\ref{2.23q}) we see that in place of (\ref{2.22}) we get 
\bea  
 &&\mbox{ch}_{\ell}( r)\,u (z_{a},x)u (x,z_{b})\label{2.70} \\
 &&\qquad+\int  C_{\De}(z_{a},z_{b},x,r,z_{1},  z_{\ell+1})\prod_{i=1}^{\ell}
u (rz_{i},rz_{i+1 })  \prod_{i=1}^{\ell+1}f (z_{i} ) \,dz_{i}  . \nn
\eea 
Here $z_{1}$ and $z_{\ell+1}$ are the variables at the end of the chain and $z_{a}$, $z_{b}$ and $x$ are variables that are not in this  chain.    Moreover, the variable $x$ is in the same interval that contains the variables in the chain, i.e. if the chain consists of a sequence of variables $z_{\pi(j)}$ with $\pi(j)=i$  then $x$ is an $x_{i}$, that is to be integrated, generally along with other terms,   by $d\nu_{i}$. Also note that there may be many   chains, of various length,  that   consist  of   sequences of variables in the $i$-th interval.
Similarly, in place of (\ref{2.23q}) we get 
\begin{eqnarray}
\lefteqn{ \mbox{ ch}_{\ell}( r) \( \int  u (z_{a},x+rz)u (x+rz,z_{b})  \,f(z)\,dz\)\label{2.23qq}}\\
&& \hspace{1in} -\mbox{ch}_{\ell}( r)\( \int  C_{\De}(z_{a},z_{b},x,r,z ,z ) \,f(z)\,dz\) \nonumber \\
&&\qquad+\int  C_{\De}(z_{a},z_{b},x,r,z_{1},  z_{\ell+1})\prod_{i=1}^{\ell}
u (rz_{i},rz_{i+1 })  \prod_{i=1}^{\ell+1}f (z_{i} ) \,dz_{i}  .\nn
\eea
The integrals containing the $C_{\De}$ are in the error terms. 

   Note that we can not extract  ch$_{\ell}( r)$ from the last line in (\ref{2.23qq}) because $z_{1}$ and $z_{\ell+1}$ are in $C_{\De}(z_{a},z_{b},x,r,z_{1},z_{\ell+1})$. They give rise to terms like $\De_{rz_{1}}  u (z_{a},x))$ in (\ref{2.20qq}) that contain variables that are not in  the chains. Therefore, to evaluate the error integrals we put off integrating with respect to any of the $z$ variables and first integrate with respect to the measures $\nu_{i}. $    To this end we go back to the definition of the chain functions in (\ref{1L.1}) and write 
  \begin{equation}
  \mbox{ch}_{\ell}(r') =  \int   u (r'z'_{1},r'z'_{2}) \cdots u(r'z'_{\ell},r'z'_{\ell+1})   \prod_{j=1}^{\ell+1} f(z'_{j}) \,dz'_{j} 
  \end{equation} 
We refer to the terms $\prod_{i=1}^{\ell}
u (r'z'_{i},r'z'_{i+1 })$ as   chain integrands.

   \medskip	 
We  arrange the order of integration in the error terms so that we integrate with respect to the $\{z_{i}\}_{i=1}^{n}$ last. Doing this we can write a  typical error term in the form  
\be
\int \(\int\prod_{j=1}^{n'}\wt u_{j}({\bf x,r',z'}) \prod_{i=1}^{k}\,d\nu_{i}(x_{i})\) V({\bf r,z})\prod_{j=1}^{n}  f(z_{j} )\,dz_{j}, \label{l2.5e}
\ee 
where $n'=\sum_{i=1}^{k}n'_{i}$,   $n'_{i}=n_{i}-|\si(i)| $,  $|\si(i)|=\sum_{j=1}^{\ff}  jk_{j}(i)$.
 The term  $V({\bf r,z})$ is the product of all the  chain integrands.
(Let  $\AA$ denote  the set of all the   variables $z_{i}$ in  the chain integrands. The term ${\bf z}$ in $V({\bf r,z})$ refers to these variable.  There are $\sum_{i=1}^{k} \ \sum_{j=1}^{\ff}  (j+1)k_{j}(i)$ of them.
 The term ${\bf r}$  refers to whatever values of $r_{1},\ldots,r_{k} $ are in  $V({\bf r,z})$.) 

 The functions $\wt u_{j}({\bf x,r',z'})$
  include all the terms not included in the chain integrands. In particular they include all terms of the form $C_{\De}(z_{a},z_{b},x,r,z_{1},z_{k})$
 and $C_{\De}(z_{a},z_{b},x,r,z ,z )$, (see (\ref{2.23q}) and its generalization in (\ref{2.23qq})). 
They  also include many terms that are not of the form  $C_{\De}(\cd))$. These are terms of the form
\begin{equation}
   u (x_{\pi(j+1)}-x_{\pi(j)})+r_{\ell} (  z_{j+1}-  z_{j}))\label{2.71}.
   \end{equation}
They come from the change of variables that has already taken place when we write  (\ref{l2.5e}) with $f$ rather than $f_{r_{\ell},x_{\cd}}$. (Here $r_{\ell}$ refers to one one of the  $r_{1},\ldots,r_{k} $.)   The product $\prod_{j=1}^{n'}\wt u_{j}({\bf x,r',z'})  $ contains all the variables $z_{1},\ldots z_{n'}$ and also some of the variables $\{z_{j}, j\in \AA\}$.
  
   We bound  (\ref{l2.5e})  by  
  \bea 
&&\(\sup_{\forall |z_{i} |\le M}\int\prod_{j=1}^{n'}\wt u_{j}({\bf x,r',z'}) \prod_{i=1}^{n}\,d\nu_{i}(x_{i})\) \int V(r,z)\prod_{j=1} ^{n} f(z_{j} )\,dz_{j} \label{l2.5f}
 \\
  &&\qquad=
\( \sup_{\forall |z_{j}|\leq M}   \int\prod_{j=1}^{n'}\wt u_{j}({\bf x,r',z'})  \prod_{i=1}^{k}\,d\nu_{i}(x_{i})\)\prod_{i=1}^{k}\prod_{j=1 }^{\ff} \(\mbox{ch}_{j}( r)\)^{k_{j}(i)}    .\nn
\eea
 (When we perform the integration with respect to $\prod_{j=1} ^{n} f(z_{j} )\,dz_{j}$ for the variables $z_{i}\notin\AA$ we just get 1.)
  We replace the terms in $\wt u$  of the form given in (\ref{2.71}) by $u(x_{\pi(j)},x_{\pi(j+1)})+\De_{r_{\ell} (  z_{j+1}-  z_{j})}u(x_{\pi(j)},x_{\pi(j+1)})$ and write 
\begin{equation}
   \prod_{j=1}^{n'}\wt u_{j}({\bf x,r',z'}) =\sum _{q}  \prod_{j=1}^{n'}\wt u_{j,q}({\bf x,r',z'}).
   \end{equation}
 Considering the terms relating to the $C_{\De}(\cd\cd\cd)$, we see that $\wt u_{j,q}(\cd)$ has one  of the following forms:
\begin{equation}
   u (x_{\pi(j )},x_{\pi(j+1)}), \quad  \De_{r'z' -r''z'' } u (x_{\pi(j )},x_{\pi(j+1)}),\quad \De_{\pm r'z'}u(x_{\pi(\ell)},x_{\pi(m)}),\label{2.73}
   \end{equation}
where $r',r''$ take values in  $r_{1},\ldots,r_{k} $, and $z',z''$take values in  $z_{1},\ldots z_{n }$, $x_{\pi(j )}\ne x_{\pi(j+1)}$ and $x_{\pi(\ell)}\ne x_{\pi(m)}$. 

Consider a term of the form 
\begin{equation}
 \int\prod_{j=1}^{n'}\wt u_{j,q}({\bf x,r',z'})\prod_{i=1}^{k}\,d\nu_{i}(x_{i})\label{2.78}.
   \end{equation}
Following the   proof of Lemma \ref{lem-bddhat} we can bound this by  
   \bea
  C\int    \prod_{i=1}^{k}\,\big|\hat \nu_{i}\(\mbox{$\sum_{J_{i}} \la_{j}$}\)\big| \,\,   \prod_{j=1}^{n'}|\widehat{\wt u_{j,q} ({\bf \la,r',z'})}|\,d\la_{j},\label{2.79}
\end{eqnarray}
in which the Fourier transform of $\wt u_{j,q}$ is taken with respect to the $x$ variable. 
Note that   $n_{i}'  $ of the functions $\wt u_{j,q}$ in (\ref{2.78}) are integrated by $\nu_{i}$, $1\le i\le k$. Recall that $n'=\sum_{i=1}^{k}n'_{i} $, where   $0\le n'_{i}\le n_{i}$.
To simplify the notation suppose that for $i=1,\ldots,p$, $1\le p<k$, $n_{i'}<n_{i} $, and for $i=p+1,\ldots,k$, $n_{i'}=n_{i} $. These are the intervals  that contain  chains and do not contain  chains, respectively. We now proceed as in (\ref{3.15a})--(\ref{2.65}) to see that (\ref{2.79}) is bounded by
  \begin{equation}
 \le 4^{n}\prod_{i=1}^{p} \sup_{|z'|\le M }\|\nu_{i}\|_{2, \tau_{2n_{i}'},r'z'}\prod_{i=p+1}^{k} \|\nu_{i}\|_{2, \tau_{2n_{i}}} ,\label{2.81}
   \end{equation}
    in which we set $\|\nu_{i}\|_{2, \tau_{0}}=1$. 
By (\ref{2.69}), for $i=1,\ldots,p$  
  \begin{equation}
   \sup_{|z'|\le M} \|\nu_{i}\|_{2, \tau_{2n_{i}'},r'z'}\le o\( \wt H ( 1/|r' |)  \)^{ n_{i}-  n_{i}'}  ,\quad\mbox{as $r'\to 0$}.\label{2.6q9}
   \end{equation}
By Lemma \ref{lem-chat}   applied to   the chains formed by variables in the $i$-th interval,  
\begin{equation}
   \prod_{j=1 }^{\ff} \(\mbox{ch}_{j}( r)\)^{k_{j}(i)}\le  O\( (\wt H( 1/|r  |))^{-(n_{i}-  n_{i}')}   \).
   \end{equation}
Therefore, the limit of (\ref{l2.5f}), as  $r=r'\to 0$, is zero. \qed

\begin{remark}\label{rem-2.2} {\rm   The process 
  $L_{n}(x,r)=B_{n}\(L(x,r)\)$,  $n\ge 1$,  where the polynomials $B_{n}(u)$ satisfy  
  \begin{equation}
  \sum_{n=0}^{\ff}{\(   \sum_{i=0}^{\ff}\mbox{ch}_{i}( r)\,s^{i} \)^{n} \over n!}s^{n}B_{n}(u)= e^{su},\label{ljo}
  \end{equation}
   with   $\mbox{ch}_{0}( r)=1$.  
     To see that (\ref{ljo}) agrees with (\ref{rilt.1}) and (\ref{rilt.2b}) we expand
   \bea
\lefteqn{ {\(   \sum_{i=0}^{\ff}\mbox{ch}_{i}( r)\,s^{i} \)^{n} \over n!}s^{n}\label{slum1}}\\   &&=\sum_{\sum_{i=0}^{\ff} k_{i}=n}
   {1 \over n!} {n\choose k_{0}\,k_{1}\,\cdots}\prod_{i=0}^{\ff}\(\mbox{ch}_{i}( r)\,s^{i} \)^{k_{i}}\,s^{n}\nn\\   &&=\sum_{ \sum_{i=0}^{\ff}k_{i}=n}
   {1 \over \prod_{i=0}^{\ff}k_{i}!} \prod_{i=0}^{\ff}\(\mbox{ch}_{i}( r) \)^{k_{i}}\,s^{(i+1)k_{i}}\nn\\   &&=\sum_{\sum_{i=0}^{\ff}k_{i}=n} {1 \over k_{0}!\prod_{i=1}^{\ff}k_{i}! }\prod_{i=1}^{\ff} \(\mbox{ch}_{i}( r)\)^{k_{i} }s^{\sum_{i=0}^{\ff}(i+1)k_{i}}.\nn
   \eea
  Fix $N$ and consider all $k_{0}, k_{1}, \ldots$ with 
   \be 
   N=\sum_{i=0}^{\ff}(i+1)k_{i}=k_{0}+| \si|_{+}.
   \ee 
 Consequently we can replace  $k_{0}$ with $N-| \si|_{+}$ in (\ref{slum1}).
 
 In addition when  $\sum_{i=0}^{\ff} k_{i}=n$, we   have 
   $n=k_{0}+| \si|_{+}-| \si|=N-| \si|$. We use this observation and  (\ref{slum1}) to equate coefficients of $s^{N}$  in (\ref{ljo})  to obtain
   \begin{equation}
\sum_{\{ \si\,|\,  | \si|_{+} \leq N\} }  {1 \over \prod_{i=1}^{\ff}k_{i}! (N-| \si|_{+})!}\prod_{i=1}^{\ff} \(\mbox{ch}_{i}( r)\)^{k_{i} } 
    B_{N-|\si|} (u)=   { u^{N}  \over N!}.\label{slum2}
   \end{equation}
   Setting $u=L(x,r)$ and  $L_{n}(x,r)=B_{n}\(L(x,r)\)$ we get   (\ref{rilt.1}) and (\ref{rilt.2b}).  
   }\end{remark}

\section{Continuity of   intersection local time  processes}\label{sec-iltcont}

  We begin by reviewing some  conditions for the continuity  of stochastic processes   with increments in an exponential Orlicz space. For a proof of Theorem \ref{maj} see \cite[Section 3]{MRcont}. Let 
\be
\rho_q(x)=\exp(x^q) -1\label{psi}
\ee
  for $1\le q<\ff $, and for
$0<q<1$, we define
\be
\label{
 psi-1}
\rho_q(x)= \left\{ \begin{array}{l@{\quad\quad}r}
K_q\,xÊ &0\leÊ x<\( \frac{1}{q}\)^{1/q}
\\
\exp(x^q) -1Ê & x\ge\( \frac{1}{q}\)^{1/q}Ê \end{array} \right.
\ee
where
\be K_q={\exp( x_0 ^q) -1 \over x_0}\qquad \mbox{and
$\qquad x_0:=x_0(q)=(1/q)^{1/q}$},
\ee
so that $\rho_q(x)$ is continuous.

 Let    $L^{\rho_{q}}(\Om,\FF,P)$ denote the set ofÊ random 
variables 
$\xi:\Om\to R^{1}$ such that $E\rho_{q}\( |\xi|/c \)<\ff$ for some $c>0$.  
$L^{\rho_{q}}(\Om,\FF,P)$ is a Banach space with norm given by  
\be
\|\xi\|_{\rho_{q}}=\inf\left\{c>0:E\rho_{q}\( |\xi|/c\)\le1\right\}.\label{1.4}
\ee
 
  Let $(T,\bar d)$   be aÊ   metric or pseudo-metric space.  Let $B_{\bar d}(t,u)$  denote the closed ball in $(T,\bar d)$  with radius $u$ and center $t$. For any probability measure  $\si$  on   $(T,\ov d)$ we define
\begin{equation}
  J_{T,\bar  d,\si,n}( a) =\sup_{t\in T}\int_0^a \(\log\frac1{\si (B_{\bar d}(t,u))} \)^{n} \,du.\label{tau}
   \end{equation}
  We   use the   following basic continuity theorem to obtain sufficient conditions for continuity of  permanental fields.   
 
   Ê

 \bt\label{maj}ÊÊ     Let $Y= \{Y(t): t\in T\}$Ê be a 
stochastic process  such that $Y(t,\om):T\times \Om\mapsto [-\ff,\ff] $  is $\mathcal{A}\times \mathcal{F}$ measurable for some $\si$-algebra $\mathcal{A}$ Êon $T$.   Suppose $Y(t) 
\in L^{\rho_{1/n}}(\Om,\FF,P)$,   where $n\ge 2$ is an integer , and  there exists a metric  $\bar d$ on $T$ such that
\begin{equation}
  \|Y(s)-Y(t)\|_{\rho_{1/n}}\le \bar d(s,t).\label{2.2ar}
   \end{equation}   
     (Note that the balls $B_{\bar d}(s,u)$ are $\mathcal{A}$ measurable.)

 Suppose that   $(T,\bar d)$ has finite diameter $D$, and  that there exists
a probability measure $\si$Ê  on $(T,\mathcal{A})$ such that
\be
  J_{T, \bar d,\si,n}(D)<\ff.
\label{2.10}
\ee
Then there exists a version $Y'=\{Y'(t),t\in T\}$ of $Y$ such that
\be
E\sup_{t\in T}Y'(t)\leÊ C\,  J_{T,\bar  d,\si,n}(D),
\label{2.2}
 \ee
 for some $C<\ff$. Furthermore for all $0<\de\le D$, 
 \be
Ê ÊÊ \sup_{\stackrel{s,t\in T}{\bar d(s,t)\le \de}} |Y' ( s,\om)-Y'( t,\om)|
Ê \le 2Z(\om) \,  J_{T, \bar d,\si,n}(\de)\label{2.12},
\ee
almost surely, where   
\begin{equation}
   Z(\om):=\inf\bigg\{\al>0:\int_{T}\rho_{1/n}(\al^{-1}|Y(t,\om)|)\,\si (dt) \leq 1\bigg\} 
   \end{equation}
 andÊ $\|Z\|_{\rho_{1/n}} \le
 K $, where $K$ is a constant.  
 
 In particular, if
 \begin{equation}
   \lim_{\de\to 0}  J_{T, \bar d,\si,n}(\de)=0,\label{2.5}
   \end{equation}
    $Y'$ is uniformly continuous on   $(T, \ov d)$ almost surely.
\et

\medskip	  For any   positive measures  $ \phi,\chi\in  \BB _{2}(R^{d})$ set
\begin{equation}
P^{\phi,\chi}(A)={\mu \(L_{1}(\phi)L_{1}(\chi)1_{A}\) \over \mu \(L_{1}(\phi)L_{1}(\chi)\)}. \label{tranh.6}
 \end{equation}
In the next lemma we show that $L_{n}(\nu)
\in L^{\rho_{1/n}}(\Om,\FF,  P^{\phi,\chi})$  for a probability measure $P^{\phi,\chi}$.   Consequently, we     can use  Theorem \ref{maj} to obtain continuity conditions for $\{L_{n}(\nu),\nu\in\VV'\}$,   where $\VV'\in \mathcal {B}_{2n}(R^{d})$,  $n\ge 2$.

\begin{lemma}\label{lem-hatcont} 
For  $\nu,\mu   \in   \BB_{ 2n}(R^{d})$ 
\begin{equation}
\|L_{n}(\nu  )-L_{n}(\mu)\|_{\rho_{1/n}, P^{\phi,\chi}}\leq C_{\phi,\rho,n}  \|\nu -\mu  \|_{2,\tau_{2n}}, \label{tranh.5a}
\end{equation}
where $\|\cdot \|_{\rho_{1/n}, P^{\phi,\chi}}$ denotes the Orlicz space norm with respect to the probability measure $P^{\phi,\chi}$ and $C_{\phi,\chi,n}$ is a constant depending on $\phi$, $\chi$ and $n$.
\end{lemma}

\Proof  It is obvious from (\ref{3}) that $\al_{n,\ep}(\nu_{s_{1}},t)-\al_{n,\ep}(\nu_{s_{1}},t) =\al_{n,\ep}(\nu_{s_{1}}-\nu_{s_{2}},t)$ and since $L_{1,\ep}(\nu)=\al_{1,\ep}(\nu)$, we see from (\ref{zzx}) and (\ref{rilt.9intro}) that $L_{n}(\nu )-L_{n}(\mu)=L_{n}(\nu-\mu )$. Therefore it suffices to show that $\|L_{n}(\nu  ) \|_{\rho_{1/n}, P^{\phi,\rho}}\leq C  \|\nu   \|_{2,\tau_{2n}}$.

We first note that for even integers  $m\ge 2$  
\be \mu\(    |L_{n}(\nu  ) |^{m} \)\leq  (mn)!\|\nu  \|_{2,\tau_{2n}}^{m}.\label{tranh.1} 
\ee 
To see this we   use (\ref{rilt.15intro})   to get 
 \begin{equation}
|\mu\(     L^{m}_{n} \)|\leq  (n!)^{m}|\mathcal{M}_{a}|\,\,\|\nu \|_{2,\tau_{2n}}^{m}.\label{tranh.2}
\end{equation}
 Since $|\mathcal{M}|={mn \choose n\,n\cdots n}$ we get  (\ref{tranh.1}).  
  
Next we note that since  
 \begin{equation}
   \mu \(|L_{1}(\phi)L_{1}(\chi)|\) \le    \(\mu \(L^{2}_{1}(\phi)L^{2}_{1}(\chi)\)\)^{1/2}\label{3.17aa},
   \end{equation}
it follows from (\ref{rilt.15intro}) that both sides of (\ref{3.17aa}) are finite and consequently that for $1\le p\le 2$,
\begin{equation}
      \mu \(|{L }_{1}(\phi)L _{1}(\chi)|^{p}\)\le  \(\mu \(L^{2}_{1}(\phi)L^{2}_{1}(\chi)\)\)^{1/2}+ \(\mu \(L^{2}_{1}(\phi)L^{2}_{1}(\chi)\)\) .\label{3.18k}
   \end{equation}
   It now follows from (\ref{tranh.6}), (\ref{tranh.1}), (\ref{3.17aa}) and (\ref{3.18k}) that for   all integers $k\ge 1$ and $r=2n/2n-1$,  
   \bea
   E_{P^{\phi,\chi}}\((L_{n}(\nu))^{k/n}\)&\le&\frac{ \mu \(|L_{1}(\phi)L_{1}(\chi)|^{r} \)^{1/r} \(\mu \( L_{n}^{2k} \) \)^{1/2n}}  { \mu \(L_{1}(\phi)L_{1}(\chi)\)}\label{3.18a}\\
   &\le&\nn C_{\phi,\chi}  \((2kn)!\|\nu  \|_{2,\tau_{2n}}^{2k}\)^{1/2n},
   \eea
 where $C_{\phi,\chi }$ is a constant depending on $\phi$, $\chi$.
Therefore for all all integers $k\ge 1$  
\begin{equation}
    E_{P^{\phi,\chi}}\( \(\frac{L_{n}}{(2n)^{n}\|\nu  \|_{2,\tau_{2n}}}\)^{k/n} \)\le C_{\phi,\chi}C  k!,\label{3.19}
   \end{equation}
for some constant  $C $. This gives (\ref{tranh.5a}). \qed

\begin{theorem} Let $\{L_{n}(\nu),\nu\in\VV \}$ be an $n$-fold  intersection local time process, where $\VV \in \mathcal {B}_{2n}(R^{d})$,  $n\ge 2$,  and let    $  \bar d((\nu,\mu))=  \|\nu  -\mu  \|_{2,\tau_{2n}}$.  If (\ref{2.5}) holds  $\{L_{n}(\nu),\nu\in\VV \}$ is continuous  on $(\VV,\bar d)$, $P^{\phi,\chi}$ almost surely.
\end{theorem}

\Proof This is an immediate application of Theorem \ref{maj} and {Lemma \ref{lem-hatcont}. \qed

  When $\VV=\{\nu_{x}, x\in R^{d}\}$, the set of translates of a fixed measure $\nu$,   the  simple concrete condition in (\ref{b101})  implies that (\ref{2.5}) holds.

\medskip	\noindent{\bf  Proof of Theorem \ref{theo-cont14.1b} } We show that
\be
\int_{1}^{\ff}  \(\int_{|\xi|\geq x}\tau_{2n}(\xi)|\hat{\nu}(\xi)|^{2}\,d\xi\) ^{1/2} \frac{
 (\log x)^{n-1}}{x}\,dx<\ff,\label{b101q}
\ee   
is a sufficient condition for $\{L_{n}(\nu_x),\,x\in R^d\}$ to be   continuous  $P^{y}$ almost surely, for all $y\in R^{d}$.

For all $s,h\in R^{d}$ 
\bea \| \nu_{s+h}-\nu_{s}\|_{2,\tau_{2n}}&=&\| \nu_{h}-\nu\|_{2,\tau_{2n}}\label{2.49}\\&=&
 \(\int\!\!\int |1-e^{ix  h}|^{2}||\hat\nu(x)|^{2}  \tau_{2n}(x)\,dx \)^{1/2}\nn\\
 &\le&
 \( 3 \int\!\!\int    ((|x||h|)^{2}\wedge 1)||\hat\nu(x)|^{2} \tau_{2n} (x) \,dx \)^{1/2}\nn.
 \eea 
  Using this bound, the fact that (\ref{b101q}) implies (\ref{2.5}) with  $\bar d=\| \nu_{s+h}-\nu_{s}\|_{2,\tau_{2n}}$ is routine. (See the proof  of  \cite[Theorem 1.6]{MR96}, where this is proved in a slightly different context) . 
Consequently we have that  $\{L_{n}(\nu_{x}),\,x\in R^{d}\}$ is continuous $P^{\phi,\chi}$ almost surely. 
 As explained in  the proof of  \cite[Theorem 5.1]{LMR}, this implies that $\{L_{n}(\nu_{x}),\,x\in R^{d}\}$ is continuous $P^{y}$ almost surely, for all $y\in R^{d}$. 
 \qed

\begin{example}  \label{ex-3.1}
{Using Theorem \ref{theo-cont14.1b} we give some examples of  L\'evy processes and measures $\nu$ for which $\{L_{n}(\nu_x),\,x\in R^d\}$ is continuous almost surely. As usual let $u$ denote the potential density of the L\'evy process and assume that (\ref{7.8q}) holds.  Let 
\begin{equation}
  (\wt H(|\xi|))^{-1}:= \int_{|\eta|\le  |\xi| } (\va _{\al} (|\eta|))^{-1}\,d\eta.\label{8.65q}
   \end{equation}
  It follows from  (\ref{8.5n})  that when $\al<d$
   \begin{equation}
  (\wt H(|\xi|))^{-1}\le C   |\xi|^{d} (\va_{\al}(|\xi|))^{-1}  \label{3.22}
   \end{equation}
and, for the  functions $\va_{\al}$ that we consider, (\ref{1.5}) holds.

By Theorem \ref{theo-cont14.1b} and (\ref {9.12w}), $\{L_{n}(\nu_x),\,x\in R^d\}$ is continuous almost surely when 
\be
|\hat\nu(\xi)|^{2}=O\(\frac{\va_{\al}  (|\xi|) \wt H ^{2n-1}(|\xi|)  }{|\xi|^{ d  }(\log |\xi|)^{2n +1 +\de}}\)\qquad 
\mbox{as } |\xi|\to \ff,\label{3.23}
\ee
 for any $\de>0$.

When $\va_{\al}  (|\xi|) =|\xi|^{\al}$, for $d(1-\frac{1}{2n})<\al<d$, using (\ref{3.22}) the right hand-side of (\ref{3.23}) is 
\be
O\(\frac{1}{|\xi|^{2n(d-\al) }  (\log |\xi| )^{2n +1 +\de}} \).
\ee
When   $1/\va_{\al}  (|\xi|) =O(|\xi|^{-d })$ as $|\xi|\to\ff$, the right hand-side of (\ref{3.23}) is  
\begin{equation}
   O\( \frac{1}{(\log |\xi|)^{4n   +\de}}\) .
   \end{equation}

  } \end{example}

  Considering (\ref{s7.8q}) we can replace condition (\ref{7.8q}) by 
 \be
\bar\kappa(\xi)\ge  C    \va_{\al}(|\xi|),    \label{wqw}
 \ee
for all $|\xi|$ sufficiently large, where $ \bar\kappa$ is the L\'evy exponent  of $X$.
Since $ \bar\kappa(\xi)=o(|\xi|^{2})$ as $|\xi| \to \ff$, it follows that the possible values of $\al$ must be less than or equal to 2 in any dimension.  When $d=1$,    the condition that $\al\leq 1$ in (\ref{8.34}) initially appears to be very restrictive.   However,  when $\al>1$ the L\'evy process $X$ has local times, and self intersection local times can be studied without resorting to the complicated process of  renormalization.

\section{Loop soup and permanental chaos}\label{sec-loopilt} 

  Let $\mathcal{L}_{\al}$ be  the Poisson point process  on $\Om_{\De}$ with intensity measure $\al \mu$.   Note that  $\mathcal{L}_{\al}$ is a random variable; each realization  of $\mathcal{L}_{\al}$ is a countable subset of $\Om_{\De}$. To be more specific, let\begin{equation}
N(A):=\# \{\mathcal{L}_{\al}\cap A \},\hspace{.2 in}A\subseteq \Om_{\De}.\label{po.1}
\end{equation}
Then for any disjoint measurable subsets $A_{1},\ldots ,A_{n}$ of $\Om_{\De}$, the random variables $N(A_{1}),\ldots ,N(A_{n})$, are independent, and $N(A)$ is a Poisson random variable with parameter $\al \mu (A)$, i.e.
\begin{equation}
P_{\mathcal{L}_{\al}}\(N(A)=k\)={(\al \mu (A))^{k} \over k!}e^{-\al \mu (A)}.\label{po.2}
\end{equation}
(When $\mu(A)=\infty$, this means that $P(N(A) =\infty) =1$.)
We call   the Poisson point process $\mathcal{L}_{\al}$ the     `loop soup'  of   the     Markov process $X$. See \cite{LL, LF, LW, Le Jan1}.  
  
  Let $D_{m}\subseteq D_{m+1}$ be a sequence of sets in $\Om_{\De}$ of finite $\mu$ measure, such that \be
   \Om_{\De}=\cup_{m=1}^{\ff} D_{m}. \label{dm}
   \ee
   
        For $\nu \in \BB _{2n}(R^{d})$ we define   the   `loop soup   $n$-fold  self-intersection local time', $\psi_{n}(\nu)$,  by 
\begin{equation}
\psi_{n}(\nu)=\lim_{m\rar \ff}\psi_{n,m}(\nu), \label{ls.10a}
\end{equation} 
where 
\begin{equation}
\psi_{n,m}(\nu)
=\(\sum_{\om\in \mathcal{L_{\al}}} 1_{D_{m}} L_{n }(\nu)(\om)\)-\al \mu( 1_{D_{m}}L_{n }(\nu)). \label{ls.10}
\end{equation}
 (Each realization of $\LL_{\al}$ is a countable set of  elements of $\Om_{\De}$.  The expression $\sum_{\om\in \mathcal{L_{\al}}}$ refers to the sum over this set.  Because $\LL_{\al}$ itself is a random variable, $\psi_{n,m}(\nu)$ is a random variable.)

 The factor $1_{D_{m}}$ is needed to make  $\mu( 1_{D_{m}}L_{n}(\nu))$ finite. To understand this   recall that $L_{n}(\nu)\in L^{2}(\mu)$ by Theorem \ref{theo-multiriltintro-m}, hence $\mu( 1_{D_{m}}L_{n}(\nu))$ is finite.   On the other hand, using  (\ref{ls.3})   we see that  
 \bea
\lefteqn{
 \mu \( \int_{0}^{\ff}f(X_{s})\,ds\)\label{ls.2mm}}\\
 &&=\int_{0}^{\ff}{1\over t}\int\,Q_{t}^{x,x}\(\int_{0}^{t}f(X_{s})\,ds \)\,dm (x)\,dt\nn \\
 &&=\int_{0}^{\ff}{1\over t}\int\,\int_{0}^{t}\,Q_{t}^{x,x}\( f(X_{s})\)\,\,ds\,\,dm (x) \,dt\nn \\
 &&=\int_{0}^{\ff}{1\over t}\int\,\int_{0}^{t}\, \(\int p_{s}(x,y) p_{t-s}(y,x)f(y)\,dy\)\,\,ds\,\,dm (x) \,dt\nn \\
 &&=\int_{0}^{\ff}{1\over t} \,\int_{0}^{t}\, \(\int p_{t}(y,y) f(y)\,dy\)\,\,ds\,dt\nn \\
 &&=\(\int_{0}^{\ff}p_{t}(0)\,dt\)\int  f(y)\,dm (y )=\ff,\nn
 \eea 
 whenever $u(0)=\ff$.

\medskip	  The next lemma gives a formula for the joint moments of $\{\psi_{n_{j},m_{j}}(\nu_{j})\}$.  When all the $n_{j}=1$ it is essentially the same as  \cite[(2.37)]{LMR}.

  \begin{lemma} \label{lem-4.1} For all  $\nu_{j} \in \BB _{ 2n_{j}}$,  $j=1,\ldots,k$, 
    \begin{equation}
E_{\mathcal{L}_{\al}}\( \prod_{j=1}^{k}\psi_{n_{j},m_{j}}(\nu_{j}) \)= \sum_{\cup_{i } B_{i}=[1,k],\,|B_{i}|\geq 2}\,\,\prod_{i }\,\al \,\,\mu\(\prod_{j\in B_{i}}1_{D_{m_{j}}}L_{n_{j} }(\nu_{j})\)\label{mom.11},
\end{equation}
where the  sum is over all partitions $B_{1},\ldots, B_{k}$ of $[1,k]$ with parts $|B_{i}|\geq 2$. 

 \end{lemma}
  
\Proof  For $j=1,\ldots,k$, let
 \begin{equation}
Y_{j}:=1_{D_{m_{j}}}\mbox{ sign }(L_{n_{j} }(\nu_{j}))\( |L_{n_{j} }(\nu_{j})|\wedge M\),\label{masj}
\end{equation}
for some constant $M>0$.
By  the master formula for Poisson processes, \cite[(3.6)]{K}\bea
&&
E_{\mathcal{L}_{\al}}\(e^{ \sum_{j=1}^{k}z_{j}\(\sum_{\om} Y_{j}(\om)-\al\mu(Y_{j})\)}\)\label{ls.11}\\
&&\qquad=\exp \(\al\(\int_{\Om_{\De}}\( e^{\sum_{j=1}^{k}z_{j}Y_{j} }
-\sum_{j=1}^{k}z_{j}Y_{j} -1\)\,d\mu(\om)  \)\).\nn
\eea
Differentiating  each side of (\ref{ls.11}) with respect to $z_{1},\ldots, z_{k}$ and then setting $z_{1},\ldots, z_{k}$ equal to zero, we   obtain
\be
E_{\mathcal{L}_{\al}}\( \prod_{j=1}^{k}\(\sum_{\om} Y_{j}(\om)-\al\mu(Y_{j})\)\)
\label{mask}= \sum_{\cup_{i } B_{i}=[1,k],\,|B_{i}|\geq 2}\,\,\prod_{i }\,\al \,\,\mu\(\prod_{j\in B_{i}}Y_{j}\), 
\ee
where the  sum is over all partitions $B_{1},\ldots, B_{k}$ of $[1,k]$ with parts $|B_{i}|\geq 2$.
Taking the limit as $M\rar\ff$ gives (\ref{mom.11}).

 (We did not initially define  $Y_{j}=1_{D_{k_{j}}}L_{n_{j} }(\nu_{j})$   because   it is not clear that  the right hand side of (\ref{ls.11}) is finite without the truncation at $M$.)
\qed

  It follows from Lemma \ref{lem-4.1} that 
\begin{equation}
   E_{\mathcal{L}_{\al}}\(  \psi_{n ,\de }(\nu ) \)=0.
   \end{equation}

\medskip	The next Theorem asserts that we can take the limit in (\ref{mom.11}) and consequently that  (\ref{ls.10a}) exists in  $L^{p}(\mu)$ for all  $p\ge 1$.

\begin{theorem}\label{theo-9.1ilt}    Let
$X $ be as in Theorem  \ref{theo-multiriltintro-m}. If  $\nu \in \BB_{2n}(R^{d} )$,    the limit in (\ref{ls.10a}) exists in   $L^{p}(\mu)$ for all  $p\ge 1$. 

In addition,  let $n=n_{1}+\cdots+n_{k}$,     and     $\nu_{j} \in \BB_{ 2n_{j}}(R^{d})$,  $j=1,\ldots,k$. Then  
\begin{equation}
E_{\mathcal{L}_{\al}}\( \prod_{j=1}^{k}\psi_{n_{j}}(\nu_{j}) \)= \sum_{\cup_{i } B_{i}=[1,k],\,|B_{i}|\geq 2}\,\,\prod_{i }\,\al \,\mu\!\(\prod_{j\in B_{i}} L_{n_{j}}(\nu_{j})\)\label{mom.11ilt},
\end{equation}
where the  sum is over all partitions $B_{1},\ldots, B_{k}$ of $[1,k]$ with parts $|B_{i}|\geq 2$.

 \end{theorem}

  \Proof  
  We take the limit as the  $m_{j}\rar \ff$ in (\ref{mom.11}) and use Theorem \ref{theo-multiriltintro-m} to see that the right hand side of (\ref{mom.11}) converges to the right hand side of  (\ref{mom.11ilt}). 
 Applying this with $\prod_{j=1}^{k}\psi_{k_{j},\de_{j}}(\nu_{j})$ replaced by $\(\psi_{n,m}(\nu_{j})-\psi_{n,m'}(\nu_{j})\)^{k}$   shows that the limit  in (\ref{ls.10a}) exists in   $L^{p}(\mu)$ for all  $p\ge 1$. 
 \qed
    
\noindent{\bf Proof of  Theorem \ref{theo-1.4} } \label{forpageref1.4}   We show that for $n=n_{1}+\cdots+n_{k}$ and $c(\pi)$ equal to the number of cycles in the permutation $\pi$,  
 \be s
E_{\mathcal{L}_{\al}}\( \prod_{i=1}^{k}\psi_{n_{i}}(\nu_{i}) \)\label{20.15qintroq}= \sum_{\pi\in \mathcal{P}_{0}}\al^{c(\pi)}  \int  \prod_{j=1}^{n}u (z_{j},z_{\pi(j)}) \prod_{i=1}^{k}\,d\nu_{i}  (x_{i} ) ,
\ee 
 where $z_{1},\ldots,z_{n_{1}}$ are all equal to $x_{1}$,  the next $n_{2}$ of the $\{z_{j}\}$ are all equal to $x_{2}$, and so on, so that the last   $n_{k}$ of the $\{z_{j}\}$ are all equal to $x_{k}$    
and  $\mathcal{P}_{0}$ is the set of permutations $\pi$ of $[1,n]$ with cycles that alternate the variables $\{x_{i}\}$; (i.e., for all $j$, if $z_{j}=x_{i}$ then $z_{\pi(j)}\neq x_{i}$), and in addition, for each $i=1,\ldots,k$, all the $ \{z_{j}\}$ that are equal to $x_{i}$ appear in the same cycle.

The relationship in (\ref{20.15qintroq})
 is simply  a restatement of Theorem  \ref{theo-9.1ilt} that allows it to be easily compared to Theorem \ref{theo-1.3}.   To begin we use Theorem \ref{theo-1.3}  on the right-hand side of  (\ref{mom.11ilt}). One can see from Theorem \ref{theo-1.3} that the $\mu$ measure acts on cycles;  (we have $\pi(n+1)=\pi( 1  )$). Therefore, the term $\mu\!\(\prod_{j\in B_{i}} L_{n_{j}}(\nu_{j})\)$ requires that all $n_{j}$   terms, say  $x_{j}$, that are in the   variables in $L_{n_{j}}(\nu_{j})$, $j\in B_{i}$, are in the same cycle. The sum over $ \cup_{i } B_{i}=[1,k],\,|B_{i}|\geq 2$ gives all the terms one gets from permutations listed according to their cycles. Thus 
 there are $c(\pi)$ factors of $\al$. The fact that $|B_{i}|\ge 2$ means that the permutations must alternate the variables $x_{i}$, a property that is required in both Theorems  \ref{theo-1.3} and  \ref{theo-1.4}. \qed

 \begin{remark} {\rm 
 It follows from  Theorem  \ref{theo-9.1ilt} and Lemma  \ref{lem-bddhat} that
  \be
E_{\mathcal{L}_{\al}}\(  |\psi^{2k}_{n}(\nu)|^{1/ n} \)\leq (2k)! C_{\al, n}^{2k}\|  \nu \|_{2,\tau_{2n}}^{2k/n}.\label{20.16} 
\ee
where $C_{\al, n}$ is a constant, depending on $\al$ and $n$. Then, as in the proof of Lemma \ref{lem-hatcont},
for  $\nu,\mu   \in   \BB_{ 2n}(R^{d})$ 
\begin{equation}
\|\psi_{n}(\nu  )-\psi_{n}(\mu)\|_{\rho_{1/n},P_{\mathcal{L}_{\al}}}  \leq C'_{\al, n} \|\nu -\mu  \|_{2,\tau_{2n}}.\label{tranh.5ab}
\end{equation}
  Consequently all the continuity results given in Section \ref{sec-iltcont} for intersection local times, including Theorem \ref{theo-cont14.1b}, also hold for the loop soup intersection local times.
}
\end{remark}

 \section{ Isomorphism Theorem  I}\label{sec-IT}
 
Theorem \ref{theo-ljriltintro},  Isomorphism Theorem  I, contains the measure $Q ^{x,y}$ defined by
\begin{equation}
Q ^{x,y}\(1_{\{\ze>s\}}\,F_{s}\)=P^{x}\( F_{s}  \,u(X_{s},y) \)\qquad \forall  F_{s} \in b\mathcal{F}^{0}_{s},\label{10.1q}
\end{equation}
where $u$ is the potential density of the L\'evy process $X_{s}$. It is related to the  measure  $Q_{t} ^{x,y}$, introduced in (\ref{10.1}), by the equation
  \begin{equation}
Q^{x,y} (F)=\int_{0}^{\ff} Q_{t}^{x,y}\(F\circ k_{t}\) \,dt, \hspace{.2 in}F\in b\mathcal{F}^{0};\label{ls.3n}
\end{equation}
 see \cite[Lemma 4.1]{LMR}.

The next  lemma  is proved by a slight modification of the proof of  \cite[(4.10)]{LMR}.   
 
\begin{lemma} 
For any   positive measures $\rho,\phi\in \BB_{2}(R^{d})$,   and all  bounded measurable functions $f_{j}$,  $j=1,\ldots,k$, 
 \bea
 &&
\int Q^{x,x}\( L_{1}(\phi)F\(\int_{0}^{\ff}f_{1}(X_{t})\,dt,\ldots, \int_{0}^{\ff}f_{k}(X_{t})\,dt\)\)\,d\rho(x)\label{palmers}\\
&&\hspace{.5 in}=\mu\( L_{1}(\rho) L_{1}(\phi)\,\,F\(\int_{0}^{\ff}f_{1}(X_{t})\,dt,\ldots, \int_{0}^{\ff}f_{k}(X_{t})\,dt\)\) \nn,
\eea
for any bounded measurable function $F$ on $R^{k}$.
 \end{lemma}
 
  We proceed to the proof of Theorem \ref{theo-ljriltintro}.  We use a special case of the Palm formula for a Poisson process $\mathcal{L}$ with intensity measure  $\vartheta$ on a measurable space $\mathcal{S}$, see \cite[Lemma 2.3]{Bertoin}, that states that for any positive function $f$ on $\mathcal{S}$ and any measurable functional $G$ of $\mathcal{L}$
 \begin{equation}
E_{\mathcal{L}} \(\sum_{\om\in \mathcal{L}}f(\om)G(\mathcal{L})\)=
\int E_{\mathcal{L}} \(G(\om'\cup\,\mathcal{L})\)f(\om')\,d\vartheta(\om'). \label{palm.0}
 \end{equation} 
  
 \noindent{\bf Proof of Theorem \ref{theo-ljriltintro} } We show that
  for any   positive measures $\rho,\phi\in  \BB  _{2}(R^{d})$ there exists $\th^{\rho,\phi}$ such that for any finite measures $\nu_{j} \in \BB_{2n_{j}}(R^{d})$ ,  $j=1,2,\ldots$, and bounded 
measurable functions  $F$ on $R^\ff_+$,  
\begin{equation}
E_{\mathcal{L}_{\al}}\int Q^{x,x}\(L_{1}(\phi)\, F \( \psi_{n_{j}}(\nu_{j})+L_{n_{j}}(\nu_{j})\) \)\,d\rho(x)={1 \over \al}E_{\mathcal{L}_{\al}} \(\th^{\rho,\phi} \,F\(\psi_{n_{j}}(\nu_{j})\)\).\label{6.1introq}
\end{equation}

 Note that (\ref{6.2}) follows from (\ref{6.1introq}) with $F\equiv 1$.
 Set
 \begin{equation}
\th^{\rho,\phi}
=\sum_{\om\in \mathcal{L_{\al}}}L_{1}(\rho)(\om)  L_{1}(\phi)(\om).\label{j8.5}
\end{equation} 
By Theorem \ref{theo-1.3} we see that $L_{1}(\rho) L_{1}(\phi)$ is integrable with respect to $\mu$.
We use   (\ref{palm.0}) with  $\vartheta=\mu$,
\begin{equation}
f(\om)=L_{1}(\rho)(\om)  L_{1}(\phi)(\om)\label{palm.2}
\end{equation}
and  
\begin{equation}
G(\mathcal{L}_{\al})=F\(\psi_{n_{j},\de}(\nu_{j})(\mathcal{L}_{\al})\),\label{palm.2a}
\end{equation}
see (\ref{ls.10}), where  $F$ is a bounded continuous function on  $R^{k}$. Note that   since $\mu$ is non-atomic, for any fixed  $\om'$ it follows that $\om'\not\in \LL_{\al}$ a.s. Hence 
\bea
\lefteqn{
\psi_{n,\de}(\nu)(\om'\cup\mathcal{L_{\al}})
=\(\sum_{\om\in \mathcal{\om' \cup L_{\al}}} 1_{\{\ze (\om)>\de\}} L_{n}(\nu)(\om)\)-\al \mu( 1_{\{\ze>\de\}}L_{n}(\nu))\nn}\\
&&=1_{\{\ze (\om')>\de\}} L_{n}(\nu)(\om')+\(\sum_{\om\in \mathcal{L_{\al}}} 1_{\{\ze (\om)>\de\}} L_{n}(\nu)(\om)\)-\al \mu( 1_{\{\ze>\de\}}L_{n}(\nu))\nn\\
&&=1_{\{\ze (\om')>\de\}} L_{n}(\nu)(\om')+\psi_{n,\de}(\nu),\label{palm.2b}
\eea
so that 
\begin{equation}
G(\om'\cup\mathcal{L}_{\al})=F\(\psi_{n_{j},\de}(\nu_{j})+1_{\{\ze (\om')>\de\}}L_{n_{j}}(\nu_{j})(\om')\).\label{palm.3}
\end{equation}
By (\ref{palm.0}), and the fact that $\al\mu$ is the intensity measure of $\mathcal{L}_{\al}$, we see that
 \bea
&&
E_{\mathcal{L_{\al}}} \(\th^{\rho,\phi}F\(\psi_{n_{j},\de}(\nu_{j})\)\)\label{palm.4}\\
&&\qquad=\al
\int E_{\mathcal{L_{\al}}} \(F\(\psi_{n_{j},\de}(\nu_{j})+1_{\{\ze (\om')>\de\}}L_{n_{j}}(\nu_{j})(\om')\)\)\nonumber\\
&&\hspace{2.5 in}L_{1}(\rho)(\om')  L_{1}(\phi)(\om')\,d\mu(\om').\nn
 \eea
We now use (\ref{palmers}) and  take the limit as  $\de\rar 0$ to obtain   (\ref{6.1introq})  when $F$  is a bounded continuous function $R^{k}$.   The extension to general bounded 
measurable functions
$F$ on $R^\ff$ is routine.\qed

\section{ Permanental Wick powers}\label{sec-powers}

  In (\ref{rilt.1q}) we define the approximate local time of $X$ as 
\begin{equation}
L_{1}(x,r) =\int_{0}^{\ff}f_{x,r}(X_{t})\,dt. \label{rilt.1qa}
\end{equation}

  Let  $\{D_{m}\}$ be as defined in (\ref{dm}).
 Set  
 
 \be 
 \psi (x,r)  =\lim_{m\rar \ff}\(\sum_{\om\in \mathcal{L_{\al}}} 1_{D_{m}}L_{1}(x,r)(\om)\)-\al \mu( 1_{D_{m}}L_{1}(x,r)) .\label{sxr}
 \ee
 Since $L_{1}(x,r)\in L^{2}(\mu)$ by  (\ref{ls.4}), $\mu( 1_{D_{m}}L_{1}(x,r))<\ff$. (The factor $ 1_{D_{m}}$ is needed since $\mu(L_{1}(x,r))=\ff$, see (\ref{2.4}).)

 	In this section we  construct $\wt\psi_{n}(\nu)$, the analogue of the higher order Wick powers, 
  $   2^{-n} :G^{2n}: (\nu)$, from the $\al$-permanental field $\psi(\nu)$. We call $\wt\psi_{n}(\nu)$ a permanental Wick power. In Section \ref{sec-lspower} we obtain a loop soup interpretation of $\wt\psi_{n}(\nu)$.
  
\medskip	The definition of  permanental Wick powers is similar to the definition of intersection local times in Section \ref{sec-iltac}. We use chain functions, defined in (\ref{1L.1}),  and a similar quantity which we call    circuit functions
    \begin{equation}
 \mbox{ci}_{k}(r) :=    \int   u (ry_{1},ry_{2}) \cdots u (ry_{k-1},ry_{k})u (ry_{k},ry_{1})  \prod_{j=1}^{k} f(y_{j}) \,dy_{j}.\label{1L.2}
  \end{equation} 
  For any $\si=(k_{1}, k_{2},\ldots ; m_{2}, m_{3},\ldots)$  we   extend the definitions in (\ref{rilt.2a}) and set
  \begin{equation}
 |\si|=\sum_{i=1}^{\ff}  ik_{i}+\sum_{j=2}^{\ff}  jm_{j},\hspace{.2 in}  |\si|_{+}=\sum_{i=1}^{\ff}   (i+1)k_{i}+\sum_{j=2}^{\ff}  jm_{j}.\label{11.2a}
  \end{equation}
 We abbreviate  the notation by not indicating the values of $k_{i+1}$, $i\ge 1$, when it and all subsequent $k_{\cdot}$ are equal to 0, and similarly for $m_{i+2}$. For example,   $\si=(2;0,1)$ indicates  that $k_{1}=2$,  $m_{2}=0$ and $m_{3}=1$ and all other $k_{i}$ and $m_{j}=0$. In this case 
 $ |\si|=5 $ and   $  |\si|_{+}=7.$
  
  In what follows $k_{i}$ indicates the number of chains ch$_{i}$ and $m_{j}$  indicates the number of circuits ci$_{j}$.
  
  Set $\wt \psi_{1}(x,r)=\psi (x,r)$ and $\wt \psi_{0}(x,r)=1$.   Analogous to  (\ref{rilt.1}) and (\ref{rilt.2b}), in the construction  of $L_{n}(x,r)$, we define   recursively
  \begin{eqnarray}
&&
\wt\psi_{n}(x,r)  =\psi^{n} (x,r)  -\sum_{\{\si\,|\, 1\leq |\si| \leq |\si|_{+} \leq n\} } I_{n}(\si,r) \label{11.2j},
  \end{eqnarray}
  where
\begin{equation}
I_{n}(\si,r)={n! \over \prod_{i=1}^{\ff}k_{i}!\prod_{j=2}^{\ff}m_{j}!}\prod_{i=1}^{\ff} \(\mbox{ch}_{i}( r)\)^{k_{i}} 
  \prod_{j=2}^{\ff} \({\al\,\mbox{ci}_{j}( r) \over j}\)^{m_{j}}   {   \wt\psi_{n-|\si| } (x,r)\over (n-|\si|_{+})!}.\label{11.2b}
\end{equation}  
  (Here we use the notation under the summation sign in (\ref{11.2j}) to indicate that we sum over all $k_{1}, k_{2},\ldots ; m_{2}, m_{3},\ldots$ such that $\{\si\,|\, 1\leq |\si| \leq |\si|_{+} \leq n\} $.)

 It is easy to check that   
 \bea
\wt \psi_{2}(x,r)  &=&\psi^{2} (x,r)-I_{2}((1;0),r)-I_{2}((0;1),r))\label{11.4}\\
 &=&\psi^{2} (x,r)-2 \mbox{ch}_{1}( r) \psi (x,r)-\al  \mbox{ci}_{2}(r),\nn 
 \eea 
 and
\bea
 \wt \psi_{3}(x,r)&=&\psi^{3} (x,r)-I_{3}((1;0),r))\label{11.10}\\
  &&\qquad -I_{3}((0,1;0),r)) -I_{3}((0;1),r)) -I_{3}((0;0,1),r))\nn\\
  &=& \psi^{3} (x,r)-6 \mbox{ch}_{1}( r) \wt \psi_{2}(x, r)\nn\\
  && \qquad-\lc 6 \mbox{ch}_{2}( r) +3\al \mbox{ci}_{2}( r) \rc \psi (x,r) -2\al \mbox{ci}_{3}( r).\nn
\eea
Using (\ref{11.4}) we can write $ \wt \psi_{3}(x,r)$ as a third degree polynomial  in $\psi(x,r)$,
\bea
 \lefteqn{
 \wt \psi_{3}(x,r)
=\psi^{3} (x,r)-6 \mbox{ch}_{1}( r) \psi^{2} (x,r) \label{11.10y}}\\
  &&  -\lc 6 \mbox{ch}_{2}( r) +3\al \mbox{ci}_{2}( r)-12 \mbox{ch}^{2}_{1}( r)  \rc \psi (x,r) -2\al \mbox{ci}_{3}( r)+6\al \mbox{ch}_{1}( r)\mbox{ci}_{2}(r).\nn
\eea

 We show in Remark \ref{rem-6.1} that  $\wt\psi^{n}(x,r)$ can also be defined by a generating function.

\medskip	
Analogous to Theorems \ref{theo-multiriltintro-m},    \ref{theo-1.3}  and \ref{theo-1.4} we have the following results about $\wt\psi_{n}$ and its joint moments.

 \bt\label{theo-multivar}  Let
$X=\{X(t),t\in R^+\}$ be a L\'evy process in $R^d$, $d=1,2$,  that is killed at the end of an  independent exponential time, with potential density $u$ that satisfies (\ref{7.8q})  and   (\ref{8.34}) and let $\nu\in \BB_{2n}(R^{d})$, $n\ge 1$.  Then
\begin{equation}
  \wt\psi_{n}(\nu):= \lim_{r\rar 0}\int     \wt\psi_{n}(x,r) \,d\nu(x)\label{11.9}
\end{equation}
    exists in  $L^{p}(P_{\LL_{\al}})$ for all $p\geq 1$.  
\et
\bt    \label{theo-6.2}
 Let $X$ be as in Theorem \ref{theo-multivar}   and let $n=n_{1}+\cdots+n_{k}$,  and 
  $ \nu_{i}\in \BB_{2n_{i}}(R^{d} )$.  Then 
 \be 
E_{\mathcal{L}_{\al}}\(  \prod_{i=1}^{k}   \wt\psi_{n_{i}}(\nu_{i} ) \)\label{20.15} = \sum_{\pi\in \mathcal{P}_{n,a}}\al^{c(\pi)}  \int  \prod_{j=1}^{n}u (z_{j},z_{\pi(j)}) \prod_{i=1}^{k}\,d\nu_{i}  (x_{i} ),  
\ee 
 where $z_{1},\ldots,z_{n_{1}}$ are all equal to $x_{1}$,  the next $n_{2}$ of the $\{z_{j}\}$ are all equal to $x_{2}$, and so on, so that the last   $n_{k}$ of the $\{z_{j}\}$ are all equal to $x_{k}$ and  $\mathcal{P}_{n,a}$ is the set of permutations $\pi$ of $[1,n]$ with cycles that alternate the variables $\{x_{i}\}$; (i.e., for all $j$, if $z_{j}=x_{i}$ then $z_{\pi(j)}\neq x_{i}$).    \et
  
  It is interesting to note that a simple version of (\ref{20.15}) appears  in   \cite[( 1.1)]{LMR} for    permanental fields and  in \cite[Proposition 4.2]{VJ} for permanental random variables, (where they are referred to as multivariate gamma distributions).  
   
\medskip	The following corollary shows that the processes   $\{\wt\psi_{n}(\nu) \}_{n=1}^{\ff}$ are orthogonal with respect to $P_{\LL_{\al}}$.
   
 \begin{corollary}\label{cor-var}  Let $X$ be as in Theorem \ref{theo-multivar}. 
  If $\nu,\nu'\in \BB_{2n}(R^{d})$, then for any $k\leq n$,
 \be
E_{\mathcal{L}_{\al}}\(    \wt \psi_{n}(\nu ) \,\,    \wt\psi_{k} (\nu' )\)\label{0.19}=\de_{n,k} \,\,\, K(\al,n)\,\int \(  u (x,y)u (y,x)  \)^{n}  \,d\nu  (x )\,d\nu'  (y) 
\ee 
 where
 \begin{equation}
 K(\al,n)=n!\al(\al+1)\cdots (\al+n-1).\label{0.20}
 \end{equation} 
 \end{corollary}
 
 \begin{remark}  {\rm 
Theorem \ref{theo-6.2} holds when $k=1$, in which case the right-hand side of   (\ref{20.15}) is empty. Therefore,   
\begin{equation}
   E\( \wt\psi_{n_{i}}(\nu_{i} )\)=0. 
   \end{equation}
   When $k=n$, the right-hand side of (\ref{0.19}), with $\de_{n,k}=1$, is the covariance of $\{\psi_{n}(\nu),\nu\in \VV\}$. To see that (\ref{0.19}) has the same form as many of our other results we note that
   \begin{equation}
   \int \(  u (x,y)u (y,x)  \)^{n}  \,d\nu  (x )\,d\nu  (y)=\frac{1}{(2\pi)^{nd}}\int \th(\xi)|\hat{\nu}(\xi)|^{2}\,d\xi
   \end{equation}
   where $\th(\xi)$ is the Fourier transform of $\(  u (x,y)u (y,x)  \)^{n} $. Also note that $u (x,y)\newline u (y,x)\ge 0$ and is symmetric.
 }\end{remark}

\noindent  {\bf  Proof  of  Theorems \ref{theo-multivar} and  \ref{theo-6.2}  }  The proof is very similar to the proof of Theorem \ref{theo-1.3}, so we shall not go through all the details.  The main difference is the possibility of circuits.   They introduce the factors $ \mbox{ci}_{k}(r)$.
 
   In analogy with (\ref{20.15s}) in the proof of Theorem \ref{theo-1.3}
 we first show that for   any  $n\ge 2$ and $m\ge 1$, 
  \bea
\lefteqn{
E_{\mathcal{L}_{\al}}\(\wt \psi_{n }(x,r)   \prod_{i=1}^{m}\psi (f_{i})\)\label{20.15sp}}\\
&&\qquad =\sum_{\pi\in \mathcal{P}_{m,n}}\al^{c(\pi)}  \int  \prod_{j=1}^{m+n}u (z_{j},z_{\pi(j)})      \prod_{i=1}^{m} f_{i}(z_{i} )\,dz_{i} \prod_{i=m+1}^{m+n}f_{r,x}  (z_{i } )\,dz_{i }\nn\\
&&  \hspace{2in}+\EE_{{\bf f }}(E_{r}(x,  {\bf z })) \nn,
\eea
  where,  $\mathcal{P}_{m,n}$ is the set of permutations $\pi$ of [1,m+n] such that 
  $\pi:[m+1,m+n]\mapsto [1,m]$,  (when $m<n$ there are no such permutations), and 
 the last term is as given in (\ref{2.13}), although the terms in      $ E_{r}(x,  {\bf z }) $ are not the same. 
 We show  below that   for   $\nu\in\BB_{2n}(R^{d})$,   
  \begin{equation}
 \lim_{r\rar 0} \sup_{\forall |z_{i}|\le M}\int E_{r}(x, {\bf z }) \,d\nu(x)=0,\label{err0a}
  \end{equation}
  which implies that 
  \begin{equation}
    \lim_{r\rar 0}\int  \EE_{{\bf f }}(E_{r}(x,  {\bf z }))\,d\nu(x)=0.
   \end{equation} 
 
 The proof of  (\ref{20.15sp}) when $n=1$  
follows from (\ref{mom.11ilt}) and (\ref{rilt.15intro}), and in this case the error term $\EE_{{\bf f }}(E_{r}(x,  {\bf z }))=0$.

  For the proof of   (\ref{20.15sp}) when $n=2$ we note that by Theorem \ref{theo-1.4}
   \bea
\lefteqn{
E_{\mathcal{L}_{\al}}\(  \psi^{2}(x,r)   \prod_{i=1}^{m}\psi (f_{i})\)\label{20.15spw}}\\
&&\qquad =\sum_{\pi\in \mathcal{P}_{0}}\al^{c(\pi)}  \int  \prod_{j=1}^{m+2}u (z_{j},z_{\pi(j)})      \prod_{i=1}^{m} f_{i}(z_{i} )\,dz_{i} \prod_{i=m+1}^{m+2}f_{r,x}  (z_{i } )\,dz_{i }\nn,
\eea  
where, because all the $n_{i}=1$, $\PP_{0}$ is the set of permutations on $[1,m+2]$. Consider the  permutations in $\PP_{0}$ that do not take  $\pi:[m+1,m+2]\mapsto [1,m]$. They can be divided into three sets. Those that have the cycle   $(m+1,m+2)$,  those that take only $m+1\mapsto [1,m]$ and  those that take only $m+2\mapsto [1,m]$. Taking these into consideration  and considering (\ref{11.4}), we get (\ref{20.15sp}) when $n=2$. This is obvious for the terms involving the cycle. To understand  this for the terms involving the two chains consider 
(\ref{2.17a})--(\ref{2.23q}) and note that an extra density function is introduced. It is this that gives the factor $\psi(x,r)$.

 \medskip	
 Assume that (\ref{20.15sp}) is proved for $\psi_{n'}(x,r)$, $n'<n$. For any  $\si=(k_{1}, k_{2},\ldots ; \newline	m_{2}, m_{3},\ldots)$  let $\PP_{ 0}  (\si)$ denote the set of permutations $\bar\pi\in\PP_{ 0}$  
that contain     $m_{j}$ circuits of order $j =2,3,\ldots$, and $k_{i}$
chains of order $i=1,2,3,\ldots$, in $[m+1,m+n]$.  Note that  $\mathcal{P}_{0} -\mathcal{P} _{m,n }=\cup_{|\si|\geq 1}\mathcal{P}_{0}(\si).$ In this rest of this proof we make these definitions more explicit by  writing $\mathcal{P}_{0}(\si)$ as 
$\mathcal{P}_{0,m+n}(\si)$.

Any  term $\bar\pi\in\PP_{ 0,m+n}  (\si)$ in the evaluation of 
\be
E_{\mathcal{L}_{\al}}\(  \psi   ^{n}(x,r)   \prod_{i=1}^{m}\psi (f_{i})\)\label{2.20s}
 \ee 
is the same as   the term in (\ref{20.15sp}) for a particular permutation $\pi'\in \mathcal{P} _{m,n-|\si| }$     in the evaluation of 
\begin{equation}
E_{\mathcal{L}_{\al}}\(   
\prod_{i=1 }^{\ff} \(\mbox{ch}_{i}( r)\)^{k_{i}} 
  \prod_{j =2}^{\ff}  \(\al\,\mbox{ci}_{j}( r)  \)^{m_{j}}  \wt  \psi_{n-|\si|}(x,r)\prod_{i=1}^{m}\psi (f_{i})\). \label{art}
\end{equation}
  This follows by the same argument given in (\ref{2.20})--(\ref{2.21q}) except that here we use  Theorem \ref{theo-1.4} and (\ref{20.15sp}).

 Similar to the analysis on page \pageref{page11}  the permutation $\pi' $ is obtained from  $\bar\pi$ by the   remove and relabel  technique.  However there is a significant difference in this case. In Theorem \ref{theo-1.3} we work with the loop measure $\mu$ so that the permutation $\bar\pi$ consists of a single cycle. (We have $\pi(n+1)=\pi(1)$.) In this theorem we take the expectation with respect to 
 $E_{\mathcal{L}_{\al}}$.  As one can see in (\ref{mom.11ilt}),  the permutation $\bar\pi$   generally  has many cycles, each of order greater than or equal to 2.  
 
  We illustrate this with an example similar to (\ref{am1}). Consider the case when   $m=10 $, $n=14$ and the permutation on $[1,24]$ given by
\begin{equation}
\bar\pi=(6, 7, 11,13, 8, 9,10)(1,14, 12,16, 2, 3, 4, 15, 5, 17, 18)(19,20,21),(22,23,24).\label{4am1}
\end{equation}
in which we use standard  cycle notation.
$\bar\pi$ has   4   cycles.  We are primarily concerned with the effects of $\bar\pi$ on  $[11,24]$. There are three chains 
\be (11,13)\qquad(14,12,16)\qquad(17, 18),\label{42.30q}
\ee 
and two circuits 
\be (19,20,21)\qquad(22,23,24),\label{42.30qc}
\ee
so that $\si=(2,1;0,2)$, $|\si|=10$ and $|\si|_{+}=13$.

We first remove all circuits and all but the first element in each chain   in (\ref{4am1}) to obtain
\begin{equation}
 (6, 7, 11, 8, 9,10)(1,14,  2, 3, 4, 15, 5,17).\label{4am2}
\end{equation}
  The permutation  $\pi' $ is obtained from (\ref{4am2}) by relabeling the remaining elements in $(11,\ldots,24)$ in increasing order from left to right, i.e.,
\begin{equation}
 \pi'=(6, 7, 11, 8, 9,10)(1,12,  2, 3, 4, 13, 5, 14).\label{4am2q}
\end{equation}

Let  $  \mathcal{P}_{0,24}(\si)_{\pi'}  $ denote the   permutations in $\mathcal{P}_{0,24}$   that have the number of chains and circuits designated by $\si$  and  that  give rise to $\pi'$ by the above procedure. We   compute $| \mathcal{P}_{0,24}((2,1 ;0, 2) )_{\pi'}|$. 
Each of the  $14!$ permutations of the elements $(11,12,\ldots, 24)$   give rise to distinct permutations in $ \mathcal{P}_{0,24}((2,1   ;0, 2) )_{\pi'}$, except for the $3$ rotations in each circuit and the interchange of the two circuits. We call these  internal permutations. There are $14! / 2!3^{2}$ of them. Furthermore, we consider     the single integer $15$, that  is not contained in a chain or cycle  of  the permutations in $ \mathcal{P}_{0,24}((2,1 ;0, 2) )_{\pi'} $,    to be a chain of order zero. We add this chain to the three in (\ref{42.30q}) and  consider that the permutations in $ \mathcal{P}_{0,24}((2,1   ;0, 2) )_{\pi'}$ contain four chains. Clearly, each of the $4!$ arrangements  of these four chains  correspond to distinct permutations in $ \mathcal{P}_{0,24}((2,1 ;0, 2 ))_{\pi'}$. However, we do not want to count the  interchanges of the two chains of order one, since they have already been counted in the  internal permutations. Consequently  
 \be
| \mathcal{P}_{0,24}((2,1 ;0,2 )_{\pi'}|={14!4! \over 2!2!3^{2}} .
 \ee

\medskip	
For general $\si$ and $\pi'\in \mathcal{P}_{0,m+n-|\si| }$, in which,  as in the  example above, the integers $m+1,\ldots, m+n-|\si| $ appear in increasing order,  
\begin{eqnarray}
&&| \mathcal{P}_{0,m+n }(\si)_{\pi'}| ={n! \over \prod_{j=2}^{\ff}\(m_{j}!  j^{m_{j}}\) }\,  { (n-|\si|)!\over \prod_{i=1}^{\ff}k_{i}!\,(n-|\si|_{+})!}.\label{41.20q}
\end{eqnarray} 
 To see this first note that there are $n!/\prod_{j=2}^{\ff}\(m_{j}!j^{m_{j}}\)$ internal permutations.  Here we  divide   by the $m_{j}!$ interchanges of circuits of order $j$, and the $j$ rotations  in each circuit of order $j$, for each $j$.  Note that for any $\bar\pi\in \mathcal{P}_{0,m+n }(\si)_{\pi'}$  there are $|\si|_{+}$  integers from $\{m+1,\ldots, m+n\}$ in the circuits and chains  of order $1,2,\ldots$.  Consequently, there are   $n-|\si|_{+}$ remaining integers in  $\{m+1,\ldots, m+n\}$ which, as in the example above,  we consider to be   chains of order $0$. Therefore, total number of   chains, including those of order $0$, in $\{m+1,\ldots, m+n\}$ is   
$n-|\si|_{+}+\sum_{i=1}^{\ff}k_{i}=n-|\si|$. Thus any of the $(n-|\si|)!$ permutations  of these chains 
in  $\bar\pi$ are in $\mathcal{P}_{0,m+n }(\si)_{\pi'}$. However, we do not want to  count the $(n-|\si|_{+})!\prod_{i=1}^{\ff}k_{i}! $ interchanges of chains of the same order, since this is counted in the  internal permutations.
  Putting all this  together gives 
(\ref{41.20q}). 

  Consider  (\ref{art}) again and the particular permutation $\pi'\in \mathcal{P}_{m+n-|\si| }$. We have  pointed out before  that there are $(n-|\si|)!$ different permutations,  the internal permutations,  in $\mathcal{P}_{m+n-|\si| } $,  whose contribution to (\ref{art})
is the same as it is for  $\pi'$. Therefore  up to the error terms, the contribution to (\ref{20.15sp}) from $\mathcal{P}_{m+n }(\si)$ is equal to 
\bea
&& 
{n! \over \prod_{i=1}^{\ff}k_{i}!  \prod_{j=2}^{\ff}m_{j}! j^{m_{j}}(n-|\si|_{+})!}   \nn\\
&&
\qquad \times \,\,E_{\mathcal{L}_{\al}}\(   
\prod_{i=1 }^{\ff} \(\mbox{ch}_{i}( r)\)^{k_{i}} 
  \prod_{j =2}^{\ff}  \(\al\,\mbox{ci}_{j}( r)  \)^{m_{j}}  \wt  \psi_{n-|\si|}(x,r)\prod_{i=1}^{m}\psi (f_{i})\) \nn\\
&&\qquad=E_{\mathcal{L}_{\al}}\(  I_{n}(\si,r)\prod_{i=1}^{m}\psi (f_{i})\).\label{42.21s}
   \eea 
 
 The rest of the proof follows as in the proof of Theorems \ref{theo-multiriltintro-m}  and \ref{theo-1.3}. (In controlling the error terms we also use    Lemma \ref{lem-chat}.)\qed

\medskip	\noindent{\bf  Proof of Corollary \ref{cor-var}  } It follows easily from the proof  of Theorem \ref{theo-multivar}
that $E_{\mathcal{L}_{\al}}\(  \wt \psi_{n}(\nu ) \,\,    \wt \psi_{k} (\nu )\)=0$ when $n\ne k$. When $n=k$ we have   
\be
E_{\mathcal{L}_{\al}}\(     \wt \psi_{n}(\nu ) \,\,   \wt  \psi_{n} (\nu )\)\label{20.14} = \sum_{\pi\in \mathcal{P}_{2n,a}}\al^{c(\pi)}\int    \prod_{j=1}^{2n}u (z_{j},z_{\pi(j)}) \,d\nu  (x )\,d\nu  (y). 
\ee  
 To evaluate this, we first note that since the $x$  and $y$ terms alternate  in (\ref{20.14})  it is equal to 
 \begin{equation}
 \sum_{\pi\in \mathcal{P}_{2n,a}}\al^{c(\pi)} \,\int \(  u (x,y)u (y,x)  \)^{n}  \,d\nu  (x )\,d\nu  (y)\label{20.14a}.
 \end{equation}
It remains to show that
  \begin{equation}
 \sum_{\pi\in \mathcal{P}_{2n,a}}\al^{c(\pi)}=n!\al(\al+1)\cdots (\al+n-1).\label{20.14b}
 \end{equation}
   To see this,  consider an arrangement of \[\{x_{1},\ldots, x_{n}, y_{1},\ldots, y_{n}\}\]
 into (oriented) cycles, such that each cycle contains an equal number of $x$ and $y$ terms in an alternating  arrangement. For each such arrangement we define a permutation $\si$ of $[1,n]$ by setting $\si(i)=j$   if   $x_{i}$ is followed by $y_{j}$.   We refer to the $n$ ordered pairs $(x_{1}, y_{\si(1)}),\ldots, (x_{n}, y_{\si(n)})$ as the pairs generated by $\si$.

  Let $\phi(\si,l)$ denote the number of permutations in $\mathcal{P}_{2n,a}$  with $l$ cycles that are obtained by a rearrangement   of the pairs  generated by  $\si$.  In the next paragraph we show that
 \begin{equation}
 \sum_{l}\al^{l}\phi(\si,l)=\al(\al+1)\cdots (\al+n-1).\label{0.100}
 \end{equation}
 Since each of the $n!$ permutations  of $[1,n]$ gives a different $\si$, we get(\ref{20.14b}).
 
	To prove (\ref{0.100}) we construct the  rearrangements of the pairs  generated by  $\si$ and consider how many cycles each one contains. We begin with  the ordered pair $(x_{1}, y_{\si(1)})$ which we consider as an incipient cycle,  $x_{1}\rar y_{\si(1)}\rar x_{1}$. (We say incipient because as we construct the rearrangements of the pairs  generated by  $\si$,   $(x_{1},y_{\si(1)})$ is sometimes  a cycle and sometimes part of a larger cycle.)  We next take  the ordered pair $(x_{2},y_{\si(2)})$ and use it to write
 \begin{equation}
 (x_{1},y_{\si(1)}) (x_{2},y_{\si(2)})\quad\mbox{and }\quad  (x_{1},y_{\si(1)},x_{2},y_{\si(2)}).\label{6.37}
   \end{equation} 
  We consider that the first of these terms contains two cycles  and the second one cycle. Therefore,  so far, we have accumulated  an $\al^{2}$ and an $\al$  towards the factor $\al^{c(\pi)}$ for the cycles in the rearrangements of the pairs  generated by  $\si$ that we are constructing  and   write $\al^{2}+\al=\al(\al+1)$.
  Similarly,  we use $(x_{3},y_{\si(3)})$  to extend the terms in (\ref{6.37}) as follows:
  \be  
   (x_{1},y_{\si(1)}) (x_{2},y_{\si(2)})(x_{3},y_{\si(3)})\qquad   (y_{\si(1)},x_{1},x_{2},y_{\si(2)})(x_{3},y_{\si(3)})\label{6.38}
 \ee 
 \[
   (x_{1},y_{\si(1)},x_{3},y_{\si(3)}) (x_{2},y_{\si(2)})\qquad(x_{1},y_{\si(1)})(x_{2},y_{\si(2)},x_{3},y_{\si(3)})
 \]
 \[
  (x_{1},y_{\si(1)},x_{3},y_{\si(3)}, x_{2},y_{\si(2)})\qquad (x_{1},y_{\si(1)}, x_{2},y_{\si(2)},x_{3},y_{\si(3)})
 \]
The accumulated $\al $ factors are now equal to $\al(\al+1)(\al+2)=\al^{3}+3\al^{2}+2\al$. The $\al^{3}$ comes from the first term in (\ref{6.38}) which has $3$ cycles, the $3\al^{2}$ comes from  the next three terms which have have two cycles each, and the $2\al$ comes from the last  two terms each of which has a single cycle. It should be clear now that when we add the term $(x_{4},y_{\si(4)})$ to this we multiply the 
 previous accumulated $\al $ factor by $(\al+3)$. The $\al$ in $(\al+3)$ because we can add $(x_{4},y_{\si(4)})$ to each of the terms in (\ref{6.38}) as a separate cycle.  The factor $3$ because we can place $(x_{4},y_{\si(4)})$ to the right of each $(x_{i},y_{\si(i)})$, $i=1,2,3$, in each of the other terms in (\ref{6.38}), without changing the number of cycles they contain.  Proceeding in this way gives (\ref{0.100}).\qed

 \begin{corollary} \label{cor-hyp}   
Let   $ \nu\in\BB_{ 2n}$, then for any $k$ and any $\al>0$,  
 \be
\big | E_{\mathcal{L}_{\al}}\(  \wt\psi^{k}_{n}(\nu) \)\big |\leq (kn)!C^{kn}_{\al}\|  \nu \|_{2,\tau_{2n}}^{k}.\label{20.15a} 
\ee
 \end{corollary}

\noindent	{\bf  Proof   } In this case the integral in  (\ref{20.15}) is equal to
    \begin{equation}
  \int  \prod_{j=1}^{kn}u (z_{j},z_{\pi(j)}) \prod_{i=1}^{k}\,d\nu  (x_{i} )  \end{equation}
in which there are exactly  $2n$ factors of $u$ containing $x_{i}$, for each $i=1,\ldots,k$. 
It follows from Lemma \ref{lem-bddhat}  that 
    \begin{equation}
\Big | \int  \prod_{j=1}^{kn}u (z_{j},z_{\pi(j)}) \prod_{i=1}^{k}\,d\nu  (x_{i} )\Big | \leq   \|  \nu \|_{2,\tau_{2n}}^{k}.\label{20.0}
 \end{equation}
  Since there are $(kn)!$  unrestricted permutations of  $[1,kn]$  we get (\ref{20.15a}).\qed
  
  \begin{theorem} Let $\{\wt\psi_{n}(\nu),\nu\in\VV \}$ be an $n$-th order permanental Wick power, where $\VV \in \mathcal {B}_{2n}(R^{d})$,  $n\ge 2$,  and let     $ \bar d((\nu,\mu))=  \|\nu  -\mu  \|_{2,\tau_{2n}}$.  If (\ref{2.5}) holds  $\{\wt\psi(\nu),\nu\in\VV \}$ is continuous  on $(\VV,\bar d)$, $P_{ \LL_{\al}}$ almost surely.
\end{theorem}

\Proof   We use (\ref{20.15a}) and the inequalities given in the  transition from (\ref{3.18a}) to (\ref{3.19}) to see that 
\begin{equation}
\|\wt \psi_{n}(\nu  )-\wt \psi_{n}(\mu)\|_{\rho_{1/n},  P_{ \LL_{\al}}}\leq C   \|\nu -\mu  \|_{2,\tau_{2n}}. \label{tranh.5ac}
\end{equation}
The theorem now follows from Theorem \ref{maj}.  \qed

 In  Corollary \ref{cor-var} we give a  formula for    $E_{\mathcal{L}_{\al}}( ( \wt \psi_{n}(\nu ))^{2} )$.   
We now give an alternate expression for this expectation  which we use in the proof of Theorem \ref{theo-multiloop}.  Let $\{A _{p}\}_{p=1}^{k} $  be a partition   of $[1,n]$.  Let  $m_{i}(\{A _{p}\}_{p=1}^{k})$ be the number of sets   in this partition with $|A_{p}|=i$. We define the degree of the partition to be  
\begin{equation}
d \( \{A _{p}\}_{p=1}^{k}  \):=\(m_{1}\(\{A _{p}\}_{p=1}^{k}\),  \ldots,m_{k}\(\{A _{p}\}_{p=1}^{k} \)\) .\label{adj.1}
\end{equation}

\begin{lemma} \label{lem-6.1}  For $\nu\in\BB_{2n}(R^{d})$,
 \begin{eqnarray}
\lefteqn{E_{\mathcal{L}_{\al}}( ( \wt \psi_{n}(\nu ))^{2} ) \label{adj.3q}}\\
&& = \sum_{\{A_{1}\cup\cdots\cup A_{k}=[1,n]\}} \sum_{\{B_{1}\cup\cdots\cup B_{k'}=[1,n]\}}\de_{\{ d  ( \{A _{p}\}_{p=1}^{k}   ) ,d  ( \{B _{p}\}_{p=1}^{k'}   )\}} \prod_{i=1}^{k }  m_{i}( \{A _{p}\}_{p=1}^{k} )!  
\nn \\
&& \qquad\qquad\lim_{r\rar 0}\int \prod_{l=1}^{k}\al \mu\(L_{|A_{l}|}(x,r)L_{|A_{l}|}(y,r)\) \,d\nu  (x )\,d\nu  (y),\nn
\end{eqnarray}
where  $\de_{\{ d  ( \{A _{p}\}_{p=1}^{k}   ) ,d  ( \{B _{p}\}_{p=1}^{k'}   )\}}$ is equal to one when the  vectors $\ d  ( \{A _{p}\}_{p=1}^{k}   )=  d  ( \{B _{p}\}_{p=1}^{k'}   ) $, (so that  $k=k'$), and is equal to zero otherwise.  

Moreover, the third line of (\ref{adj.3q}) is equal to 
\begin{equation}
    \al ^{k}\prod_{l=1}^{k} |A_{l}|! ( |A_{l}|-1)!       \(\int \(  u (x,y)u (y,x)  \)^{|A_{l}|}  \,d\nu  (x )\,d\nu  (y)\).\label{6.46}
   \end{equation}
 \end{lemma}
 
 \Proof    
To understand (\ref{adj.3q})  consider  Theorem \ref{theo-6.2} with $n_{1}=n_{2}=n$, and replace the products of the potentials on the right-hand side of  (\ref{20.15}) by the limit,  as $r$ goes to zero, of (\ref{2.27w}). The sum in (\ref{20.15})   is a sum over permutations with cycles that alternate two pairs of  $n$ variables, which we denote by $x$ and   $y$. These permutations divide 
$[1,n]$ and $[n+1,2n]$,  which we identify with another copy of $[1,n]$,  into two partitions $\mathcal{A}$ and $\mathcal{B}$ of $[1,n]$ which must have the same degree.  Moreover  each set $A_{l}$ in $\mathcal{A}$ must be paired with a set $B_{l'}$ of $\mathcal{B}$ of the same cardinality.   There are $e(\mathcal{A})$ ways to make such a pairing. If $\mathcal{A}$ consists of  $k$ sets then the permutation $\pi$ has $k$  cycles,  with the $l$-th  cycle alternating the elements of $A_{l}$ and $B_{l'}$, for some $1\le l'\le k'$.

We use (\ref{2.27w}) to get (\ref{6.46}). \qed

\begin{remark}\label{rem-6.1}  {\rm 
 The process   
  $\wt\psi_{n}(x,r)=A_{n}\(\psi (x,r)\)$,  $n\ge 1$,  where the polynomials $A_{n}(u)$ satisfy 
  \begin{equation}
  \sum_{n=0}^{\ff}{\(   \sum_{i=0}^{\ff}\mbox{ch}_{i}( r)\,s^{i}+ \sum_{j=2}^{\ff}\al\,(\mbox{ci}_{j}( r)/j)\,s^{j-1} \)^{n} \over n!}s^{n}A_{n}(u)= e^{su},\label{ljoy}
  \end{equation}
   with   $\mbox{ch}_{0}( r)=1$.  
  To prove that this  this agrees with (\ref{11.2j}) and (\ref{11.2b}) we need only minor modifications of the proof in Remark \ref{rem-2.2}.  As in (\ref{slum1}) we have
   \bea
&& \(   \sum_{i=0}^{\ff}\mbox{ch}_{i}( r)\,s^{i}+ \sum_{j=2}^{\ff}\frac{\al\,\mbox{ci}_{j}( r)}{j} \,s^{j-1} \)^{n} {s^{n}\over n!}\label{slum3} 
\\
     &&\qquad=\sum_{\sum _{i=0}^{\ff}k_{i}+\sum_{j=2}^{\ff} m_{j}=n}
   {1 \over n!} {n\choose k_{0}\,\,k_{1}\,\cdots\,m_{2}\,\,m_{3}\,\cdots}\nn\\
   && \hspace{.9 in}\prod_{i=0}^{\ff} \(\mbox{ch}_{i}( r)\,s^{i}\)^{k_{i} }\prod_{j=2}^{\ff} \(\frac{\al\,\mbox{ci}_{j}( r)}{j}\,s^{j-1}\)^{m_{j}}\,s^{n}\nn
   \eea
   \bea
      &&\qquad=   \sum_{ \sum _{i=0}^{\ff}k_{i}+\sum_{j=2}^{\ff} m_{j}=n}
   \,{1 \over k_{0}!\prod_{i=1}^{\ff}k_{i}! \prod_{j=2}^{\ff}m_{j}!}\nn\\
   &&\hspace{.9 in}\prod_{i=1}^{\ff} \(\mbox{ch}_{i}( r)\)^{k_{i} }\prod_{j=2}^{\ff} \(\al\,(\mbox{ci}_{j}( r)/j)\)^{m_{j}}s^{\sum_{i=0}^{\ff}(i+1)k_{i}+\sum_{j=2}^{\ff}jm_{j}}.\nn
   \eea

  Fix $N$ and consider all $k_{0}, k_{1}, \ldots, m_{2},\,m_{3},\,\ldots$ with
\be
    N=\sum_{i=0}^{\ff}(i+1)k_{i}+\sum_{j=2}^{\ff}jm_{j}=k_{0}+| \si|_{+}.
\ee
Consequently we  can replace  $k_{0}$ with $N-| \si |_{+}$ in (\ref{slum3}). 

In addition, when  $\sum_{i=0}^{\ff} k_{i}+\sum_{j=2}^{\ff} m_{j}=n$, we have 
   $n=k_{0}+| \si|_{+}-| \si |=N-| \si |$.    We use this observation and (\ref{slum3}) to equate the coefficients of $s^{N}$  in (\ref{ljoy}) to  obtain
    \bea
   &&
\sum_{\{ \si\,|\,  | \si|_{+} \leq N\} }   {1 \over \prod_{i=1}^{\ff}k_{i}! \prod_{j=2}^{\ff}m_{j}!(N-| \si|_{+})!}\label{ljoy2x}\\
&&\hspace{1 in}\prod_{i=1}^{\ff} \(\mbox{ch}_{i}( r)\)^{k_{i} }   \prod_{j=2}^{\ff} \(\al\,(\mbox{ci}_{j}( r)/j)\)^{m_{j}}    A_{N-|\si|} (u)=   { u^{N}  \over N!},\nn
   \eea
  in which we use the expanded definitions of $|\si|$ and $|\si|_{+}$ in (\ref{11.2a}).   Setting $u=\psi (x,r)$ and  $\wt\psi _{n}(x,r)=A_{n}\(\psi (x,r)\)$ we get  (\ref{11.2j}) and  (\ref{11.2b}). 
 } \end{remark}

\section{Poisson chaos decomposition, I}\label{sec-Poisson}

 We continue our study of $\mathcal{L}_{\al}$,  a   Poisson point process  on $\Om_{\De}$ with intensity measure $\al \mu$,  and obtain a decomposition of $L^{2}(P_{\LL_{\al}})$ into orthogonal function spaces (in (\ref{11.0r})) which we refer to as the Poisson chaos decomposition of  $L^{2}(P_{\LL_{\al}})$. This is used in Section    \ref{sec-lspower}  to obtain a definition of $\wt\psi_{n}(\nu)$ in terms of loop soups. Some of the results in this section generalize results in \cite[Chapter 5.4]{Le Jan1}  which considers 
processes with   finite state spaces.  (Also, although we do not pursue this further, we mention that the results of this section as well as those of Section \ref{sec-expP} are valid for any Poisson process with diffuse intensity measure.)
 
 As explained at the beginning of Section \ref{sec-loopilt}, each realization of $\LL_{\al}$ is a countable set of  elements of $\Om_{\De}$.  The expression $\sum_{\om\in \mathcal{L_{\al}}}$ refers to the sum over this set.  Because $\LL_{\al}$ itself is a random variable such a sum is also a random variable.

   For any set $A$, we  use $S_{n}(A)\subset A^{n}$ to denote the  subset of $A^{n}$ with distinct entries. That is, if   $(a_{i_{1}},\ldots, a_{ i_{n}})\in S_{n}(A)$ then $a_{i_{j}}\neq a_{i_{k}}$ for $i_{j}\neq i_{k}$.

Let $\DD_{k} (\mu)=\cap_{p\geq k}  L^{p}(\mu)$.  
Note that $1_{B}\in \DD_{1} (\mu)$ for any $B\subseteq \Om_{\De}$ with finite $\mu$ measure. For $g_{j}\in\DD_{1} (\mu)$, $j=1,\ldots, n$, we define 
\begin{equation}
\phi_{n}(g_{1},\ldots, g_{n})=\sum_{ (\om_{i_{1} },\ldots,\om_{i_{n} })\in S_{n}(\LL_{\al}) }\,\,\prod_{j=1}^{n}g_{j}(\om_{i_{j} }).\label{11.0}
\end{equation}
  In particular
\begin{equation}
\phi_{1}(g_{1})=\sum_{  \om  \in  \LL_{\al} } g (\om ),\label{11.1}
\end{equation}
 and   $\phi_{1}(1_{B})=N(B)$; (see (\ref{po.1})). 

For any subset $A\subseteq Z_{+}$ set $g_{A}(\om) =\prod_{j\in A }g_{j}(\om)$.  In this notation the factors $g_{j}$, $j\in A$,  are functions of the same variable.   Thus, in contradistinction  to (\ref{11.0}) we have,     \begin{equation}
\phi_{1}(g_{[1,n]})= \sum_{  \om  \in  \LL_{\al} } \,\,\prod_{j=1}^{n}g_{j}(\om 
).\label{11.0a}
\end{equation}
Note that 
\begin{eqnarray}
\prod_{j=1}^{n}\phi_{1}(g_{j})&=&\sum_{ (\om_{i_{1} },\ldots,\om_{i_{n} })\in  \LL_{\al}^{n} }\,\,\prod_{j=1}^{n}g_{j}(\om_{i_{j} })
\label{11.0b}\\
&=& \sum_{k=1}^{n}\,\, \sum_{\cup_{l=1}^{k}A_{l}=[1,n]}\,\, \sum_{ (\om_{i_{1} },\ldots,\om_{i_{k} })\in S_{k}(\LL_{\al}) }\,\,\prod_{l=1}^{k}g_{A_{l}}(\om_{i_{l} })\nonumber\\
&=& \sum_{k=1}^{n}\,\, \sum_{\cup_{l=1}^{k}A_{l}=[1,n]}\,\,\phi_{k}(g_{A_{1}},\ldots, g_{A_{k}}),\nonumber
\end{eqnarray}
in which all $\{A_{l}\}_{l=1}^{k}$, $1\le k \le n$, are partitions of $[1,n]$. 
For example,  
\begin{equation}
\prod_{j=1}^{2}\phi_{1}(g_{j})=\phi_{2}(g_{1},g_{2})+\phi_{1}(g_{\{1,2\}}).\label{11.0bb}
\end{equation}

The set of random variables  $\phi_{1}(1_{B})=N(B)$, for  $B\subseteq \Om_{\De}$ with finite $\mu$ measure, generate the $\si$-algebra for $P_{\mathcal{L}_{\al}}$ By the master formula for Poisson processes, \cite[(3.6)]{K}, the random variables $N(B)$ are exponentially integrable. Hence the polynomials in $\phi_{1}$ are dense in $L^{2}\(P_{\mathcal{L}_{\al}}\)$.   Let $\phi_{0}\equiv 1$ and  let $\HH_{n}$ denote the closure in $L^{2}\(P_{\mathcal{L}_{\al}}\)$ of all linear combinations of $\phi_{j}$, $j=0,1,\ldots, n$. Then it follows from (\ref{11.0b})    that
\begin{equation}
L^{2}\(P_{\mathcal{L}_{\al}}\)=\overline{\cup_{n=0}^{\ff } \HH_{n}}.    \label{11.0c}
\end{equation}

  For  integers $n\leq p$, let $\{A_{l}\}$ be a partition of $\{1,\ldots p\}$ with the property that no two integers in $\{1,2,\ldots,n\}$ are  in the same  set $A_{l}$. That is, the partition $\{A_{l}\}$ separates $\{1,2,\ldots,n\}$.
 We   denote such a partition by $(\cup A_{l}=[1,p])\perp \{1,2,\ldots,n\}$.  We use the notation $\sum_{ \cup A_{l}=[1,p])\perp \{1,2,\ldots,n\}}$ to indicate the sum over all such partitions.

\begin{lemma}\label{lem-7.1}   Let $\mathcal{L}_{\al}$ be  a Poisson point process  on $\Om_{\De}$ with intensity measure $\al \mu$. Let  $\{n_{j}\}_{j=1}^{N}$ be a set of integers and define $s_{m} =\sum_{j=1}^{m-1}n_{j}$, $m=1,\ldots, N+1$.    Then for any functions  $g_{j}\in \DD_{1} (\mu)$, $j=1,\ldots, s_{N+1} $, 
 \bea
&&
E_{\mathcal{L}_{\al}}\(\prod_{m=1}^{N}   \phi_{n_{m}}  \(   g_{s_{m}+1},\ldots, g_{s_{m}+n_{m}} \)   \)\nonumber\\
&&\qquad  =
\sum_{\stackrel{\cup_{l}A_{l}=[1,s_{N+1}]}{\perp \{s_{m}+1,\ldots, s_{m}+n_{m}\},\,m=1,\ldots,N }}\prod_{l} \al \mu \(  g_{A_{l}} \).
\label{11.2}
\eea 
  where the    notation under the summation sign indicates that no sets in the partitions   can    contain more than one element from each of the sets $\{s_{m}+1,\ldots, s_{m}+n_{m}\}$, m=1,\ldots,N.
 \end{lemma}

\Proof  As in (\ref{ls.11}), by  the master formula for Poisson processes, \cite[(3.6)]{K},   for any finite  set of positive integers  $B$ and random variables $h_{j}$,  \bea
&&
E_{\mathcal{L}_{\al}}\(e^{ \sum_{j\in B}z_{j} \sum_{\om\in\LL_{\al}} h_{j}(\om)}\)\label{10.21w}\\
&&\qquad=\exp \(\al\(\int_{\Om_{\De}}\( e^{\sum_{j\in B} z_{j}h_{j}(\om) }
 -1\)\,d\mu(\om)  \)\).\nn
\eea
Differentiating  each side of (\ref{10.21w}) with respect to   each $z_{j}$, $j\in B$, and then setting all  the $z_{j}$ equal to zero, we   obtain
\be
E_{\mathcal{L}_{\al}}\(\prod_{j\in B}\( \sum_{\om\in\LL_{\al}}h_{j}(\om) \)\)=\sum_{\cup A_{l}=B }\prod_{l} \al \mu\(\prod_{j\in A_{l}} h_{j} \),\label{10.22j} 
\ee
  in which we sum over all partitions of $B$. We can write this as 
  \be
E_{\mathcal{L}_{\al}}\(\prod_{j\in B}\phi_{1}(h_{j})\)=\sum_{\cup A_{l}=B }\prod_{l} \al \mu\(  h_{A_{l}} \).\label{10.22r} 
\ee

The lemma follows from   (\ref{10.22r}). To see this we first consider first the case   $N=1$  and show that 
\begin{equation}
E_{\mathcal{L}_{\al}}\(  \phi_{n}(g_{1},\ldots, g_{n}) \)=\prod_{j=1}^{n}\al\mu(g_{j}).\label{N=1}
\end{equation}
This equation is stated in \cite[(3.14)]{K}. We present a detailed proof. 
When $n=1$, (\ref{N=1}) follows immediately from (\ref{10.21w}) and is well know.     It is  obviously the same as   (\ref{10.22r}). The proof for general $n$ proceeds by induction.
By  (\ref{10.22r})
\be
E_{\mathcal{L}_{\al}}\(\prod_{j=1}^{n}  \phi_{1}(g_{j}) \)\label{10.20c} = \sum_{k=1}^{n}\,\, \sum_{\cup_{l=1}^{k}A_{l}=[1,n]}\prod_{l=1}^{k} \al \mu \(  g_{A_{l}} \). 
\ee
Using (\ref{11.0b}) and   (\ref{10.20c}) we have
\begin{eqnarray}
E_{\mathcal{L}_{\al}}\(  \phi_{n}(g_{1},\ldots, g_{n}) \)&=&E_{\mathcal{L}_{\al}}\(\prod_{j=1}^{n}  \phi_{1}(g_{j}) \)
\label{7.10}\\
&&\qquad - \sum_{k=1}^{n-1}\,\, \sum_{\cup_{l=1}^{k}A_{l}=[1,n]}\,\,E_{\mathcal{L}_{\al}}\(\phi_{k}(g_{A_{1}},\ldots, g_{A_{k}}) \) \nonumber\\
&=& \sum_{k=1}^{n}\,\, \sum_{\cup_{l=1}^{k}A_{l}=[1,n]}\prod_{l=1}^{k} \al \mu \(  g_{A_{l}} \)
\nn\\
&&\qquad - \sum_{k=1}^{n-1}\,\, \sum_{\cup_{l=1}^{k}A_{l}=[1,n]}\,\,E_{\mathcal{L}_{\al}}\(\phi_{k}(g_{A_{1}},\ldots, g_{A_{k}}) \) .\nonumber
\end{eqnarray}
  Assuming that (\ref{N=1}) holds for $k\leq n-1$ we get
\begin{eqnarray}
E_{\mathcal{L}_{\al}}\(  \phi_{n}(g_{1},\ldots, g_{n}) \)
&=& \sum_{k=1}^{n}\,\, \sum_{\cup_{l=1}^{k}A_{l}=[1,n]}\prod_{l=1}^{k} \al \mu \(  g_{A_{l}} \)
\label{7.11}\\
&&\qquad - \sum_{k=1}^{n-1}\,\, \sum_{\cup_{l=1}^{k}A_{l}=[1,n]} \prod_{l=1}^{k} \al \mu \( g_{A_{l}} \). \nonumber
\end{eqnarray}
Since the only partition of $[1,n]$ with $n$ parts is given by $A_{j}=\{j\}, j=1,\ldots,n$, we obtain (\ref{N=1}) for $n$.

  We show below that 
 \begin{eqnarray}
&&
 \prod_{m=1}^{N}   \phi_{n_{m}}  \(   g_{s_{m}+1},\ldots, g_{s_{m}+n_{m}} \) \label{11.0bn} \\ 
&&\qquad = \sum_{k }\sum_{\stackrel{\cup_{l=1}^{k}A_{l}=[1,p]}{\perp \{s_{m}+1,\ldots, s_{m}+n_{m}\},\,m=1,\ldots,N }}\,\,\phi_{k}(g_{A_{1}},\ldots, g_{A_{k}}).\nonumber
\end{eqnarray}
Taking the expectation of (\ref{11.0bn}) and using (\ref{N=1}), we get (\ref{11.2}). This follows because
\begin{equation}
    \sum_{k }\sum_{\stackrel{\cup_{l=1}^{k}A_{l}=[1,p]}{\perp \{s_{m}+1,\ldots, s_{m}+n_{m}\},\,m=1,\ldots,N }}\prod_{j=1}^{k}\al\mu(g_{A_{j}})
   \end{equation}
is the same as the right-hand side of  (\ref{11.2}).

 To obtain (\ref{11.0bn}) we note that by (\ref{11.0})
 \begin{eqnarray}
\lefteqn{
 \prod_{m=1}^{N}   \phi_{n_{m}}  \(   g_{s_{m}+1},\ldots, g_{s_{m}+n_{m}} \) \label{11.0bnq}}\\ 
&&= \prod_{m=1}^{N}  \sum_{ (\om_{i_{s_{m}+1} },\ldots,\om_{i_{s_{m}+n_{m}} })\in S_{n_{m}}(\LL_{\al}) }\,\,\prod_{l=1}^{n_{m}}g_{s_{m}+l }(\om_{s_{m}+l })   \nn.
\eea
We rearrange this product of sums,   taking into account the points that occur in $S_{n_{m}}$ for more than one   $m$, to see that the second line in (\ref{11.0bnq}) is the same as
\be
\sum_{k }\sum_{\stackrel{\cup_{l=1}^{k}A_{l}=[1,p]}{\perp \{s_{m}+1,\ldots, s_{m}+n_{m}\},\,m=1,\ldots,N }}\,\, \sum_{ (\om_{i_{1} },\ldots,\om_{i_{k} })\in S_{k}(\LL_{\al}) }\,\,\prod_{l=1}^{k}g_{A_{l}}(\om_{i_{l} })
\ee
By the definition (\ref{11.0}) this is the same as the second line of (\ref{11.0bn}). \qed

The relationship in 
 (\ref{11.0bn})  
is   a  Poissonian   analogue of the Wick expansion formula in \cite[Lemma 2.3]{MRmem}. We restate it as a lemma.

\begin{lemma}\label{lem-11.5} Let $\phi_{n}$ be as defined in (\ref{11.0}),   $s_{m} =\sum_{j=1}^{m-1}n_{j}$, $m=1,\ldots, N+1$ and   $p=   s_{N+1} $.   Then for any functions  $g_{j} $, $j=1,\ldots, p$,   
\bea
&&
\prod_{m=1}^{N}   \phi_{n_{m}}  \(   g_{s_{m}+1},\ldots, g_{s_{m}+n_{m}} \)\label{11.4a}\\
&& \qquad
=\sum_{k }  \sum_{\stackrel{\cup_{l=1}^{k}A_{l}=[1,p]}{\perp \{s_{m}+1,\ldots, s_{m}+n_{m}\},\,m=1,\ldots,N }}\,\,\phi_{k}(g_{A_{1}},\ldots, g_{A_{k}}). \nn
\eea
 \end{lemma}
 
 \begin{example} {\rm 
  When all $n_{m}=1$, $m=1,\ldots,N$, in (\ref{11.4a}) it gives the equation in  (\ref{11.0b}).
Here are  some other examples,  
\be
\phi_{2}(g_{1},g_{2}) \phi_{1}(g_{3})=\phi_{3}(g_{1},g_{2},g_{3})+\phi_{2}(g_{[1,3]},g_{2})
\label{11.4b}+\phi_{2}(g_{1},g_{[2,3]}),
\ee
\begin{eqnarray}
&&\hspace{-.3 in}\phi_{2}(g_{1},g_{2})\,\phi_{2}(g_{3},g_{4})=\phi_{4}(g_{1},g_{2},g_{3},g_{4})
\label{11.5}\\
&&+ \phi_{3}(g_{[1,3]},g_{2}, g_{4})+\phi_{3}(g_{[1,4]},g_{2} ,g_{3})+\phi_{3}(g_{1},g_{[2,3]},g_{4}) + \phi_{3}(g_{1},g_{3},g_{[2,4]})
 \nonumber\\
&&\hspace{1.2 in}+ \phi_{2}(g_{[1,3]},g_{[2,4]})+\phi_{2}(g_{[1,4]},g_{[2,3]}).
 \nonumber
\end{eqnarray}
In addition when all $n_{m}=1$, and $g_{m}=g$, $m=1,\ldots,n$ \begin{eqnarray}
  \phi^{n}_{1}(g) \label{11.5a} &= &\sum_{k=1}^{n}\,\, \sum_{\cup_{l=1}^{k}A_{l}=[1,n]}\,\,\phi_{k}(g_{A_{1}},\ldots, g_{A_{k}}) \\
 &=&\sum_{k=1}^{n}{1 \over k!}\sum_{ n_{1}+ \cdots +n_{k}=n}{n \choose n_{1}\hspace{.1 in} \cdots \hspace{.1 in}n_{k} }\phi_{k}(g_{[1,n_{1}]},\ldots, g_{[1,n_{k}]}).\nonumber
\end{eqnarray} 
 This is a special case of (\ref{11.0b}). 

 }\end{example}

 The random variables $\{\phi_{n}(g_{1},\ldots, g_{n})\}_{n=1}^{\ff}$   are not necessarily orthogonal with respect to $E_{\LL_{\al}}$. 
Let 
\bea
&&
I_{n }(g_{1},\ldots, g_{n}) \label{11.21}\\
&&\qquad:=\sum_{D=\{i_{1},\ldots,i_{|D|}\}\subseteq [1,n]}(-1)^{|D^{c}|}
\phi_{|D|}(g_{i_{1}},\ldots,g_{i_{|D|}} ) \,\prod_{j\in D^{c}}\al\mu(g_{j }).\nn
	 \eea
	 We refer to the random variables $I_{n }(g_{1},\ldots, g_{n})$ as Poisson Wick\label{p56} products and  show in Corollary \ref{cor-7.1q} that they are 	 orthogonal with respect to $E_{\LL_{\al}}$.
 
Considering (\ref{11.21}) we note that  $I_{n }(g_{1},\ldots, g_{n})-\phi_{n}(g_{1},\ldots, g_{n})\in \HH_{n-1}$.

  \medskip	 Let $p\ge k$ and let $\{A_{l}\}$ be a partition of $[1,p]$. Let  $ \{\cup A_{l}=[1,p], \,[1,k]_{ns}  \}$ be the subset of partitions of $[1,p]$
  in which all partitions that contain  a set consisting of only one member of $[1,k]$ are eliminated. (The symbol $ns$ indicates no singleton.)

\begin{lemma}\label{lem-7.3}   Let $\mathcal{L}_{\al}$ be  a Poisson point process  on $\Om_{\De}$ with intensity measure $\al \mu$. Let $s_{m} =\sum_{j=1}^{m-1}n_{j}$, $m=1,\ldots, N+1$ and let   $p\ge    s_{N+1} $.  Then for any functions  $g_{j}\in \DD_{1} (\mu)$, $j=1,\ldots, p$, 
 \bea
&&
E_{\mathcal{L}_{\al}}\(\prod_{m=1}^{N} I_{n_{m} }(g_{s_{m}+1},\ldots, g_{s_{m}+n_{m}}) \prod_{j=s_{N+1}+1}^{p} \,\,\phi_{1}( g_{j } )  \) \label{11.23}\\
&&\qquad  =
\sum_{\stackrel {\cup_{l}A_{l}=[1,p],\,  [1,s_{N+1}]_{ns} }{\perp \{s_{m}+1,\ldots, s_{m}+n_{m}\},\,m=1,\ldots,N }}\prod_{l} \al \mu \(  g_{ A_{l} } \).
\nn
\eea 
 \end{lemma}
 
  (The following proof   is direct, but involves a good deal of combinatorics. In Section \ref{sec-expP} we give an alternate proof using   Poissonian exponentials.)

 \medskip	 \Proof 
   Consider first the case in which $N=1$. Using (\ref{11.21}), with $n_{1}=n$, we see that the left-hand side of (\ref{11.23}) is equal to    
     \be 
  \sum_{D=\{i_{1},\ldots,i_{|D|}\}\subseteq [1,n]}(-1)^{|D^{c}|}E_{\mathcal{L}_{\al}}\( \phi_{|D|}(g_{i_{1}},\ldots,g_{i_{|D|}} )\prod_{j=n+1}^{p} \,\,\phi_{1}( g_{j } )\)  
 \prod_{j\in D^{c}}\al\mu(g_{j}).  
 \ee 
 By (\ref{11.2}) this is equal to 
 \be
\sum_{D=\{i_{1},\ldots,i_{|D|}\}\subseteq [1,n]}(-1)^{|D^{c}|} \sum_{(\cup A_{l}= D\cup [n+1,p])\perp D }\,\prod_{l} \al \mu \(  g_{A_{l} } \)\prod_{j\in D^{c}}\al\mu(g_{j }),
 \label{alt.5}
 \ee
 which we  write as
  \be
\sum_{D\subseteq [1,n]}(-1)^{|D^{c}|} \sum_{\stackrel{(\cup A_{l}=[1,p])\perp \{1,2,\ldots,n\}}{ (D^{c})_{s}}}\,\prod_{l} \al \mu \(  g_{A_{l}} \),
 \label{alt.6}
 \ee
where $(D^{c})_{s}$ indicates that each element in $D^{c}$ is a singleton of the partition $\cup A_{l}=[1,p]$.    (This is simple, we just take a partition in $(\cup A_{l}= D\cup [n+1,p])\perp D $ and add the singletons in $D^{c}$.)  
 
 Consider a partition with singletons.  Let $C\subseteq [1,n]$ denote the   singletons in such a partition.  The terms in (\ref{alt.6}) which give rise to such a partition are precisely those with $D^{c}\subseteq C$.   Therefore, if $|C|=m$,
\begin{equation}
\sum_{D^{c}\subseteq C}(-1)^{|D^{c}|}=\sum_{k=0}^{m}{m \choose k}(-1)^{k}=0.\label{10.22s}
\end{equation}
  This shows that when $N=1$, (\ref{alt.6}) equals the right hand side of  (\ref{11.23}).  

   The proof for general $N$ consists of a straight forward, albeit tedious,  generalization of the  arguments used for $N=1$. 
  Using (\ref{11.21}) and expanding the left-hand side of (\ref{11.23})  we get that it is equal to  
\begin{eqnarray}
 &&  \sum_{\stackrel {D_{m}=\{i_{m,1},\ldots,i_{m,|D_{m}|}\}\subseteq [s_{m}+1, s_{m}+n_{m}]}   {m=1,\ldots,N}  }(-1)^{|\cup D_{m}^{c}|} \prod_{j\in \cup D_{m}^{c}}\al\mu(g_{j})  \label{alt.4n}\\
 &&\hspace{1in} E_{\mathcal{L}_{\al}}\( \prod_{m=1}^{N}\phi_{|D_{m}|}(g_{m,1},\ldots,g_{m,|D_{m}|} ) \prod_{j=s_{N+1}+1}^{p} \,\,\phi_{1}( g_{j } )\).\nn
 \end{eqnarray}
 By (\ref{11.2}) this is 
 \bea
 &&
\sum_{\stackrel {D_{m}\subseteq [s_{m}+1, s_{m}+n_{m}]}   {m=1,\ldots,N}  }(-1)^{|\cup D_{m}^{c}|}  \prod_{j\in \cup D_{m}^{c}}\al\mu(g_{j}) \label{alt.5n}\\
&& \hspace{1in}\sum_{ (\cup A_{l}=\cup D_{m}\cup [s_{N+1}+1,p])\perp   D_{j},\,j=1,\ldots,N }\,\prod_{l} \al \mu \( g_{A_{l}} \).
\nn
 \eea
  As in (\ref{alt.6}) this is
  \be
\sum_{\stackrel {D_{m}\subseteq [s_{m}+1, s_{m}+n_{m}]}   {m=1,\ldots,N}  }(-1)^{|\cup D_{m}^{c}|} \!\sum_{\stackrel{(\cup A_{l}=[1,p])\perp \{s_{m}+1,\ldots, s_{m}+n_{m}\}, m=1,\ldots, N}{ (\cup D_{m}^{c})_{s}}}\,\prod_{l} \al \mu \(  g_{A_{l}} \).
 \label{alt.6n}
 \ee
 
  Consider a partition with singletons, and let  $C$ denote the collection of singletons in such a partition.  We  write $C=\cup C_{m}$ where $C_{m}\subseteq [s_{m}+1, s_{m}+n_{m}]$. The terms in (\ref{alt.6n}) which give rise to such a partition are precisely those with $D_{m}^{c}\subseteq C_{m}$  for all $m$. Note then  if $|C_{m}|=s_{m}$,
\begin{equation}
\sum_{\stackrel{D_{m}^{c}\subseteq C_{m}}{m=1,\ldots, N}}(-1)^{|\cup D_{m}^{c}|}= \prod_{m=1}^{N}\(\sum_{D_{m}^{c}\subseteq C_{m}}(-1)^{| D_{m}^{c}|}\)= \prod_{m=1}^{N}\(\sum_{k=0}^{s_{m}}{s_{m} \choose k}(-1)^{k}\)=0.\label{10.22sn}
\end{equation}
This shows that (\ref{alt.6n}) equals the right hand side of (\ref{11.23}).\qed

An important consequence of Lemma \ref{lem-7.3} is that $I_{n}(g_{1},\ldots, g_{n})$ extends to a continuous multilinear map from $\(\DD_{2}(\mu)\)^{n}$ to $\DD_{1}(P_{\mathcal{L}_{\al}})$, and that (\ref{11.23}) holds for   the extension.     This observation is critical  in the proofs in Section \ref{sec-lspower}.   We state it as a lemma.

\begin{lemma} 
 \label{cor-11.14}   Let $\mathcal{L}_{\al}$ be  a Poisson point process  on $\Om_{\De}$ with intensity measure $\al \mu$. Let $s_{m} =\sum_{j=1}^{m-1}n_{j}$, $m=1,\ldots, N+1$.  Then for any functions  $g_{j}\in \DD_{ 2} (\mu)$, $j=1,\ldots, s_{N+1}$,  
 \bea
&&
E_{\mathcal{L}_{\al}}\(\prod_{m=1}^{N} I_{n_{m} }(g_{s_{m}+1},\ldots, g_{s_{m}+n_{m}})  \) \label{11.23a}\\
&&\qquad=
\sum_{\stackrel {\cup_{l}A_{l}=[1,s_{N+1}],\,  |A_{l}|\geq 2 }{\perp \{s_{m}+1,\ldots, s_{m}+n_{m}\},\,m=1,\ldots,N }}\prod_{l} \al \mu \(    g_{ A_{l} } \).
\nn
\eea 
 \el

 \Proof    Let $\{D_{i}\}$ be as defined in (\ref{dm}).
 If $g\in \DD_{ 2} (\mu)$ then $1_{\{D_{i}\}}g \in \DD_{ 1} (\mu)$ and $1_{\{D_{i}\}}g \to g$ as $i\rar\ff$ in all $L^{p} (\mu)$, $p\geq 2$. By Lemma \ref{lem-7.3} we see that (\ref{11.23a})
 holds with $g_{j}$ replaced by $1_{\{D_{i}\}}g_{j}$. Because of the condition that all    the sets $A_{l}$ have  $|A_{l}|\geq 2$,
 \be  \lim_{i\rar\ff}\prod_{l} \al \mu \(\prod_{j\in A_{l}} 1_{\{D_{i}\}}g_{j } \) =\prod_{l} \al \mu \(\prod_{j\in A_{l}} g_{j } \).
 \ee
In this way we can  extend the left-hand side of  (\ref{11.23a}) to random variables   $g_{j}\in \DD_{2}(\mu)$, and   get that (\ref{11.23a}) holds for the extension. 
  \qed

 \begin{corollary} \label{cor-7.1q}The random variables $\{I_{n}(g_{1}\ldots,g_{n} )\}_{n=0}^{\ff}$,   with  $I_{0}\equiv 1$   and  $g_{j}\in \DD_{ 2} (\mu)$, $j=1,\ldots, n$, are orthogonal in $L^{2}\(P_{\mathcal{L}_{\al}}\)$. In particular,   for $f_{i}\in \DD_{ 2} (\mu)$, $i=1,\ldots, m $,
 \begin{equation}
E_{\mathcal{L}_{\al}}\(  I_{n  }(g_{ 1},\ldots, g_{ n} )  I_{m  }(f_{  1},\ldots, f_{ m} )   \)=\de_{n,m} \sum_{\pi\in  \mathcal{P}_{n}}\prod_{j=1}^{n} \al \mu \(  g_{j }f_{\pi (j)} \).\label{11.23p}
 \end{equation}
\end{corollary}
 
 \Proof By (\ref{11.23}) 
  \bea
&&
E_{\mathcal{L}_{\al}}\(  I_{n  }(g_{ 1},\ldots, g_{ n} )  I_{m  }(g_{n+ 1},\ldots, g_{n+ m} )   \) \label{11.23q}\\
&&\qquad=
\sum_{\stackrel{\{\cup_{l}A_{l}=[1,n+m],\,  [1, n+m]_{ns}\}}{\perp \{  1,\ldots,  n \}, \{n+1,\ldots,  n +m \}  }}\prod_{l} \al \mu \(\prod_{j\in A_{l}} g_{j } \).
\nn
\eea 
The restrictions on the partitions that are summed over require that all the sets $\{A_{l}\}$ contain two elements, one from $[1,n]$ and one from $[n+1,n+m]$. Therefore we must have $m=n$ to get any contribution from the right-hand side of (\ref{11.23q}).   The statement in (\ref{11.23p}) follows easily from this observation.\qed

In the paragraph containing (\ref{11.0c}) we define  $\HH_{n}$ to be the closure in $L^{2}\(P_{\mathcal{L}_{\al}}\)$ of all linear combinations of $\phi_{j}$, $j=0,1,\ldots, n$. By Corollary \ref{cor-7.1q} and the remarks in the paragraph containing (\ref{11.0c}) we see that
$I_{n}\in \HH_{n}  \ominus \HH_{n-1}=\HH_{n}\cap \HH^{\perp}_{n-1}$  and that linear combinations of the $\{I_{n}\}$  are dense in  $H_{n}:=\HH_{n}\ominus \HH_{n-1}$.
Using (\ref{11.0c}) we then have the orthogonal decomposition
\begin{equation}
L^{2}\(P_{\mathcal{L}_{\al}}\)=   \oplus _{n=0}^{\ff }\,\, H_{n}.  \label{11.0r}
\end{equation}
We call this the   Poisson chaos decomposition   of $L^{2}\(P_{\mathcal{L}_{\al}}\)$ .

\begin{remark}\label{rem-75} {\rm  Using the   Poisson chaos decomposition   of $L^{2}\(P_{\mathcal{L}_{\al}}\)$ we obtain the following analogue of Lemma \ref{lem-11.5}:

\begin{lemma}\label{lem-11.35} Let $I_{n}$ be as defined in (\ref{11.21}),   $s_{m} =\sum_{j=1}^{m-1}n_{j}$, $m=1,\ldots, N+1$.   Then for any functions  $g_{j} $, $j=1,\ldots,  s_{N+1}$,
\bea
&&
\prod_{m=1}^{N} I_{n_{m} }(g_{s_{m}+1},\ldots, g_{s_{m}+n_{m}}) 
\label{11.34}\\
&&\qquad=\sum_{j,k} \sum_{\stackrel{\cup_{i=1}^{j}B_{i}\,\,\cup_{l=1}^{k}C_{l}=[1,  s_{N+1}],\,|C_{l}|\geq 2 }{\perp \{s_{m}+1,\ldots, s_{m}+n_{m}\},\,m=1,\ldots,N }}I_{j }(g_{B_{1}},\ldots, g_{B_{j}})\,\prod_{l=1}^{k}\al\mu (g_{C_{l}}). \nn
\eea
 \end{lemma}
 
 \Proof   Considering the observations in   the paragraph preceding this remark,  it suffices to show that both sides of (\ref{11.34}) have the same inner product in $L^{2}\(P_{\mathcal{L}_{\al}}\)$ with all functions of the form
  \be J:=I_{n_{N+1} }(g_{s_{N+1}+1},\ldots, g_{s_{N+1}+n_{N+1}}). 
  \ee
  We first compute this for the left-hand side of (\ref{11.34}). It follows from    Lemma \ref{cor-11.14} that
  \bea
&&
E_{\mathcal{L}_{\al}}\(\prod_{m=1}^{N} I_{n_{m} }(g_{s_{m}+1},\ldots, g_{s_{m}+n_{m}})\,\, J \) \label{11.33a}\\
&&
\qquad=E_{\mathcal{L}_{\al}}\(\prod_{m=1}^{N+1} I_{n_{m} }(g_{s_{m}+1},\ldots, g_{s_{m}+n_{m}})  \) \nn\\
&&\qquad=
\sum_{\stackrel {\cup_{l}A_{l}=[1,s_{N+2}],\,  |A_{l}|\geq 2 }{\perp \{s_{m}+1,\ldots, s_{m}+n_{m}\},\,m=1,\ldots,N+1 }}\prod_{l} \al \mu \(    g_{ A_{l} } \).
\nn
\eea 

We divide the sets $A_{l}$ into two groups. We denote by $C_{1},\ldots, C_{k}$ the collection of sets $A_{l}$ that are disjoint from $[s_{N+1}+1, s_{N+1}+n_{N+1}]$, and by $B'_{1},\ldots, B'_{j}$ the collection of sets $A_{l}$ that contain an element from $[s_{N+1}+1, s_{N+1}+n_{N+1}]$.  Note that the restriction indicated by the notation  $\perp \{s_{N+1}+1,\ldots, s_{N+1}+n_{N+1}\}$ implies that there is  exactly one element  $\rho(i)\in [s_{N+1}+1, s_{N+1}+n_{N+1}]$ in each $B'_{i}$.  Consequently   $j$, the number of sets in $B'_{1},\ldots, B'_{j}$ , is $n_{N+1}$. In addition, the condition $|A_{l}|\geq 2$ implies that  each $B'_{i}$ contains   at least one element which is not in $[s_{N+1}+1, s_{N+1}+n_{N+1}]$.   Set  $B_{i}=B'_{i}-\rho(i)$. Clearly $B_{i}\subseteq [1,  s_{N+1}]$. 

Considering the remarks in the previous paragraph we see that   (\ref{11.33a}) can be rewritten as
  \bea
\lefteqn{
E_{\mathcal{L}_{\al}}\(\prod_{m=1}^{N} I_{n_{m} }(g_{s_{m}+1},\ldots, g_{s_{m}+n_{m}})\,\, J \) \label{11.33b}}\\
&&=
\sum_{k} \sum_{\stackrel{\cup_{i=1}^{n_{N+1}}B_{i}\cup_{l=1}^{k}C_{l}=[1,  s_{N+1}],\,|C_{l}|\geq 2 }{\perp \{s_{m}+1,\ldots, s_{m}+n_{m}\},\,m=1,\ldots,N }}\,\,\sum_{\rho\in  \mathcal{M}_{n_{N+1}}}\prod_{i=1}^{n_{N+1}} \al \mu \( g_{B_{i}} \,\, g_{\rho(i)}\)\,\prod_{l=1}^{k}\al\mu (g_{C_{l}}),
\nn
\eea 
where
  $   \mathcal{M}_{n_{N+1}} $ is the set of bijections from $ [ 1, n_{N+1}]$ to $ [s_{N+1}+1, s_{N+1}+n_{N+1}]$.  

  Using Corollary \ref{cor-7.1q}
we see that we obtain the same value for the
 inner product, in $L^{2}\(P_{\mathcal{L}_{\al}}\)$,  of the right-hand side of (\ref{11.34}) with $J$. \qed }
\end{remark}

\section{Loop soup decomposition of  permanental Wick powers}\label{sec-lspower}

The $n$-th order permanental Wick power,
  $\wt\psi_{n}(\nu)$ is defined in Section \ref{sec-powers} as a renormalized sum of powers of a permanental field,  $\psi(\nu)$.  In this section we define $\wt\psi_{n}(\nu)$ using loop soups. This enables us to easily obtain   an  isomorphism theorem involving $\wt\psi_{n}(\nu)$ and $n$-fold intersection local times $L_{n}(\nu)$. 
  
  The loop soup    description of  $\wt\psi_{n}(\nu)$ involves the permanental chaos   $\psi_{n}(\nu)$,  which we think of as involving self-intersections of each path in the loop soup, and  new random variables, $\II_{l_{1},\ldots, l_{j}}(\nu)$, which we think  of as involving intersections of different paths.   The representation of $\wt\psi_{n}(\nu)$ is actually its Poisson chaos decomposition into  random variables  $\II_{l_{1},\ldots, l_{j}}(\nu)$, for different indices $l_{1},\ldots, l_{j}$.

 We proceed to  develop the material needed to define $\II_{l_{1},\ldots, l_{j}}(\nu)$.  
 Note that for fixed $r>0$, $L_{j}(x,r)$ is a polynomial in $L(x,r)$ without a constant term. Therefore by (\ref{ls.4}),  with $g_{j}=f_{x_{j},r}$ we see that    
 $L_{j}(x,r)\in \DD_{2}(\mu)$  for each $j\geq 1$.  It follows from this and Corollary \ref{cor-11.14} that  
 \begin{equation}
 \II_{l_{1},\ldots, l_{n}}(x,r)  :=  I_{n}(L_{l_{1}}(x,r),\ldots,L_{l_{n}}(x,r))\label{defil}
   \end{equation}
   is well defined.

   \begin{theorem} \label{theo-7.1}
Let $l=l_{1}+\cdots+l_{n}$ and $ \nu\in \BB_{2l}(R^{d})$. Then
\begin{equation}
\II_{l_{1},\ldots, l_{n}}(\nu):=\lim_{r\rar 0}\int I_{n}(L_{l_{1}}(x,r),\ldots,L_{l_{n}}(x,r))\,d\nu (x)\label{10.24}
\end{equation}
exists in all $L^{p}(P_{  \LL_{ \al}})$.
 \end{theorem}
 
 Note that $\II_{l_{1},\ldots, l_{n}}(\nu)\in H_{n}$.
 
 \medskip	\Proof
 It follows from (\ref{11.23a})  that for integers $v\ge 1$, 
\bea
&&
E_{\mathcal{L}_{\al}}\(\prod_{i=1}^{ v}  \II_{l_{1},\ldots, l_{n}}(x_{i},r_{i})  \)\label{10.23sq}\\
&&\qquad=
\sum_{\stackrel{\cup_{l}A_{l}=[1,nv], |A_{l}|\geq 2}{\perp \{(m-1)n+1,\ldots, mn\},\,m=1,\ldots,v}}\prod_{l} \al \mu \(\prod_{j\in A_{l}} g_{j} \),
\nn
\eea 
in which each  $g_{j}$ is one of the random variables  $L_{l_{k}}(x_{i},r_{i})$.     The lemma follows from  (\ref{2.27w}) and the work that follows it in Section \ref{sec-iltac}.  \qed

It should be clear from the definitions   that $\II_{l_{1},\ldots, l_{n}}(\nu)=\II_{l_{\si(1)},\ldots, l_{\si (n)}}(\nu)$ for any permutation $\si$ of $[1,n]$.

\medskip	 Let ${\ell}=(l_{1},\ldots, l_{k})$ and let  $m_{i}(\ell)$ be the number of indices $j$ with  $l_{j}=i$.  We define the degree of $\ell$ to be 
\begin{equation}
d(\ell )=(m_{1}(\ell ),m_{2}(\ell ),\ldots, m_{k}(\ell)),\label{adj.11}
\end{equation}
and set
\begin{equation}
e(\ell)=\prod_{i}^{k}m_{i}(\ell_{k} )!\label{adj.12}
\end{equation}

The next lemma shows that the processes $\II_{l_{1},\ldots, l_{k}}(\nu) $   with different  degrees $d(\ell)$ are orthogonal.

\bl\label{lem-7.11}
  Let ${\ell}=(l_{1},\ldots, l_{k})$ and ${\ell'}=(l'_{1},\ldots, l'_{k'})$. Then
\begin{eqnarray}
&&E_{\mathcal{L}_{\al}}\( \II_{l_{1},\ldots, l_{k}}(\nu)  \,\,  \II_{l'_{1},\ldots, l'_{k'}}(\nu) \) 
\label{adj.13}\\
&&\qquad = \de_{\{d(\ell),d(\ell') \}  }\,\,e(\ell)\,\,\lim_{r\rar 0}\int \prod_{j=1}^{k}\al \mu\(L_{l_{j}}(x,r)L_{l_{j}}(y,r)\) \,d\nu  (x )\,d\nu  (y). \nonumber
\end{eqnarray}
\el
\Proof This proof is similar to the proof of Lemma \ref{lem-6.1}. By (\ref{11.23}),
\bea
\lefteqn{
E_{\mathcal{L}_{\al}}\( \II_{l_{1},\ldots, l_{k}}(\nu)  \,\,  \II_{l'_{1},\ldots, l'_{k'}}(\nu) \) \label{adj.15}}\\
&&=
\sum_{\stackrel{\cup_{m}A_{m}=[1,k+k'], |A_{m}|\geq 2}{\perp [1,k], [k+1,k+k']}}\,\,\lim_{r\rar 0}\int\prod_{m} \al \mu \(\prod_{j\in A_{m}} \wt L_{a_{j}}(\wt x,r) \) \,d\nu  (x )\,d\nu  (y),
\nn
\eea
where $\wt L_{a_{j}}(\wt x,r)= L_{l_{j}}(  x,r)$ if $j\in [1,k]$ and   $\wt L_{a_{j}}(\wt x,r)= L_{l'_{j-k}}(  y,r)$ if $j\in [k+1,k+k']$. By the separation condition, $\perp [1,k], [k+1,k+k']$, we see that we must have $|A_{m}|=2$ for each $m$, and each $A_{m}$ consists of one  element $j_{m}$ from $[1,k]$
and one   element $j'_{m}$ from    $[k+1,k+k']$. However, by the alternating condition in (\ref{2.35}), we see that in the limit $\mu \(L_{l_{j_{m}}}(x,r)L_{l'_{j'_{m}-k}}(y,r) \)$ would give a contribution of $0$ unless $l_{j_{m}}=l'_{j'_{m}-k}$. This forces $d(\ell )=d(\ell' )$ and (\ref{adj.13}) follows since there are 
$e(\ell )$ ways to pair the elements of $\ell $ with those of $\ell' $ (with identical integers in each pair).\qed

 	 \begin{lemma}  \label{lem-7.7}
 Let  $\,l_{m}=\sum_{j=1}^{n_{m}}l_{m(j)}$  and $l=\sum_{m=1}^{N}l_{m}$.  Furthermore, let $ \nu_{m}\in \BB_{2l_{m}}$, $m=1,\ldots ,N$.   Then  
  \bea
&&
E_{\mathcal{L}_{\al}}\(\prod_{m=1}^{N}\II_{l_{m(1)},\ldots, l_{m( n_{m})}}(\nu_{m})\, \psi^{p} (x,r)\)\label{10.25}\\
&&\qquad= \sum_{\pi\in \mathcal{P}_{l+p, a^{-}}}\al^{c(\pi)}  \int  \prod_{j=1}^{l+p}u (z_{j},z_{\pi(j)})
\prod_{m=1}^{N}\,d\nu_{m}  (x_{m} )\prod_{i=l+1}^{l+p}f_{x,r}(x_{i} ) \,dx_{i}   ,\nn\nn
\eea
 where  $z_{1},\ldots,z_{l_{1}}$ are all equal to $x_{1}$,  the next $l_{2}$ of the $\{z_{j}\}$ are all equal to $x_{2}$, and so on, so that the last   $l_{N}$ of the $\{z_{j}\}$ are all equal to $x_{N}$.
 
In addition $\mathcal{P}_{l+p,a^{-}}$ is the set of permutations $\pi$ of $[1,l+p]$ with cycles that alternate the variables $\{x_{i}\}$; (i.e., for all $j$, if $z_{j}=x_{i}$ then $z_{\pi(j)}\neq x_{i}$), 
and, for each $m=1,\ldots,N$, $\pi$ contains distinct cycles $C_{m(1)},\ldots, C_{m(n_{m})}$, such that $C_{m(1)}$
 contains the first $l_{m(1)}$ of the $x_{m}$'s, $C_{m(2)}$
 contains the next $l_{m(2)}$ of the $x_{m}$'s, etc. ( If $z_{a+1}=\cdots=z_{a+l_{m}}=x_{m}$, then the first  $l_{m(1)}$ of the $x_{m}$'s are $z_{a+1},\ldots,z_{a+l_{m(1)}}$,   the next $l_{m(2)}$ of the $x_{m}$'s are $z_{a+l_{m(1)}+1},\ldots,z_{a+l_{m(1)}+l_{m(2)}}$, etc.)

 \end{lemma}

 \Proof     
  Let $s_{m} =\sum_{j=1}^{m-1}n_{j}$, $m=1,\ldots, N+1$.  Then for any functions  $g_{j}\in \DD_{ 2} (\mu)$, $j=1,\ldots, s_{N+1}+p$ it follows from    (\ref{11.23a})  that
  \bea
&&
E_{\mathcal{L}_{\al}}\(\prod_{m=1}^{N} I_{n_{m}   }(g_{s_{m}+1},\ldots, g_{s_{m}+n_{m}})\prod_{j=1}^{p}I_{1}(g_{s_{N+1}+j})\)\label{10.23o}\\
&&\qquad=
\sum_{\stackrel{\cup_{l}A_{l}=[1,s_{N+1}+p], |A_{l}|\geq 2}{\perp \{s_{m}+1,\ldots, s_{m}+n_{m}\},\,m=1,\ldots,N }}\prod_{l} \al \mu \(\prod_{j\in A_{l}} g_{j } \).
\nn
\eea  
 
\noindent  Since $1_{D_{m}}L_{1}(x,r)\to L_{1}(x,r)$ in $\DD_{2}(\mu)$ it follows from (\ref{sxr}) that
\begin{equation}
\psi (x,r)=I_{1}(L_{1}(x,r)).\label{ls.10aay}
\end{equation}   
We use (\ref{10.23o})  with $g_{s_{m}+j}=L_{m(j) }(x_{m},r' )$, $1\le j\le n_{m}$ and $g_{j}=  L_{1}(x,r )$,    $s_{N+1}+1\le j \le p$ and integrate both sides  with respect to the measures $\nu_{m}(x_{m})$, $1\le m\le N$   and use Theorem \ref{theo-7.1}  to get 
 \bea
&&
E_{\mathcal{L}_{\al}}\(\prod_{m=1}^{N}\II_{l_{m(1)},\ldots, l_{m( n_{m})}}(\nu_{m})\, \psi^{p} (x,r)\)\label{10.25e}\\
&&\qquad= \sum_{\stackrel{\cup_{l}A_{l}=[1,s_{N+1}+p], |A_{l}|\geq 2}{\perp \{s_{m}+1,\ldots, s_{m}+n_{m}\},\,m=1,\ldots,N }}\prod_{l} \al \mu \(\prod_{j\in A_{l}} \wh L_{j } \).\nn
\eea
 where    $\wh L_{s_{m}+j}=L_{m(j) }(\nu_{m}  )$, $1\le j\le n_{m}$ and $\wh L_{j }=  L_{1  }(x,r )$,  $ l+1\le j \le l+ p$.

 The right-hand side of (\ref{10.25e}) is the same as the right-hand side of (\ref{10.25}).
 To see this we use  Theorem \ref{theo-1.3} to analyze each factor $\mu \(\prod_{j\in A_{l}} \wh L_{j } \)$ as in the proof of  Theorem \ref{theo-1.4} on page \pageref{forpageref1.4}. Each of these factors corresponds to a cycle in a permutation.  The  condition in the second line under the summation sign in (\ref{10.25e}) requires that in each   factor $\mu \(\prod_{j\in A_{l}} \wh L_{j } \)$, for each $m=1,\ldots ,N$,  there is at most one  term $L_{m(j) }(\nu_{m}  )$, $1\le j\le n_{m}$.   
 Therefore   distinct terms  $L_{m(j) }(\nu_{m}  )$, $1\le j\le n_{m}$,  are contained in   distinct cycles.   \qed

\bt\label{theo-multiloop} For $\nu\in\BB_{2n}(R^{d})$  
\be
\wt\psi_{n}(\nu)= \sum_{D_{1}\cup\cdots\cup D_{l}=[1,n]}
\II_{|D_{1}|,\ldots, |D_{l}|}(\nu),\label{10.26}
\ee
where the sum is over all partitions of $[1,n]$.  
\et

  \Proof  Set
\be
\phi_{n}(\nu)= \sum_{D_{1}\cup\cdots\cup D_{l}=[1,n]}
\II_{|D_{1}|,\ldots, |D_{l}|}(\nu).\label{10.26s}
\ee
It follows from (\ref{adj.13}) and   (\ref{adj.3q}) and then from  (\ref{20.14}) and (\ref{20.14a}),  that
\bea
 E_{\mathcal{L}_{\al}}\(    \phi^{2}_{n}(\nu )\) & =&E_{\mathcal{L}_{\al}}\(  \wt  \psi^{2}_{n}(\nu )\)\label{7.47}\\
 &=&\sum_{\pi\in \mathcal{P}_{2n,a}}\al^{c(\pi)} \,\int\!\!\int \(  u (x,y)u (y,x)  \)^{n}  \,d\nu  (x )\,d\nu  (y)\nn,
   \eea 
  where, as above, $\mathcal{P}_{2n,a}$ is the set of permutations of $[1,2n]$ which alternate $ 1,\ldots, n $ and $ n+1,\ldots, 2n $.
  
It follows from  (\ref{10.25}) with $N=1$ and $k=n$,  that  
 \bea 
 &&
E_{\mathcal{L}_{\al}}\(  \II_{|D_{1}|,\ldots, |D_{l}|}(\nu) \psi^{n} (x,r) \)\label{20.14xh}\\
&&\qquad= \sum_{\pi\in \mathcal{P}_{2n,a^{-}}}\al^{c(\pi)}\int     \prod_{j=1}^{2n}u (z_{j},z_{\pi(j)})\,d\nu  (x_{1} )\prod_{i=2}^{n+1}f_{x,r}(x_{i} ) \,dx_{i}\nn. \eea
Consider the restrictions on $\pi\in \mathcal{P}_{2n,a^{-}}$ for this partition, $|D_{1}|,\ldots, |D_{l}|$, given in Lemma \ref{lem-7.7}. Note that each different partition, say $|D_{1}|,\ldots, |D_{l'}|$ of $[1,n]$ give rise to a disjoint set of partitions and adding over all of them as in (\ref{10.26s}) we get 
 \bea 
 &&
E_{\mathcal{L}_{\al}}\(    \phi_{n}(\nu )  \psi^{n} (x,r) \)\label{20.14x}\\
&&\qquad= \sum_{\pi\in \mathcal{P}_{2n,a}}\al^{c(\pi)}\int     \prod_{j=1}^{2n}u (z_{j},z_{\pi(j)})\,d\nu  (x_{1} )\prod_{i=2}^{n+1}f_{x,r}(x_{i} ) \,dx_{i}.\nn
\eea

 We write the left-hand side of (\ref{20.14x}) as 
\begin{equation}
   E_{\mathcal{L}_{\al}}\(    \phi_{n}(\nu ) \(\wt  \psi_{ {n}} (x,r) + H(x,r)\)\)
   \end{equation}
   and use (\ref{20.14xh}) to estimate $   E_{\mathcal{L}_{\al}}\(    \phi_{n}(\nu )   H(x,r) \)$. Exactly as in the proof of Theorem \ref{theo-multivar} and    \ref{theo-6.2} we see that these are what we called error terms and   that  after  integrating with respect to $\nu$ and taking the limit as $r$ goes to zero we get
 \bea 
 &&
E_{\mathcal{L}_{\al}}\(    \phi_{n}(\nu )\,\,\wt\psi_{n}(\nu)\)\label{20.14y}\\
&&\qquad= \sum_{\pi\in \mathcal{P}_{2n,a}}\al^{c(\pi)}\int    \prod_{j=1}^{2n}u (z_{j},z_{\pi(j)})\,d\nu  (x )\,d\nu  (y).\nn
\eea 
It follows from (\ref{7.47}) and (\ref{20.14y})   that \be
E_{\mathcal{L}_{\al}}\( \(\phi_{n}(\nu )-\wt\psi_{n}(\nu)\)^{2}\)=0.
\ee \qed

\begin{corollary} \label{cor-7.1} For $\nu\in \BB^{2n}$
\be
\II_{n}(\nu)=\psi_{n}(\nu)\label{7.49}
\ee
in $L^{p}(P_{\LL_{\al} )}$.
\end{corollary}

\Proof As in the proof of (\ref{ls.10aay}) we have 
 \begin{equation}
\psi_{n}(\nu)=I_{1}(L_{n}(\nu)).\label{7.49a}
\end{equation} 
By  Theorem \ref{theo-7.1} 
\begin{equation}
\II_{n}(\nu)=\lim_{r\rar 0}\int I_{1}(L_{n}(x,r))\,d\nu(x).\label{7.49b}
\end{equation}
We show below that
\begin{equation}
 \int I_{1}(L_{n}(x,r))\,d\nu(x)= I_{1}\(\int L_{n}(x,r)\,d\nu(x)\).\label{7.49c}
\end{equation}
Since $\int L_{n}(x,r)\,d\nu(x)\to L_{n}(\nu)$ in $\DD_{2}(\mu)$ by (\ref{rilt.9introq}) and $I_{1}$ is continuous, we get 
(\ref{7.49}).

To prove (\ref{7.49c}) we   show that 
\be
E_{\mathcal{L}_{\al}}\( \(\int I_{1}(L_{n}(x,r))\,d\nu(x)-I_{1}\(\int L_{n}(x,r)\,d\nu(x)\)\)^{2}\)=0.\label{7.49da}
\ee 
Using Fubini's theorem  and   (\ref{11.23a}) we get
\bea
&&
E_{\mathcal{L}_{\al}}\( \(\int I_{1}(L_{n}(x,r))\,d\nu(x)\)^{2}\) \label{7.49db}\\
&&\qquad=\int\int E_{\mathcal{L}_{\al}}\(  I_{1}(L_{n}(x,r))I_{1}(L_{n}(y,r)) \)\,d\nu(x)\,d\nu(y)\nn\\
&&\qquad=  \al \int\int \mu\(   L_{n}(x,r) L_{n}(y,r)  \)\,d\nu(x)\,d\nu(y),\nn
\eea 
and
\bea
&&
E_{\mathcal{L}_{\al}}\( \(I_{1}\(\int L_{n}(x,r)\,d\nu(x)\)\)^{2}\)\label{7.49d}\\
&&\qquad=\al\mu  \(\int L_{n}(x,r)\,d\nu(x)\int L_{n}(y,r)\,d\nu(y)\).     \nn
\eea 
Also
\bea
&&
E_{\mathcal{L}_{\al}}\(  \int I_{1}(L_{n}(x,r))\,d\nu(x)\,I_{1}\(\int L_{n}(x,r)\,d\nu(x)\) \)\label{7.49dc}\\
&&\qquad=\int E_{\mathcal{L}_{\al}}\(  I_{1}(L_{n}(x,r))I_{1}\(\int L_{n}(y,r)\,d\nu(y)\) \) \,d\nu(x)\, \nn\\
&&\qquad=\al\int \mu  \(L_{n}(x,r)\int L_{n}(y,r)\,d\nu(y)\)\,d\nu(x).     \nn
\eea 
The statement in (\ref{7.49c}) follows by Fubini's theorem. This completes the proof of  the lemma
\qed

 	\begin{remark} {\rm

To illustrate  (\ref{10.26}), by Theorem \ref{theo-multiloop} and Corollary \ref{cor-7.1}
\begin{equation}
\wt\psi_{2}(\nu)=\psi_{2}(\nu)+\II_{1,1}(\nu)\label{10.262}
\end{equation}
and
\begin{equation}
\wt\psi_{3}(\nu)=\psi_{3}(\nu)+3\II_{2,1}(\nu)+\II_{1,1,1}(\nu).\label{10.263}
\end{equation}
Consider 
\be 
\II_{1,1}(\nu)=\lim_{r\rar 0}\int I_{2}(L_{1}(x,r),L_{1}(x,r))\,d\nu (x).
\label{} 
\ee
Ignoring all limits and renormalization terms,  this involves the process
\bea
&&
\int L_{1}(x,r) (\om)L_{1}(x,r) (\om') \,d\nu (x)\label{}\\
&&\qquad=\int \int_{0}^{\ff}\int_{0}^{\ff} f_{r}(\om_{s}-x)\,f_{r}(\om'_{t}-x) \,ds\,dt\,d\nu (x). \nonumber
\eea
In the  limit, as $r$ goes to zero, this is an intersection local time for the two  paths $\om$ and $\om'$. Thus we see that  $ \II_{1,1}(\nu)$
involves intersections between different paths in the loop soup. Similar remarks apply to all $\II_{n_{1},\ldots, n_{l}}(\nu)$.   This should be investigated further. }
\end{remark}

\begin{remark} {\rm  We have often  commented that $\wt\psi_{2n}(\nu)$ is a Gaussian chaos, the 2$n$-th Wick power, when the potential of the underlying L\'evy process is symmetric and $\al=1/2$. A similar characterization for  $ \psi_{2n}(\nu)$ would be interesting.

}\end{remark}

We next obtain an isomorphism theorem  for $\wt\psi_{n}$. To recall the meaning of the notation it may be useful to look at Section \ref{sec-IT}.

\begin{lemma} For  $\nu\in \mathcal{B}_{2(j+k)}(R^{d})$
\begin{equation}
\(\wt\psi_{j}\times  L_{k}\) (\nu)=\lim_{r\rar 0}\int  \wt\psi_{j}(x,r) L_{k}(x,r) \,d\nu(x)\label{6.0}
\end{equation}
exists  in $L^{p}(\mu)$ for all $p\ge 2$. 
\end{lemma}

\Proof This follows as in the proofs of similar results in 
Theorems \ref{theo-1.3}   and Theorem \ref{theo-multivar}.\qed

\begin{lemma}\label{lem-loppand} For  $\nu\in \mathcal{B}_{2n}(R^{d})$
\be
\wt \psi_{n}(\nu)(\mathcal{L}_{\al}\cup \bar\om)=\sum_{l=0}^{n}{n \choose l}\(\wt\psi_{n-l}\times L_{l})(\nu)(\mathcal{L}_{\al},\bar\om\).\label{7.62}
\ee
\end{lemma}

\Proof By (\ref{11.21}), for $g_{1},\ldots, g_{l}\in \DD_{1}(\mu)$
\bea
&&
I_{l  }(g_{1},\ldots, g_{l})(\LL_{\al}\cup \bar\om) =\sum_{D=\{i_{1},\ldots,i_{|D|}\}\subseteq [1,l]}\label{10.212}\\
&&\hspace{.8 in}(-1)^{|D^{c}|}
	 \( \sum_{  (\om_{i_{1} },\ldots,\om_{i_{|D|} })\in S_{|D|}(\LL_{\al}\cup \bar\om)  }\,\,\prod_{j\in D}g_{j }(\om_{i_{j}})\)\prod_{j\in D^{c}}\al\mu(g_{j }).\nn
	 \eea
We have   
 \begin{eqnarray}
 &&\hspace{-.5 in} \sum_{  (\om_{i_{1} },\ldots,\om_{i_{|D|} })\in S_{|D|}(\LL_{\al}\cup \bar\om)  }\,\,\prod_{j\in D} g_{j }(\om_{i_{j}})= \sum_{  (\om_{i_{1} },\ldots,\om_{i_{|D|} })\in S_{|D|}(\LL_{\al})  }\,\,\prod_{j\in D} g_{j }(\om_{i_{j}})
\nn\\
 &&\hspace{-.3 in}+\sum_{m=1}^{|D|} \,\, \(\sum_{  (\om_{i_{1} },\ldots,\wh\om_{i_{m} },\ldots,\om_{i_{|D|} })\in S_{|D|-1}(\LL_{\al})  }\,\,\prod_{j\in D-\{   m\}}g_{j }(\om_{i_{j}}) \)g_{m }(\bar \om),\label{10.213}
 \end{eqnarray}
where   the notation  $(a_{1},\ldots,\wh a_{m},\ldots,a_{k})$ means that we remove 
the $m$-th entry $a_{m}$ from the vector $(a_{1},\ldots ,a_{k})$.

Using (\ref{10.213}) in (\ref{10.212}) and rearranging we obtain
\bea
&&
I_{l  }(g_{1},\ldots, g_{l})(\mathcal{L}_{\al}\cup \bar\om)=I_{l }(g_{1},\ldots, g_{l})(\mathcal{L}_{\al})\label{10.214}\\
&&\hspace{2 in}+\sum_{j=1}^{l} I_{l }(g_{1},\ldots,\wh g_{j}, \ldots,  g_{l})(\mathcal{L}_{\al})\,\,g_{j }(\bar \om).\nn
\eea
The first term on the right-hand side comes from the first term on the right hand side of (\ref{10.213}).   The  term $ I_{l  }(g_{1},\ldots,\wh g_{j}, \ldots,  g_{l})(\mathcal{L}_{\al})\,\,g_{j }(\bar \om)$, on the second line of (\ref{10.214}),   contains  all terms with $g_{j }(\bar \om)$ which arise after substituting (\ref{10.213})
in (\ref{10.212}). 
 
Proceeding as in the proof of Theorem \ref{theo-7.1} this implies that
\begin{equation}
\II_{n_{1},\ldots, n_{l}}(\nu)(\mathcal{L}_{\al}\cup \bar\om)=\II_{n_{1},\ldots, n_{l}}(\nu)(\mathcal{L}_{\al})+\sum_{j=1}^{l} 
(\II_{n_{1},\ldots,\wh n_{j}, \ldots, n_{l}}\times L_{n_{j}})(\nu)(\mathcal{L}_{\al},\bar\om),\label{13.1}
\end{equation}
where $(\II_{n_{1},\ldots,\wh n_{j}, \ldots, n_{l}}\times L_{n_{j}})(\nu)$ is defined in a manner similar to (\ref{6.0}). 
The lemma now follows from (\ref{10.26}) and the fact that there ${n \choose n_{j}}$ ways to choose a part of size $n_{j}$.\qed

\bt[Isomorphism Theorem  II]  \label{theo-pjr}  
 For any   positive measures $\rho,\phi\in  \BB  _{2}(R^{d})$ and all finite measures $\nu_{j} \in \BB_{2n }(R^{d})$ , $j=1,\ldots,n$, and bounded 
measurable functions  $F$ on $R^\ff_+$,  
\bea
&&
E_{\mathcal{L}_{\al}}\int Q^{x,x}\(L_{1}(\phi)\, F \( \sum_{k=0}^{n}{n \choose k}\(\wt\psi_{(n-k)}\times  L_{k}\) (\nu_{i}) \) \)\,d\rho(x)  \nn\\
&&\qquad={1 \over \al}E_{\mathcal{L}_{\al}}\(\th^{\rho,\phi} \,F\(\wt\psi_{n}(\nu_{i})\)\).\label{6.1}
\eea
  (The notation $F( f(x_{ i}))$ is explained in the statement of Theorem \ref{theo-ljriltintro}.) 
\et

\Proof The proof is the same as the proof of   Theorem \ref{theo-ljriltintro} except we use (\ref{7.62}) in place of (\ref{palm.2b}).\qed

An alternate, combinatorial proof of this isomorphism theorem can be obtained following the techniques in \cite{MRmem}. See also \cite{Dynkin}.

 \begin{example} {\rm 
Consider (\ref{6.1}). In the simplest case, $n=2$, it is
\bea
&&
E_{\mathcal{L}_{\al}}\int Q^{x,x}\(L_{\ff}^{\phi}\, F \Big (   \wt\psi_{2}(\nu_{i})+2\(\wt \psi_{1}\times   L_{1}\) (\nu_{i}) 
+  L_{2} (\nu_{i}) \Big ) \)\,d\rho(x)\nn\\
&&\hspace{2 in}={1 \over \al}E_{\mathcal{L}_{\al}}\(\th^{\rho,\phi} \,F\(\wt\psi_{2}(\nu_{i})\)\).\label{it6.2}
\eea
}\end{example}

\begin{remark} {\rm   Starting with (\ref{defil}) we construct the processes considered in this section from $I_{n}(g_{1},\ldots,g_{n})$ which is defined in (\ref{11.21}). The relationship (\ref{11.23a}) in   Lemma \ref{cor-11.14}  plays a critical role in the results we obtain. Lemma \ref{cor-11.14} is a simple consequence of   Lemma \ref{lem-7.3}, which contains difficult combinatorics. If we give up the need of having an explicit definition of   $I_{n}(g_{1},\ldots,g_{n})$, we can easily obtain a variant of Lemma \ref{cor-11.14}   that does not require Lemma \ref{lem-7.3}. This variant determines a continuous multilinear map   $\wt I_{n}(g_{1},\ldots, g_{n})$     from $\DD^{n}_{2}(\mu)$ to $\DD_{1}(P_{\mathcal{L}_{\al}})$ that we can use in place of $I_{n}(g_{1},\ldots,g_{n})$ in all the results in this section. 	 In particular we have the following lemma:
 
\begin{lemma} \label{lem-8.5}
   Let $\mathcal{L}_{\al}$ be  a Poisson point process  on $\Om_{\De}$ with intensity measure $\al \mu$. Then for each $n$ there exists a continuous multilinear map   $\wt I_{n}(g_{1},\ldots, g_{n})$     from $\DD^{n}_{2}(\mu)$ to $\DD_{1}(P_{\mathcal{L}_{\al}})$
such that for any $N$ and $n_{j}, j=1,\ldots, N$, and  any functions  $g_{j}\in \DD_{ 2} (\mu)$, $j=1,\ldots, s_{N+1}$,   where $s_{m} =\sum_{j=1}^{m-1}n_{j}$, $m=1,\ldots, N+1$,  
\bea
&&
E_{\mathcal{L}_{\al}}\(\prod_{m=1}^{N}\wt  I_{n_{m} }(g_{s_{m}+1},\ldots, g_{s_{m}+n_{m}})  \) \label{a11.23aq}\\
&&\qquad=
\sum_{\stackrel {\cup_{l}A_{l}=[1,s_{N+1}],\,  |A_{l}|\geq 2 }{\perp \{s_{m}+1,\ldots, s_{m}+n_{m}\},\,m=1,\ldots,N }}\prod_{l} \al \mu \(    g_{ A_{l} } \).
\nn
\eea 
\el

\Proof
 Let
\begin{equation}
\DD_{1,0} =\{f\in \DD_{1}\,|\, \mu(f)=0\}.\label{11.50}
\end{equation}
We obtain  $\wt I_{n}(g_{1},\ldots, g_{n})$ as an extension of  $ \phi_{n}(g_{1},\ldots, g_{n})$  from $ \DD^{n}_{1,0}(\mu) $ to $\DD^{n}_{2}(\mu)$.  Lemma \ref{lem-7.1} implies (\ref{a11.23aq}) for $g_{1},\ldots, g_{n}\in \DD_{1,0}$.  We get that $|A_{l}|\ge 2$ because  $g_{j}\in \DD_{1,0}$.

The key observation  that allows us to make the extension is that
$\DD_{1,0}$ is dense in $\DD_{2} $.  We note in the proof  of Lemma \ref{cor-11.14}  that $\DD_{1}$ is dense in $\DD_{2} $.   Therefore, to show that  $\DD_{1,0}$ is dense in $\DD_{2} $ it suffices to show that for any $f\in \DD_{1}$ we can find a sequence $f_{n}\in \DD_{1,0}$ which converges in $L^{p}(\mu)$ for any $p\geq 2$. 
To see this let $\{D_{n}\}$ be as defined in (\ref{dm}) and let  
\begin{equation}
f_{n}=f-{\mu(f) \over \mu(D_{n})}1_{D_{n}}.\label{11.51}
\end{equation}
Clearly, $f_{n}\in \DD_{1,0}$.
Then for any $p>1$
\begin{equation}
\int |f_{n}-f|^{p}\,d\mu={\mu^{p}(f) \over \mu^{p}(D_{n})}\int 1_{D_{n}}\,d\mu={\mu^{p}(f) \over \mu^{p-1}(D_{n})}\rar 0\label{}
\end{equation}
since $\mu(D_{n})\to\ff$ as $n\rar \ff$.\qed

We can   use $\wt I_{n}$ in place of $I_{n}$ in the definition, (\ref{defil}) of $ \II_{l_{1},\ldots, l_{n}}(x,r)$.  In the proof of  Lemma \ref{lem-7.7}    instead of (\ref{ls.10aay}), we define 
\begin{equation}
\psi (x,r)=\wt I_{1}(L_{1}(x,r)).\label{als.10aay}
\end{equation} 
We continue with the obvious substitutions.  The proof of Lemma \ref{lem-loppand} is simplified since we can replace $I_{l}(g_{1},\ldots,g_{l})$  in (\ref{10.212}) by $ \phi_{l}(g_{1},\ldots,g_{l}) $.
    
    	It is not a surprise  that $\wt I_{n}$ is the same as $  I_{n}$ defined in (\ref{11.21}). To see this we first note that the first step in the proof of Lemma \ref{lem-8.5}   extends $\phi_{n}(g_{1},\ldots, g_{n})$ from $ \DD^{n}_{1,0}(\mu) $ to $\DD^{n}_{1}(\mu)$.   By Lemma \ref{lem-7.3} we know that 
$ I_{n}(g_{1},\ldots, g_{n})$ is a continuous multilinear map  from $\DD^{n}_{1}(\mu)$ to $\DD_{1}(P_{\mathcal{L}_{\al}})$. Since $ I_{n}(g_{1},\ldots, g_{n})= \phi_{n}(g_{1},\ldots, g_{n})$ for $g_{1},\ldots, g_{n}\in \DD_{1,0}$, it follows that the extension  of $\phi_{n}(g_{1},\ldots, g_{n})$ from $ \DD^{n}_{1,0}(\mu) $ to $\DD^{n}_{1}(\mu)$ is $ I_{n}(g_{1},\ldots, g_{n})$.

}\end{remark}

  \section{Poisson chaos decomposition, II}\label{sec-expP}
 
 We introduce several new processes because we find them interesting and think that they point out a direction for future research. We continue the development of the Poisson chaos decomposition considered in  Section \ref{sec-Poisson}. We define an exponential version of the Poisson chaos and an enlargement of it that gives Poissonian chaos martingales. We give  an alternate  proof of  
Lemmas \ref{lem-7.3} and \ref{lem-11.35}  using the exponential  Poisson chaos.

\subsection{Exponential Poisson chaos} 

Let  
\begin{equation}
\phi_{\exp}(g)=\sum_{n=0}^{\ff}{1 \over n!}\phi_{n}(g,\ldots, g),\label{11.6aa}
\end{equation}
where, as before, $\phi_{0}\equiv 1$. Using (\ref{N=1}) we see that
\begin{equation}
E_{\LL_{\al}}\(|\phi_{\exp}(g)|\)\leq E_{\LL_{\al}}\(\phi_{\exp}(|g|)\)=\sum_{n=0}^{\ff}{(\al\mu(|g|))^{n} \over n!}=e^{\al\mu(|g|)},\label{expmx}
\end{equation}
  so that $\phi_{\exp}(g)$ is well defined for all $g\in \DD_{1} (\mu)$ and a similar calculation without absolute values shows that
\begin{equation}
E_{\LL_{\al}}\(\phi_{\exp}(g)\)=\sum_{n=0}^{\ff}{(\al\mu(g))^{n} \over n!}=e^{\al\mu(g)}.\label{expm}
\end{equation}

  Let $g\in \DD_{1} (\mu)$. By  (\ref{N=1})
 \begin{equation}
E_{\mathcal{L}} \(\phi_{1}(|g|) \)=\al
\int  |g(\om )|\,d\mu(\om )<\ff. \label{10.4}
 \end{equation}
Therefore,   
\begin{equation}
\phi_{1}(|g|)=\sum_{  \om  \in  \LL_{\al} } |g (\om )|\label{11.1r}
\end{equation}
converges almost surely.  Consequently
\begin{equation}
 \prod_{  \om  \in  \LL_{\al} }\(1+ g (\om )\)\label{11.1s}
\end{equation}
also converges almost surely.

\begin{lemma}\label{lem-11.6} For  $f,g\in \DD_{1}(\mu)$ 
\begin{equation}
\phi_{\exp}(g)= \prod_{  \om  \in  \LL_{\al} }\(1+ g (\om )\) \label{11.1t}
\end{equation}
and
\begin{equation}
\phi_{\exp}(f)\phi_{\exp}(g)=\phi_{\exp}(f+g+fg).\label{11.7w}
\end{equation}
 \end{lemma}
 
 \Proof By distinguishing between ordered and unordered $n$-tuples,
\bea
\hspace{-.4 in}\phi_{\exp}(g)=\sum_{n=0}^{\ff}{1 \over n!}\phi_{n}(g,\ldots, g)&=&\sum_{n=0}^{\ff}{1 \over n!}\sum_{ (\om_{i_{1} },\ldots,\om_{i_{n} })\in S_{n}(\LL_{\al}) }\,\,\prod_{j=1}^{n}g (\om_{i_{j} })\label{11.6}\\
&=&\sum_{n=0}^{\ff} \sum_{ \{\om_{i_{1} }\neq \ldots \neq \om_{i_{n} }\in  \LL_{\al}\}  }\,\,\prod_{j=1}^{n}g (\om_{i_{j} }),\nn
\eea
which is exactly what one obtains by expanding the product on the right hand side of (\ref{11.1t}).

The relationship in (\ref{11.7w}) follows from (\ref{11.1t}) and the fact that 
  $(1+f)(1+g)=(1+f+g+fg)$.
\qed

 For random variables $g\in\DD_{1}(\mu)$ we define the  exponential Poisson chaos
\begin{equation}
   \EE(g):=I_{exp}(g)=\sum_{n=0}^{\ff}{1 \over n!}I_{n }(g ,\ldots, g ).\label{10.10}
   \end{equation}

\begin{lemma} Let $g\in\DD_{1}(\mu)$,  then  
\begin{equation}
     \EE(g) =\phi_{exp}(g)e^{-\al\mu(g)}=  \prod_{  \om  \in  \LL_{\al} }\(1+ g (\om )\)e^{-\al \mu(g)}.\label{10.11w}
   \end{equation}
When $g>-1$, 
\begin{equation}
\EE(g)=e^{\sum_{  \om  \in  \LL_{\al} } \log (1+g)(\om)-\al \mu(g)}.\label{11.26}
\end{equation}
In addition  
\begin{equation}
  \EE(f) \,\EE(g) =\EE(f+g+fg)e^{\al\mu(fg)}\label{10.13}
   \end{equation}
and  
 \begin{equation}
E_{\LL_{\al}}\(\EE(f)\)=1.\label{expm2}
\end{equation}
 \end{lemma}
 
 \Proof
By the definition (\ref{11.21}),
\begin{equation}
I_{n }(g ,\ldots, g )=\sum_{k=0}^{n} {n \choose k}\phi_{k}(g,\ldots, g)  (-\al \mu(g))^{n-k}.  \label{11.25qq}
\end{equation}
Therefore
\begin{eqnarray}
\EE(g)& =&  \sum_{n=0}^{\ff}{1 \over n!}\sum_{k=0}^{n} {n \choose k}\phi_{k}(g,\ldots, g)  (-\al \mu(g))^{n-k}\label{11.25}\\
&=&  \sum_{n=0}^{\ff}\sum_{k=0}^{n} {\phi_{k}(g,\ldots, g)  \over k!} \,\, {(-\al \mu(g))^{n-k}  \over (n-k)!}=\phi_{exp}(g)e^{-\al\mu(g)}. \nn
\end{eqnarray}
 The second equality in  (\ref{10.11w}) follows from this and (\ref{11.1t}).
 
Note that when $g>-1$,  this can be written as
\begin{equation}
\phi_{\exp}(g)=e^{\sum_{  \om  \in  \LL_{\al} } \log (1+g)(\om)}.\label{11.13}
\end{equation}
Using this and (\ref{11.25}) we get (\ref{11.26}).
The product formula in (\ref{10.13}) follows from (\ref{11.7w}) and (\ref{expm2})  follows from
(\ref{expm}).\qed
 
\begin{remark}\label{rem-10.1} {\rm

   We   give proofs of    Lemmas \ref{cor-11.14} and  Lemma \ref{lem-11.35} using exponential Poisson chaoses.

For the proof of Lemmas \ref{cor-11.14} we note that, as in the proof of Lemma \ref{lem-11.6},  for   $h_{1},\ldots, h_{N}\in \DD_{1}(\mu)$,  
\begin{equation}
\prod_{m=1}^{N}\phi_{\exp}(h_{m})=\phi_{\exp}\(\prod_{m=1}^{N}(1+h_{m})-1\)=\phi_{\exp}\(\sum_{B\subseteq [1,N],\,|B|\geq 1}\prod _{m\in B}h_{m}\).\label{11.57}
\end{equation}
Therefore,  by (\ref{expm})
\begin{equation}
E_{\LL_{\al}}\(\prod_{m=1}^{N}\phi_{\exp}(h_{m})\) =e^{\sum_{B\subseteq [1,N],\,|B|\geq 1}\al\mu\(\prod _{m\in B}h_{m}\)}.\label{11.58qq}
\end{equation}
Using  (\ref{11.25})
\begin{equation}
E_{\LL_{\al}}\(\prod_{m=1}^{N}\EE(h_{m})\) =e^{\sum_{B\subseteq [1,N],\,|B|\geq 2}\,\al\mu\(h_{B}\)}=\prod_{B\subseteq [1,N],\,|B|\geq 2}e^{\al \mu\(\prod_{m\in B}h_{m}\)}.\label{11.58}
\end{equation}
Set 
 $h_{m}=\sum_{j=1}^{n_{m}} z_{s_{m}+j}g_{s_{m}+j}$ where, as   in Lemma \ref{cor-11.14}, $s_{m}:=\sum_{j=1}^{m-1}n_{j}$. We consider the    coefficients of $\prod_{m=1}^{N}\prod_{j=1}^{n_{m}} z_{s_{m}+j}$ in the expansion of the  first   and   third term in (\ref{11.58}).  Equating them 
 we obtain
 \bea
&&
E_{\mathcal{L}_{\al}}\(\prod_{m=1}^{N} I_{n_{m} }(g_{s_{m}+1},\ldots, g_{s_{m}+n_{m}}) \)\nonumber\\
&&\qquad=
\sum_{\stackrel{\cup_{l}A_{l}=[1,s_{N+1}],\,  |A_{l}|\geq 2}{\perp \{s_{m}+1,\ldots, s_{m}+n_{m}\},\,m=1,\ldots,N }}\prod_{l}  \al \mu \(\prod_{j\in A_{l}} g_{j } \).
\label{11.60}
\eea 
This is fairly easy to see for the third term and the right-hand side of (\ref{11.58}). 
The coefficient of  $\prod_{m=1}^{N}\prod_{j=1}^{n_{m}} z_{s_{m}+j}$ is a sum of products of $\prod_{l}  \al \mu \(\prod_{j\in B_{l}} g_{j } \)$ where $\{B_{l}\}$ is a partition of $[1,s_{N+1}]$. The fact that each $h_{k}$ appears at most once in each term
$\mu\(\prod_{m\in B}h_{m}\)$ gives the condition denoted by $\perp \{s_{m}+1,\ldots, s_{m}+n_{m}\},\,m=1,\ldots,N $. 

 Consider the expansion of first term in (\ref{11.58}). It follows from the fact that the  functions $I_{n_{m}}$ are multilinear,  that the coefficients  of $\prod_{m=1}^{N}\prod_{j=1}^{n_{m}} z_{s_{m}+j}$ are of the form 
\begin{equation}
   \prod_{m=1}^{N} \frac{I_{n_{m} }(g_{s_{m}+1},\ldots, g_{s_{m}+n_{m}}) }{n_{m }!}.\label{10.24a}
   \end{equation}
We add up all the terms of this form taking into account the fact that the variables $g_{s_{m}+1},\ldots, g_{s_{m}+n_{m}}$ can be arranged in $n_{m}!$ and that (\ref{10.24a}) remains the same in all these arrangements, we get the first term in (\ref{11.60}).  
  
Note that Lemma \ref{cor-11.14} is Lemma \ref{lem-7.3} for $p= s_{N+1} $. This is all we need in Section \ref{sec-lspower}. The argument just given can be extended to give a proof of  Lemma \ref{lem-7.3} for $p>s_{N+1}$. 
 
 \medskip	
   For the proof of   Lemma \ref{lem-11.35} we note that by (\ref{10.11w}) and (\ref{11.57})  
 \bea
\lefteqn{
\prod_{m=1}^{N}\EE(h_{m})=\prod_{m=1}^{N}\phi_{\exp}(h_{m})e^{-\al \mu (h_{m})}\label{11.57m}}\\
&&=\phi_{\exp}\(\sum_{B\subseteq [1,N],\,|B|\geq 1}\prod _{m\in B}h_{m}\)e^{-\al \sum^{N}_{m=1}\mu (h_{m})}\nn\\
&&=\EE\(\sum_{B\subseteq [1,N],\,|B|\geq 1}\prod _{m\in B}h_{m}\)
e^{\al \sum_{B\subseteq [1,N],\,|B|\geq 1}\mu (\prod _{m\in B}h_{m})}       e^{-\al \sum^{N}_{m=1}\mu (h_{m})}\nn\\
&&=\EE\(\sum_{B\subseteq [1,N],\,|B|\geq 1}\prod _{m\in B}h_{m}\)
e^{\al \sum_{B\subseteq [1,N],\,|B|\geq 2}\mu (\prod _{m\in B}h_{m})}. \nn
\eea
Set 
 $h_{m}=\sum_{j=1}^{n_{m}} z_{s_{m}+j}g_{s_{m}+j}$ where, as   in Lemma \ref{lem-11.35},  $s_{m} =\sum_{j=1}^{m-1}n_{j}$. We consider the    coefficients of $\prod_{m=1}^{N}\prod_{j=1}^{n_{m}} z_{s_{m}+j}$ in the expansion of the  first   and last term in (\ref{11.57m}).  Equating them 
and arguing as in the first proof in this remark,
we obtain  Lemma \ref{lem-11.35}. 
  }\end{remark}

 \subsection{Extensions to martingales }

Until this point  we have considered the variable $\al$ in the Poisson intensity measure $\al\mu$ to be  fixed. We show here  that many of the processes considered in this monograph can be extended in such a way that they become martingales  in $\al$. Consider the Poisson process $\LL$ with values in 
 $R^{1}_{+}\times \Om_{\De}$ and intensity measure $\la\times\mu$, where $\la$ is Lebesgue measure.   Clearly  
 \[\LL\cap ([0,\al]\times \Om_{\De}) \qquad\mbox{and}\qquad\LL\cap ((\al,\al+\al']\times \Om_{\De})\]
 are independent Poisson processes with intensity measures $\al\mu$ and $\al'\mu$ respectively.
 
We define 
\begin{equation}
 \LL_{\al}=\LL\cap ([0,\al]\times \Om_{\De}), \label{}
\end{equation}
 and  
 \begin{equation}
\bar \LL_{\al'}  = \LL\cap ((\al,\al+\al']\times \Om_{\De}) \label{eee}.
\end{equation}
Clearly   $ \bar \LL_{\al'}\stackrel{law}{=} \LL_{\al'} $ and  $ \LL_{\al}$ and $\bar \LL_{\al'} $ are independent.  We consider  the definition in (\ref{11.0}) to depend on $\al$ and write
 \begin{equation}
\phi^{(\al)}_{n}(g_{1},\ldots, g_{n}) =\sum_{ (\om_{i_{1} },\ldots,\om_{i_{n} })\in S_{n}(\LL_{\al}) }\,\,\prod_{j=1}^{n}g_{j}(\om_{i_{j} }).\label{11.40a}
\end{equation}
  and define 
\begin{equation}
\bar\phi^{(\al')}_{n}(g_{1},\ldots, g_{n}) =\sum_{ (\om_{i_{1} },\ldots,\om_{i_{n} })\in S_{n}(\bar \LL_{\al'}) }\,\,\prod_{j=1}^{n}g_{j}(\om_{i_{j} }).\label{11.40b}
\end{equation}
   Since $\mu$ is non-atomic, $ \LL_{\al}\cap \bar \LL_{\al'}=\emptyset $.
Therefore
\be 
\phi_{n}^{(\al+\al')}(g_{1} ,\ldots, g_{n} )\label{11.40} =\sum_{A\subseteq [1,n]}\phi_{|A|}^{(\al)}(g_{j} ,j\in A)\,\,\bar \phi_{|A^{c}|}^{(\al')}(g_{j} ,j\in A^{c}). 
\ee 
We define  $I^{(\al)}_{n}$ and $\bar I^{(\al')}_{n}$ with respect to $\phi_{j}^{(\al)}$ and $\phi_{j}^{(\al')}$, $j=1,\ldots,n$ as in (\ref{11.21}).

\bl\label{lem-mart} Let $g_{j}\in\DD_{1}(\mu)$, $j=1,\ldots,n$. Then
\be 
I_{n}^{(\al+\al')}(g_{1} ,\ldots, g_{n} )\label{11.41}\\
 =\sum_{A\subseteq[1,n]}I_{|A|}^{(\al)}(g_{j} ,j\in A)\,\,\bar I_{|A^{c}|}^{(\al')}(g_{j} ,j\in A^{c}) . 
\ee 
\el

For example
\begin{equation}
I_{2}^{(\al+\al')}(g_{1} , g_{2} )=I_{2}^{(\al)}(g_{1} , g_{2} )+I_{1}^{(\al)}(g_{1}  )\bar I_{1}^{(\al')}(g_{2}  )+I_{1}^{(\al)}(g_{2}  )\bar I_{1}^{(\al')}(g_{1}  )+\bar I_{2}^{(\al')}(g_{1} , g_{2} ).\label{11.41exam}
\end{equation}

\Proof  We define  $\EE^{(\al)}_{n}$   with respect to $I_{n}^{(\al)}$   $n=1,\ldots $ as in (\ref{10.10}). By  (\ref{10.11w})
\bea
\EE^{(\al+\al')}(f)&=& \prod_{  \om  \in  \LL_{\al+\al'} }\(1+ f (\om )\)e^{-(\al+\al')\mu(f)}\label{11.61}\\
&=& \prod_{  \om  \in  \LL_{\al} }\(1+ f (\om )\) \prod_{  \om  \in  \bar \LL_{\al'} }\(1+ f (\om )\)e^{-(\al+\al')\mu(f)}\nn\\
&=& \EE^{(\al)}(f)\EE^{(\al')}(f).\nn
\eea

 We get (\ref{11.41}) by writing $f=\sum_{j=1}^{n}z_{j}g_{j}$ and matching coefficients of $\prod_{j=1}^{n}z_{j}$;  (see Remark \ref{rem-10.1}).  \qed

 We also give a direct proof of Lemma \ref{lem-mart} that only uses the material developed in Section  \ref{sec-Poisson}. Using (\ref{11.40}) and the definition of $I_{n}$ we have 
\bea
&& 
I_{n }^{(\al+\al')}(g_{1},\ldots, g_{n}) \label{11.42}\\
&&\qquad=\sum_{D\subseteq [1,n]}
\phi_{|D|}^{(\al+\al')}(g_{j} ,j\in D ) \,(-1)^{|D^{c}|}\prod_{j\in D^{c}}(\al+\al')\mu(g_{j })\nn\\
&&\qquad=\sum_{D\subseteq [1,n]} 
\sum_{A\subseteq D}\phi_{|A|}^{(\al)}(g_{j} ,j\in A)\bar \phi_{|D-A|}^{(\al')}(g_{j} ,j\in D-A) \nn\\
&&\hspace{2in}(-1)^{|D^{c}|}\prod_{j\in D^{c}}(\al+\al')\mu(g_{j }).\nn
	 \eea
Writing
\begin{eqnarray}
&&(-1)^{|D^{c}|}\,\prod_{j\in D^{c}}(\al+\al')\mu(g_{j })
\label{11.43}\\
&&=\sum_{R\subseteq D^{c}}   \((-1)^{|R|}\,\prod_{j\in R}\al \mu(g_{j })\)\((-1)^{|D^{c}-R|}\,\prod_{j\in D^{c}-R} \al'\mu(g_{j })\),   \nonumber
\end{eqnarray}
we have 
	\begin{eqnarray}
	&& \sum_{ D\subseteq [1,n] }\sum_{A\subseteq D} 
\phi_{|A|}^{(\al)}(g_{j} ,j\in A)\bar \phi_{|D-A|}^{(\al')}(g_{j} ,j\in D-A)\label{11.44} \\&&\hspace{2in}(-1)^{|D^{c}|}\prod_{j\in D^{c}}(\al+\al')\mu(g_{j })
	\nn\\
	&&  \qquad=\sum_{ D\subseteq [1,n] } 
 \sum_{\stackrel{A\subseteq D} {R\subseteq D^{c}}  } \(\phi_{|A|}^{(\al)}(g_{j} ,j\in A)\,(-1)^{|R|}\prod_{j\in R}\al \mu(g_{j })\)\nn\\
	&&\hspace{1 in}\(\bar \phi_{|D-A|}^{(\al')}(g_{j} ,j\in D-A)\,(-1)^{|D^{c}-R|}\prod_{j\in D^{c}-R} \al'\mu(g_{j })\).\nn
	\end{eqnarray}
Consider 	$T\subseteq [1,n]$ and $A\subseteq T$ and $B\subseteq T^{c}$. Let 
$D=A\cup B$ and $R=T-A$. Therefore  $T^{c}-B=D^{c}-R$  and obviously,  $B=D-A$.  
  Using this we can rewrite the   terms  on the right hand side of (\ref{11.44}) following the double sum as  
	\begin{eqnarray} 
	&&    \(\phi_{|A|}^{(\al)}(g_{j} ,j\in A)\,(-1)^{|T-A|}\prod_{j\in R=T-A}\al \mu(g_{j })\)\nn\\
	&&\hspace{.5 in}\(\bar \phi_{|B|}^{(\al')}(g_{j} ,j\in B)\,(-1)^{|T^{c}-B|}\prod_{j\in T^{c}-B} \al'\mu(g_{j })\).\label{11.44m}  
	\end{eqnarray}	
  Consequently (\ref{11.44}) is equal to
	\begin{equation}
	 \sum_{T\subseteq [1,n]} I_{|T|}^{(\al)}(g_{j} ,j\in T)\,\,\bar I_{|T^{c}|}^{(\al')}(g_{j} ,j\in T^{c}).\label{}
	\end{equation}
These equalities prove the lemma. \qed

Since $I_{|A|}^{(\al)} $ and $\bar I_{|A^{c}|}^{(\al')} $ are independent and have  zero mean, 
the next theorem  follows immediately from Lemma \ref{lem-mart}:

\bt\label{cor-mart}   For $g_{i}\in \DD_{1}(\mu)$, $i=1,\ldots,n$ the stochastic process $I_{n}^{(\al)}(g_{1} ,\ldots, \newline g_{n} )$ is an  $(E_{\LL},\FF_{\al})$ martingale. 
\et
 
The same analysis applies to many of the processes we have studied. Using the same notation as above and  Corollary \ref{cor-mart}, Theorems \ref{theo-7.1} and \ref{theo-multiloop} and Corollary \ref{cor-7.1} we get:

\bt\label{cor-mart2} For $\nu\in \BB_{2n}(R^{d})$, $\wt\psi^{(\al)}_{n}(\nu)$ and   $\psi^{(\al)}_{n}(\nu)$ are   $(E_{\LL},\FF_{\al})$ martingales and for $\nu\in \BB_{2l}(R^{d})$, $l=\sum_{i=1}^{n}l_{i}$, 
  $\II^{(\al)}_{l_{1},\ldots, l_{n}}(\nu)$ is an  $(E_{\LL},\FF_{\al})$ martingale.
\et

\section{Convolutions of regularly varying functions}  \label{sec-8}

  Let $f$ and $\bar f $ be functions on $R_{+}^{1}$.  We write   $f\approx \bar f $ if  there exists a  $C<\ff$,   
 and $0<C_{1}, C_{2}<\ff$, such that
 \begin{equation}
 C_{1}f(x)\le \bar f(x)\le C_{2}f(x), \hspace{.2 in}\forall x\geq C.\label{app.1}
 \end{equation}
We also express this as $f(|x|)\approx \bar f (|x|)$,  with $x\in R^{d}$.
We   say that $f$ is approximately regularly varying at infinity with index $\al$ if 
 $f\approx \bar f $ for some function $\bar f $ which is regularly varying at infinity with index $\al$. 
 We   say that a function  $f\geq 0$ on $R_{+}^{1}$ is   controllable if 
  \begin{equation}
   \int_{K}\( f(|\xi|)\)^{-1}\,d\xi<\ff, \label{9.4} 
   \end{equation}
 for all compact sets $K\subset R^{d}$  and 
  \begin{equation}
   \int_{R^{d}} \( f(|\xi|)\)^{-1}\,d\xi=\ff.\label{9.4w}
   \end{equation}

    \begin{lemma}\label{lem-9.1} Let $h$ and g   be   controllable functions on $R_{+}^{1}$ that are approximately regularly varying at infinity with indices $\al$ and $\bb$ respectively, where  $\al + \bb > d$ and $\max(\al,\bb) \le  d$, $d=1,2$.    
Then
 \bea 
 &&(h) ^{-1}\underset{d}{\ast}(g)^{-1}( |\xi|) =:\int_ {R^{d}} { (h (|\xi-\eta|))^{-1} (g (| \eta|))^{-1}\, { d\eta}}\label{9.25}\\
 &&\qquad\approx     (h (|\xi|))^{-1}\int_{|\eta|\le  |\xi| } (g (|\eta|))^{-1}\,d\eta +    (g  (|\xi|))^{-1}\int_{|\eta|\le  |\xi| } (h (|\eta|))^{-1}\,d\eta,  \nn
\eea
  Furthermore,    $1/((h)^{-1}\underset{d}{\ast}(g)^{-1})$,   is a   controllable  function which is approximately regularly varying at infinity with index $\al + \bb -d$.
 \end{lemma}
 
It is useful to recall some facts about the integrals of regularly varying functions. If   $f$ is a regularly varying function at infinity with index $0<\al<d$   and $(f)^{-1}$ is locally integrable, then
\begin{equation}
   \int_{|\eta|\le  |\xi| } (f (|\eta|))^{-1}\,d\eta\approx |\xi|^{d}(f (|\xi|))^{-1}\quad\mbox{as $\xi \to\ff$}.\label{8.5n}
   \end{equation}
If $f$ is  a regularly varying function at infinity with index $d$  and $(f)^{-1}$ is locally integrable, then 
\begin{equation}
    \int_{|\eta|\le  |\xi| } (f   (|\eta|))^{-1}\,d\eta\qquad\mbox{is slowly varying at infinity}\label{8.6n}
   \end{equation}
   and  
\begin{equation}
|\xi|^{d}(f  (|\xi|))^{-1}=o\( \int_{|\eta|\le  |\xi| } (f  (|\eta|))^{-1}\,d\eta\)\quad\mbox{as $\xi \to\ff$}.\label{8.7n}
   \end{equation}
   If $f$ is  a regularly varying function at infinity with index $\al>d$, then  
     \begin{equation}
   \int_{|\eta|\ge  |\xi| } (f  (|\eta|))^{-1}\,d\eta\approx |\xi|^{d}(f  (|\xi|))^{-1}\quad\mbox{as $\xi \to\ff$}.\label{cond3}
   \end{equation} 
   
\medskip	 It is well known, (see e.g.  \cite[Lemma 7.2.4]{book}), that when $\al>0$, $ ( h (|\xi|))^{-1}$ is asymptotic to a decreasing function at infinity and similarly for $  ( g (|\xi|))^{-1}$. Therefore, since we do not specify the constant in the asymptotic bounds that we give, we  assume that $  (  h (|\xi|))^{-1}$ and $ (  g(|\xi|))^{-1}$ are  decreasing for all $|\xi|\ge |\xi_{0}|$ for some $|\xi_{0}|$ sufficiently large.

  \medskip	
  \noindent{\bf  Proof of Lemma \ref{lem-9.1}}  We first assume that $h$ and $g$ are themselves  regularly varying at infinity.  
  
  Let $|\xi|\ge 2|\xi_{0}|$. To obtain    (\ref{9.25}) we take  
  \bea 
 && \int_ {R^{d}}   (h(|\xi-\eta|))^{-1}   (g (| \eta|))^{-1}\,{ d\eta}\label{7.5}\\
 &&\qquad=  \(\int_{|\eta|\le |\xi|/2}+\int_{|\xi|/2\le |\eta|\le (3/2)|\xi|}+\int_{|\eta|\ge 3/2|\xi|} \)\cdots=I+II+III\nn .
  \eea
Using the fact that $|\xi|\ge 2|\xi_{0}|$, we see that for $|\eta|\le  |\xi|/2$, 
\be 
(h (3| \xi|/2))^{-1}\le (h(|\xi-\eta|))^{-1}\le (h (| \xi|/2))^{-1}.
\ee 
Therefore, since $h^{-1}$  is regularly varying at infinity,
\begin{equation}
(h(|\xi-\eta|))^{-1} \approx   (h(|\eta|))^{-1}.
   \end{equation} 
    Consequently for $|\xi|\ge 2|\xi_{0}|$
\begin{equation} 
  I\approx \,(h(|\xi|))^{-1}  \int_{|\eta|\le  |\xi|/2} {  (g(| \eta|))^{-1}}{\, d\eta} \approx \,(h(|\xi|))^{-1}  \int_{|\eta|\le  |\xi| } {  (g (| \eta|))^{-1}}{\, d\eta} \label{7.11a}
   \end{equation}
since the integral is also  is regularly varying at infinity.  
 
 Similarly 
 \be 
   II\approx { ((g (| \xi|))^{-1}} \int_{|\xi|/2\le|\eta|\le 3 |\xi|/2} { (h(|\xi-\eta|))^{-1}} \,{  d\eta }\label{8.9}\
   \ee
   and
\bea  \int_{0\le|\eta|\le  |\xi|/2 } { (h(|\eta|))^{-1}} \,d \eta& \le & \int_{|\xi|/2\le|\eta|\le 3 |\xi|/2} {(h(|\xi-\eta|))^{-1}} \,{  d\eta }\label{7.11aa}\\
&\le  & \int_{0\le|\eta|\le 3 |\xi| } {  (h(|\eta|))^{-1}} \,d \eta .
\eea
Consequently
\be
  II\approx  { ((g (| \xi|))^{-1}} \int_{ |\eta|\le  |\xi| } { (h(|\eta|))^{-1}} \,d \eta  .\label{kko}
\ee
Using the fact that $|\xi|\ge 2|\xi_{0}|$, we see that when $|\eta|\ge  3|\xi|/2$, $(h(|\xi-\eta|))^{-1}\approx (h(|\eta|))^{-1}$. 
Therefore, by $(\ref{cond3})$,  
 \begin{equation}
   III\approx  \int_{|\eta|\ge  3|\xi|/2}  { (h(|\eta|))^{-1} (g (| \eta|))^{-1} }\,d \eta \approx |\xi|^{d}(h(|\xi|))^{-1} ((g (| \xi|))^{-1},\label{9.8f}
   \end{equation}
as $\xi \to\ff$  because  $\al+\bb>d$.  

  The  approximate   equivalence in (\ref{9.25}) follows from (\ref{7.11a}),   (\ref{kko}), and (\ref{9.8f}).  

We next show that (\ref{9.4}) holds for  $1/((h) ^{-1}\underset{d}{\ast}(g)^{-1})$, that is
  \begin{equation}
   \int_{K} \int_ {R^{d}} { (h(|\xi-\eta|))^{-1} (g (| \eta|))^{-1}\, { d\eta}}\,d\xi<\ff\quad\mbox{for all compact sets $K\subset R^{d}$}.\label{9.25ws}
   \end{equation}  
To get (\ref{9.25ws}) we first note that 
 \bea
    &&  \int_{|\xi|\le N} \int_ {|\eta|\le M} {(h(|\xi-\eta|))^{-1} (g (| \eta|))^{-1}\, { d\eta}}\,d\xi\label{8.16}\\
      &&\qquad\le  \int_{|\xi|\le N+M}   (h(|\xi|))^{-1}  \,d\xi\int_ {|\eta|\le M}(g (| \eta|))^{-1}\, { d\eta}<\ff\nn
   \eea
by (\ref{9.4}).  In addition by  (\ref{9.8f}) if $M\geq 3N/2$, for $N$ sufficiently large
  \begin{equation}
   \int_{|\eta|\ge  M} {|\eta|^{d-1}}{ (h(|\eta|))^{-1} (g (| \eta|))^{-1} }\,d|\eta|<\ff. \label{9.8}
   \end{equation}
   
   That (\ref{9.4w}) holds for $1/((h )^{-1}\underset{d}{\ast}(g)^{-1})$ follows easily from (\ref{9.25}) and (\ref{9.4w}) applied to either $(h )^{-1}$ or $(g)^{-1}$.
   This concludes the proof of our Lemma when $h$ and $g$ are themselves  regularly varying at infinity.
   
   Suppose that   $h$ and $g$ are  controllable and approximately regularly varying functions at infinity with respect to functions $\bar h$ and $\bar g$ that  regularly varying at infinity, and let   $C$  in (\ref{app.1}) be large enough so that (\ref{app.1}) holds for both pairs $h,\bar h$ and $g, \bar g$. It is clear that for $|x|\ge C$ we can pass freely between $h$ and $\bar h$ and $g$ and $ \bar g$.

 In (\ref{7.11a}), (\ref{7.11aa}) and (\ref{8.16})   the integrals are over compact sets.  These integrals remain bounded  by the hypothesis   (\ref{9.4}).  \qed

Let h  be  a controllable function that is approximately regularly varying at infinity with index $\al$. We define
   \begin{equation}
   (H (|\xi|))^{-1}= \int_{|\eta|\le  |\xi| } (h(|\eta|))^{-1}\,d\eta.\label{defH}
   \end{equation}
   When $\al=d$, $ (H (|\xi|))^{-1}$ is a   controllable function that is approximately  slowly varying, and by hypothesis we have that   $\lim_{|\xi|\to\ff}H ^{-1}(|\xi|)=\ff$.
   
\begin{remark} {\rm \label{rem-8.1}
Note that by (\ref{8.5n})  when $d/2<\al<d$, $  (H (|\xi|))^{-1}\approx |\xi|^{d} (h(|\xi|))^{-1}$ as $|\xi|\to\ff$. Therefore, if $h$ and $g$ are both  controllable functions that are approximately  regularly varying at infinity with indices less than $d$, the right-hand side of (\ref{9.25}) is asymptotic  to $C|\xi|^{d} (h(|\xi|))^{-1}((g (| \xi|))^{-1}$ as $|\xi|\to \ff$. If $h$  is a controllable function that is approximately  
regularly varying at infinity with index $d$, and $g$   is a controllable function that is approximately  regularly varying at infinity  index less than $d$, then  the right hand side of (\ref{9.25}) is asymptotic  to $C'((g (| \xi|))^{-1} (H (|\xi|))^{-1}$ as $|\xi|\to \ff$, and similarly with $h$ and $g$ interchanged. If $h$ and $g$ are both  controllable functions that are approximately   regularly varying at infinity with index  $d$ the situation is less clear. However, whenever we are in this situation in what follows, one of the terms on the right-hand side of   (\ref{9.25}) will be larger than the other, as $|\xi|\to \ff$, so, obviously,   the right-hand side of   (\ref{9.25}) is asymptotic to  the larger term.
 }\end{remark}

   For all $n\geq 2$ let
\begin{equation}
\wt\theta_{n}(|\xi|):= \overset{\mbox{$n$-times}}
 {\overbrace {(h)^{-1} \underset{d}{\ast} (h)^{-1}\underset{d}{\ast} \cdots  \underset{d}{\ast} (h)^{-1}}}  (|\xi|).    \end{equation}
 
 \begin{lemma}\label{lem-9.2}  Let h  be  a controllable function which is approximately regularly varying at infinity with index   $\,d(1-\frac{1}{k})<\al\le d$, $k\ge 2$. Then  
 \begin{equation}
       \wt\theta_{k}(|\xi|)\approx     (h(|\xi|))^{-1} (H (|\xi|))^{-(k-1)},\label{9.12w}
   \end{equation}
    and     $1/\wt\theta_{k}$,   is a controllable function which is approximately regularly varying at infinity with index $k\al -(k-1)d$.
  \end{lemma}
 
 \Proof   For $k=2$ this lemma is just Lemma \ref{lem-9.1} with $g^{-1}=h^{-1}$.  Assume that this lemma holds $j-1$, where $j\le k$.    Consider
 \begin{equation}
   \wt\theta_{j}(|\xi|)=  \int_ {R^{d}} {(h(|\xi-\eta|))^{-1} \wt \th_{j-1}} (| \eta|)\, { d\eta} .
   \end{equation}
To check that the hypotheses of  Lemma \ref{lem-9.1} are satisfied with $g=\wt\th_{j}$, we consider the 
 indices of  regular variation of $h$ and $1/\wt \th_{j-1} $.    By (\ref{9.12w}), which we assume holds for $j-1$ and  Remark \ref{rem-8.1}, the index of regular variation of  $h  H^{j-2}$ is $\wt\bb:=\al-(j-2)(d-\al)$.  Using (\ref{8.34}) it is easy to see that $\al$ and $\wt\bb$ satisfy  the hypotheses of  Lemma \ref{lem-9.1}. Therefore, by Lemma \ref{lem-9.1}
  \bea
  \wt\theta_{j}(|\xi|)&\approx&  (h(|\xi|))^{-1}  \int_{|\eta|\le |\xi|}  \wt\th_{j-1}(|\eta|)\,d\eta+  \wt\th_{j-1}(|\xi|)  \int_{|\eta|\le |\xi|} (h(|\eta|))^{-1}\,d\eta\nn \\
  &\approx& (h(|\xi|))^{-1} \int_{|\eta|\le |\xi|}  \wt\th_{j-1}(|\eta|)\,d\eta +  (h(|\xi|))^{-1} (H (|\xi|))^{-(j-1)}.\label{8.21}
  \eea
 Furthermore 
  \bea
   \int_{|\eta|\le |\xi|}  \wt\th_{j-1}(|\eta|)\,d\eta   &\approx &    \int_{|\eta|\le |\xi|} (h(|\eta|))^{-1} \( H ^{-1}(|\eta|) \)^{j-2}\,d\eta \\
      &\le&  (H (|\xi|))^{-(j-2)} \int_{|\eta|\le |\xi|} (h(|\eta|))^{-1}  \,d\eta,\nn
       \eea
since $H ^{-1}(|\xi|)$ is increasing. Using this in (\ref{8.21}) gives (\ref{9.12w}). \qed

 \begin{lemma}\label{lem-hatcontb}  Let h  be  a controllable function which is approximately regularly varying at infinity with index   $\,d(1-\frac{1}{2n})<\al\le d$, and assume  that
    \begin{equation}
  \int    \wt \th_{ 2n} (| \la|)\,|\hat \nu (\la)|^{2}\,d\la<\ff,\label{B.6}
 \end{equation} 
  for some finite   measure $\nu$.
  Then for   any $1\leq k\leq 2n$ 
     \begin{equation}
  \int    \wt \th_{ k} (| \la|)\,|\hat \nu (\la)|^{2}\,d\la<\ff,\label{B.6a}
 \end{equation} 
 and  
  \begin{equation}
  \int   |1-e^{iz\cdot \la}|^{2}\,\, \wt \th_{ k} ( |\la |)\,|\hat \nu (\la)|^{2}\,d\la = o\( H ( 1/|z|)  \)^{2n-k},\label{B.7qq}
  \end{equation}
   as $|z|\rar 0$.
  \end{lemma}

  \Proof   By   Lemma \ref{lem-9.2},  (\ref{9.4}) holds with $f^{-1}=\wt\theta_{k}$. Using this and  the fact that   $|\hat\nu(\la)|\le |\hat\nu(0)|<\ff$, we see  that for any $N$,
   \begin{equation}
  \int_{|\la|\le N}    \wt \th_{ k} (| \la|)\,|\hat \nu (\la)|^{2}\,d\la<\ff.\label{B.6aa}
 \end{equation} 
   By (\ref{9.12w}), for $N$ sufficiently large, $\wt \th_{ k} (| \la|)\le \wt \th_{ 2n} (| \la|)$, for all $|\la|\ge N$. Using   (\ref{B.6aa}) and  (\ref{B.6}) we get (\ref{B.6a}).
  
Now since (\ref{B.7qq})  is obvious when $k=2n$   we   assume that   $1\leq k\leq 2n-1$.   
For any positive number $M$ and $1/|z|>M$ we write 
  \begin{eqnarray}
  && \int   |1-e^{iz\cdot \la}|^{2}\,\,  \wt\th_{ k} (| \la|)\,|\hat \nu (\la)|^{2}\,d\la
  \label{B.2s}\\
  &&\qquad = \int_{|\la|\leq M} +\int_{M\le |\la|\leq 1/|z|}+\int_{|\la|\geq 1/|z|}  \nn\\
  &&\qquad =I + II +III. \nonumber
  \end{eqnarray}
It follows from   (\ref{B.6a})  that 
\begin{equation}
   I\le C_{M}\, |z|^{2} \label{1.4new}
   \end{equation}
  for some constant depending on $M$.
  
    By (\ref{9.12w}) 
      \begin{equation}
\frac{  \wt\th_{ k} (  |\la|)}{  \wt\th_{ 2n} ( |\la|)}\approx  \( H (|\la|)  \)^{2n-k}. \label{8.41}
   \end{equation}
   Since this is    decreasing for all $|\la|$ sufficiently large, we have that for all $|z|$ sufficiently small,
    \bea
   III&\le&4\int_{|\la|\ge 1/|z|}    \wt\th_{ k} ( |\la|)\,|\hat \nu (\la)|^{2}\,d\la\label{9.17} \\
   &\le& \nn\frac{ 4\wt\th_{ k} ( 1/|z|)}{  \wt\th_{ 2n} (1/|z|)} \int_{|\la|\geq 1/|z|}   \,\,  \wt\th_{ 2n} (| \la|)\,|\hat \nu (\la)|^{2}\,d\la=o\(  H( 1/|z|)  \)^{2n-k}.
   \eea
  For $0<\de<2$ we write
  \begin{equation}
   |1-e^{iz\cdot \la}|^{2}\le C  |1-e^{iz\cdot \la}|^{2-\de}|z|^{\de}|\la|^{\de}.
   \end{equation}
  We choose $\de$ so that 
  \begin{equation}
  \frac{  |\la|^{\de}\wt\th_{ k} (|\la|)}{  \wt\th_{ 2n} (|\la|)}\label{9.15}
   \end{equation}
is regularly varying with a strictly positive index, i.e. $\de>(d-\al)(2n-k)$.    That this is possible follows from   the bounds on $\al$ which imply that 
\be
(d-\al)(2n-k)<d(1-\frac{k}{2n}).\label{8.36}
\ee

  We choose $M$ so that $  \wt\th_{ 2n} (|\la|)>0$ for $\la\ge M$.  
Using (\ref{9.15}) we see that  for all $|z|$ sufficiently small 
  \bea 
   II&\le &\int_{M\le |\la|\le 1/|z|}   |1-e^{iz\cdot \la}|^{2-\de} \frac{ \,|z|^{\de}  |\la|^{\de}\wt\th_{ k} (|\la|)}{  \wt\th_{ 2n} (|\la|)} \wt\th_{ 2n} ( \la)\,|\hat \nu (\la)|^{2}\,d\la\label{8.33q}\\
   &\le &\frac{|z|^{\de} ( 1/|z|)^{\de}  \wt\th_{ k} ( 1/|z|)}{  \wt\th_{ 2n} (1/|z|)}\int   |1-e^{iz\cdot \la}|^{2-\de} {  \wt\th_{ 2n} (|\la|)} \,|\hat \nu (\la)|^{2}\,d\la\nn\\
   &=&o\(  H( 1/|z|)  \)^{2n-k},\nn
 \eea
  by the dominated convergence theorem. 
 
It follows from (\ref{8.36}) that 
 \begin{equation}
   |z|^{2}=o\( H ( 1/|z|)  \)^{2n-k}. \label{8.34q}
   \end{equation} 
 Combining (\ref{1.4new}), (\ref{9.17}), (\ref{8.33q}) and (\ref{8.34q}) the proof is complete.
   \qed 
   
 The next lemma is a  variation of  Lemma \ref{lem-hatcontb}.  
 
    \begin{lemma}\label{lem-hatbound}  Let   $h $   be as in Lemma 
    \ref{lem-hatcontb} and 
let   
\be \vth_{k} (z,\la):= \int |1-e^{iz\cdot \la_{1}}|^{2 }\,h^{-1} (|\la_{1}|) \wt\th_{ k-1}(|\la-\la_{1}|)\,d\la_{1}.
\ee
If (\ref{B.6}) holds, 
then for   any $2\leq k\leq 2n$  
  \begin{equation}
  \int    \vth_{k} (z,\la)  \,|\hat \nu (\la)|^{2}\,d\la= o\( H ( 1/|z|)  \)^{2n-k}.\label{B.7k}
  \end{equation}
  \end{lemma}
  
  \Proof Consider   
  \bea  
&&  \int \int_{|\la_{1}|\le M} \,h^{-1} (|\la_{1}|) \wt\th_{ k-1}(|\la-\la_{1}|)\,d\la_{1}\,|\hat \nu (\la)|^{2}\,d\la\label{8.33}\\
&&\qquad=  \int_{|\la_{1}|\le M} \,h^{-1} (|\la_{1}|)\(\int \wt\th_{ k-1}(|\la-\la_{1}|) \,|\hat \nu (\la)|^{2}\,d\la\)\,d\la_{1}\nn 
   \eea
For  $|\la_{1}|\le M$
\begin{equation}
   \int_{|\la|\le 2M} \wt\th_{ k-1}(|\la-\la_{1}|) \,|\hat \nu (\la)|^{2}\,d\la\le   |\hat \nu (0)|^{2}\int_{|\la|\le 3M} \wt\th_{ k-1}(|\la |) \,d\la=C_{M}<\ff
   \end{equation}
   because $1/\wt\th_{k-1}$ is a controllable function. In addition by   (\ref{B.6a}),  
  \begin{equation}
   \int_{|\la|\ge 2M} \wt\th_{ k-1}(|\la-\la_{1}|) \,|\hat \nu (\la)|^{2}\,d\la\le   \int  \wt\th_{ k-1}(|\la |)  |\hat \nu (\la)|^{2}\,d\la=C'_{M}<\ff.\label{8.32}
   \end{equation} 
  (Here we use the fact that $|\la-\la_{1}|\ge |\la|/2$ and choose $M$  sufficiently large  so that  $\wt\th_{ k-1}(|\la |) $ is decreasing and $\wt\th_{ k-1}(|\la |/2) \le c\,\wt\th_{ k-1}(|\la |) $.)
  Therefore, by (\ref{9.4})   \begin{equation}
     \int \int_{|\la_{1}|\le M} |1-e^{iz\cdot \la_{1}}|^{2 } \,(h(|\la_{1}|))^{-1} \,\wt\th_{ k-1}(|\la-\la_{1}|)\,d\la_{1}\,|\hat \nu (\la)|^{2}\,d\la\le C''_{M}|z^{2}|.\label{8.35}
   \end{equation}
  For $z<1/M$ we write
 \begin{eqnarray}
 &&\int_{ |\la_{1}|\ge M} |1-e^{iz\cdot \la_{1}}|^{2 }\,(h(|\la_{1}|))^{-1} \, \wt\th_{ k-1}(\la-\la_{1})\,d\la_{1}
 \label{qwer}\\
 &&\qquad= \int_{M\le|\la_{1}|\le 1/|z|}+\int_{ |\la_{1}|\ge 1/|z|}  := I(z,\la)+II(z,\la).\nonumber
 \end{eqnarray}

  Choose  $(d-\al)(2n-k)< \de<2$ as in (\ref{9.15}), and note that, as $|\la_{1}|\to \ff $,  
  \begin{equation}
\frac{|\la_{1}|^{\de}\,\,(h(|\la_{1}|))^{-1} \,}{\wt\th_{ 2n-k+1}( \la_{1})}\sim |\la_{1}|^{\de}\(H(|\la_{1}|)\)^{2n-k}=:F(|\la_{1}|)   . \label{8.37mm}\end{equation}
We have 
   \be    I(z,\la)
   \leq  C|z|^{\de}\int_{M\le |\la_{1}|\le 1/|z|} |1-e^{iz\cdot \la_{1}}|^{2 -\de} F(|\la_{1}|)\wt\th_{ 2n-k+1}( \la_{1})\wt\th_{ k-1}(\la-\la_{1})\,d\la_{1}\nn.
\ee 
  $F(|\la_{1}|)$ is a regularly varying function at infinity with a strictly positive index, and we choose $M$ above so that we can take $F(|\la_{1}|)$   to be  increasing for $|\la_{1}|\ge M$.  Consequently 
\bea
\lefteqn{ I(z,\la)}\label{8.38}\\
  && \le \(  H ( 1/|z|)  \)^{2n-k}\int  |1-e^{iz\cdot \la_{1}}|^{2 -\de}  \wt\th_{ 2n-k+1}( \la_{1})\wt\th_{ k-1}(\la-\la_{1})\,d\la_{1}.\nn
   \eea
   Consider
   \begin{equation}
   \int    I(z,\la) \,|\hat \nu (\la)|^{2}\,d\la.
   \end{equation}
   Since 
   \begin{equation}
   \int    \wt\th_{ 2n-k+1}( \la_{1})\wt\th_{ k-1}(\la-\la_{1})\,d\la_{1}=\wt\th_{2n }( \la )\label{8.40}
   \end{equation}
it follows from (\ref{B.6}) and  dominated convergence  theorem that 
\begin{equation}
 \int    I(z,\la) \,|\hat \nu (\la)|^{2}\,d\la = o \(  |z|^{d} h ( 1/|z|)  \)^{2n-k}. \label{9.31}
   \end{equation}
Also
\be
    II(z,\la)\le 4\int_{|\la_{1}|>1/|z|} (h(|\la_{1}|))^{-1} \, \wt\th_{ k-1}(\la-\la_{1})\,d\la_{1}.\label{9.31j} 
  \ee
By (\ref{9.12w}) 
\begin{equation}
   (h(|\la_{1}|))^{-1} \, \sim \(   H(|\la_{1}|)\)^{2n-k} {\wt\th_{ 2n-k+1}( \la_{1})}.\label{9.31a}
   \end{equation}
Since $  H(|\la_{1}|)$ is deceasing in $|\la_{1}|$  we see that  
\bea
   II(z,\la) 
 &\le&  \( H ( 1/|z|)  \)^{2n-k}\int_{|\la_{1}|> 1/|z|}     \wt\th_{ 2n-k+1}( \la_{1})\wt\th_{ k-1}(\la-\la_{1})\,d\la_{1} \nn\\
  &=&   \( H ( 1/|z|)  \)^{2n-k}\int 1_{|\la_{1}|> 1/|z|}     \wt\th_{ 2n-k+1}( \la_{1})\wt\th_{ k-1}(\la-\la_{1})\,d\la_{1} \nn.
   \eea
As above it follows from (\ref{B.6}) and  dominated convergence  theorem that
 
\begin{equation}
 \int    II(z,\la) \,|\hat \nu (\la)|^{2}\,d\la = o \( H ( 1/|z|)  \)^{2n-k}.\label{9.35}
   \end{equation}
The estimates in (\ref{9.31}) and (\ref{9.35}) give (\ref{B.7k}).\qed

The next lemma gives inequalities that are used to calculate the rate of growth of chain and cycle functions.

 \bl\label{lem-chat-h} Let $h$   be controllable function that is  regularly varying at infinity with index $d/2<\al\le d$ and let $f$ be a smooth function of compact support. Then for all $r>0$ sufficiently small   \begin{equation}
   \int  (h(|x| ))^{-1} |\hat f(rx)|\,dx\le C  {(H(1/r))^{-1}}\label{9.46}
   \end{equation}
   and
   \be    
 \int  \int  (h(|s-y| ))^{-1}  (h(|y| ))^{-1}  \,d y \,|\hat f(rs)|^{2}  \,d s  \le  C(H(1/r))^{-2}\label{8.70} .
   \ee  
 \el
 
 \Proof     
 The function $\hat f$ is rapidly decreasing. We choose $K$ so that for $|x|>K$, $|\hat f(x)\le |x|^{-p} $, for some   $p>d$. We write   
 \bea
  \int (h(|x| ))^{-1} |\hat f(rx)|\,dx&=&\int_{|x|\le K}  \label{8.54}  +\int_{K<|x|\le K/r}  +\int_{|x|>K/r} \\
 &:=& I+II+III.\nn
   \eea
   Since $ h^{-1}$ is locally integrable  and $|\hat f(x)|\le \hat f(0)=1$,     $I\le C  $, for some constant $C$ depending on   $K$. In addition
  \begin{equation}
  II\le \int_{K<|x|\le K/r} (h(|x| ))^{-1}\,dx\le C_{1}  (H(1/r))^{-1}  \label{8.55}
   \end{equation}
   for all  $r\le r_{0}$ sufficiently small. Lastly
   \begin{equation}
   III\le\frac{1}{r^{p}} \int_{|x|>K/r}\frac{ (h(|x| ))^{-1} }{|x|^{p}} \,dx\le C_{ 2}\frac{(h(1/r))^{-1}}{r^{d}}\le C_{ 3} (H(1/r))^{-1} \label{8.56}
  \end{equation}
by (\ref{8.7n}). Thus we get (\ref{9.46})

  To obtain (\ref{8.70}) we first write
\bea
&&  \int  \int  (h(|s-y| ))^{-1}  (h(|y| ))^{-1}  \,d y \,|\hat f(rs)|^{2}  \,d s  \label{8.70a}\\ 
&&\qquad \nn=\int_{|u|\le K} \HH(u)\,du+\int_{K<|u|\le K/r} \HH(u)\,du\int_{|u|\ge K'r} \HH(u)\,du,
   \eea 
 where 
 \begin{equation}
   \HH(u)=  \int   (h(|u-w| ))^{-1}   (h(|w| ))^{-1}  \,d w\,|\hat f(ru)|^{2} .\label{8.60q}
   \end{equation}
   By Lemma \ref{lem-9.1}
   \begin{equation}
   \int_{|u|\le K} \HH(u)\,du\le C.
   \end{equation}
   Also, by (\ref{9.12w}) and the fact that $(H(|u|))^{-1}$ is increasing and regularly varying
   \bea
   \int_{K<|u|\le K/r} \HH(u)\,du&\le & C\int_{K<|u|\le K/r}    (h(|u | ))^{-1} (H(|u|))^{-1} \,d u  \nn\\
    &\le& C'(H(1/r))^{-1}  \int_{K<|u|\le K/r}     (h(|u | ))^{-1} \nn  \,d u \\
     &\le& C'(H(1/r))^{-2} .
   \eea
  Substituting the upper bound for $|\hat f(ru)|^{2}$,     and using the fact that $2\al>d$, which implies that $ (h(|u | ))^{-1}  (H(|u|))^{-1}$ is decreasing for large $|u|$, we obtain
  \bea
   \int_{ |u|\ge K/r} \HH(u)\,du&\le &\frac{1}{r^{p}} \int_{ |u|\ge K/r}      (h(|u | ))^{-1}  (H(|u|))^{-1} \frac{1}{|u|^{p}} \,d u\label{8.63a} \\
 &\le &C\frac{(h(1/r))^{-1}(H(1/r))^{-1} }{r^{p}} \int_{ |u|\ge K/r}        \frac{1}{|u|^{p}} \,d u\nn\\
   &\le& C'  \frac{(h(1/r))^{-1}(H(1/r))^{-1} }{r^{d}}\le     C'(H(1/r))^{-2} .\nn
   \eea 
     Combining (\ref{8.70a})--(\ref{8.63a}) we get (\ref{8.70}).
\qed
 
\begin{remark}\label{rem-al}{\rm   In $R^{1}$ the condition in Lemma \ref{lem-9.1} that  $\al,\bb\leq 1$ seems unnatural and it is. We use it for our convenience since (\ref{9.4w}) does not hold when $\al>1$. We can handle these cases but it requires rewriting the proofs in this section and it doesn't seem worthwhile for reasons discussed at the end of   Example \ref{ex-3.1}.  }\end{remark}

  We now can  provide the bounds used in Sections \ref{sec-iltac}\label{sec-2.1} and  \ref{sec-powers}.
 Consider (\ref{sl.12e}) and set
   \begin{equation}
  (\wt H(|\xi|))^{-1}:= \int_{|\eta|\le  |\xi| } (\va _{\al}(|\eta|))^{-1}\,d\eta.\label{8.65}
   \end{equation}
   Note that when $\al<d$
   \begin{equation}
   (\wt H(1/r))^{-1}\approx  \frac1 { r^{d}\va_{\al}(1/r)} 
   \end{equation}
and, for the  functions $\va_{\al}$ that we consider, (\ref{1.5}) holds.

\medskip	
Finally, we obtain bounds for the chain and cycle functions.
 \bl\label{lem-chat} When $\hat u$ satisfies (\ref{7.8q})  
  \begin{equation}
 \mbox{ch}_{k}(r)= 
 O\(   (\wt H(1/r))^{-k}\)   \qquad\mbox{as }r\rar 0\label{ss}
  \end{equation}
   and
   \begin{equation}
 \mbox{ci}_{k}(r)= 
O\(   (\wt H(1/r))^{-k}\)   \qquad\mbox{as }r\rar 0.\label{ssci}
  \end{equation}
 \el
 
 \Proof  By Lemma \ref{lem-bddhat},  (\ref{7.8q})  and (\ref{l2.1})
 \bea
 \mbox{ci}_{k}(r)&=&\frac{1}{r^{dk}}\int\prod_{j=1}^{n} u (x_{\pi(j+1)}-x_{\pi(j)}) \prod_{i=1}^{n}\,f(x_{i}/r)\,dx_{i}\\
 &\leq&\frac{1}{(2\pi)^{dk} }\(\int |\hat f(r\la )|^{2}\,\,\wt\tau_{2}(\la)\,d\la\)^{k/2}  \label{l2.22am}= O\(   (\wt H(1/r))^{-k}\) ,  \nn
\eea
where, for the last inequality we use  (\ref{8.70}). This proves  (\ref{ssci}).

 \medskip	 The proof of (\ref{ss}) follows the proof of  Lemma \ref{lem-bddhat}, however, there are enough differences that it seems necessary to give details. We  have
  \bea
\lefteqn{  \mbox{ch}_{k}(r) =  \int   u (ry_{1},ry_{2}) \cdots u(ry_{k},ry_{k+1})   \prod_{j=1}^{k+1} f(y_{j}) \,dy_{j}, \label{B.1}}\\
  &&=\int\(\int  \prod_{l=1}^{k}e^{ir(y_{l+1}-y_{l})\cdot \la_{l}} \hat u(\la_{l}) \,d\la_{l} \)  \prod_{j=1}^{k+1} f(y_{j}) \,dy_{j}\nn\\
  &&= \int  \(\prod_{l=1}^{k} \hat u(\la_{l})\)\(\hat  f(r\la_{1})\hat  f(r\(\la_{2}-\la_{1}\))\cdots 
  \hat  f(r\(\la_{k}-\la_{k-1}\))\hat  f(r\la_{k})\) \,d\la_{l}.\nn
  \eea
  Here we take the Fourier transform of $u$ considered as a distribution in $\SS'$; (see the paragraph containing (\ref{2.8})).
 
 To estimate this last integral we note that for functions $v$, $u$ and $w_{j}$, by repeated use of the Cauchy-Schwarz Inequality,   
\bea
\lefteqn{\int\dots\int    |v(\la_{1})|\prod_{j=2}^{k}   |w_{j}(\la_{j }, \la_{j-1})| |u(\la_{n})|  \prod_{j=1}^{k}\,d\la_{j} \label{32.2q}} \\
  && 
\le   \(\int |v(\la )|^{2}\,d\la\)^{1/2} \( \int |u(\la )|^{2}\,d\la\)^{1/2} \(\prod_{j=2}^{k}   \int  \int |w  _{j} (\la,\la')|^{2}   \,d \la   \,d \la'\)^{1/2}.\nn
\eea
Take  
\be 
   v= \hat u^{1/2} (\la_{1})f(r\la_{1}),\quad u=\hat u^{1/2} (\la_{k})f(r\la_{k})
   \ee
   and  
    \be
  w_{j}= \hat u^{1/2} (\la_{j-1})\hat u^{1/2} (\la_{j})f(r(\la_{j}-\la_{j-1})),\qquad
  j=2,\ldots,k.\label{8.61}
   \ee 
    Using  (\ref{32.2q})--(\ref{8.61}) and (\ref{7.8q})  and (\ref{l2.1}) we see 
 that the last integral in (\ref{B.1}) 
   \bea
  &&  \le \(\int  (\va _{\al}(|x|))^{-1}    \hat f(rx)\,dx\) \label{8.63} \\
  &&  \qquad\qquad \(    \int  \int  (\va _{\al}(|s-y|))^{-1}  (\va _{\al}(|y|))^{-1}   \,d y \,|\hat f(rs)|^{2}  \,d s  \)^{(k-1)/2}\nn.
   \eea
Using  (\ref{9.46}) and (\ref{8.70}) we get (\ref{ssci}). \qed

\def\noopsort#1{} \def\printfirst#1#2{#1}
\def\singleletter#1{#1}
            \def\switchargs#1#2{#2#1}
\def\bibsameauth{\leavevmode\vrule height .1ex
            depth 0pt width 2.3em\relax\,}
\makeatletter
\renewcommand{\@biblabel}[1]{\hfill#1.}\makeatother
\newcommand{\bysame}{\leavevmode\hbox to3em{\hrulefill}\,}

 \def\wh{\widehat}
\def\ol{\overline}

{\footnotesize

\noindent
\begin{tabular}{lll} & \hskip20pt Yves Le Jan\\  & \hskip20pt Equipe Probabilit\'es et Statistiques \\    &
\hskip20pt Universit\'e Paris-Sud \\    & \hskip20pt B\^atiment 425 \\    & \hskip20pt 91405 Orsay Cedex, France \\   & \hskip20pt yves.lejan@math.u-psud.fr 
\end{tabular}

\bigskip
\bigskip

\noindent
\begin{tabular}{lll} & \hskip20pt Michael Marcus
     & \hskip20pt  Jay Rosen\\  & \hskip20pt Department of Mathematics 
     & \hskip20pt Department of Mathematics \\    &
\hskip20pt CIty College, CUNY
     & \hskip20pt College of Staten Island, CUNY \\    & \hskip20pt New York, NY 10031
     & \hskip20pt  Staten Island, NY
10314 \\    & \hskip20pt U.S.A.
     & \hskip20pt U.S.A. \\   & \hskip20pt 
mbmarcus@optonline.net
     & \hskip20pt  jrosen30@optimum.net
\end{tabular}

}


\begin{thebibliography}{10}


\bibitem{Bertoin} J.  Bertoin {\em Random fragmentation and coagulation processes}, Cambridge University Press, New
York,  (2006).  



 

\bibitem {Dynkin}
E. B. Dynkin,
\newblock Polynomials of the occupation field and related random
    fields, {\em J. Fcnl. Anal., 58,} (1984), 20--52.  


\bibitem {EK}
N. Eisenbaum  and  H. Kaspi,
\newblock On permanental processes, {\em Stochastic Processes and their Applications},
{\em 119},  (2009),  1401-1415.


\bibitem {FPY}
P. Fitzsimmons, J. Pitman,  and M. Yor,  
\newblock Markovian bridges: construction, Palm interpretation, and splicing.  
{\em Seminar on Stochastic Processes}, 1992, E. Cinlar and K.L. Chung and M.J. Sharpe editors, 101-134, BirkhŠuser, Boston (1993).

\bibitem{K} J. F. C. Kingman, {\em Poisson Processes}, Oxford Studies in Probability, Clarendon Press, Oxford, (2002).  

\bibitem{LL} G. Lawler and V. Limic, {\em Random Walk: A Modern Introduction}, Cambridge University Press, New
York,  (2009).

\bibitem{LF}
G. Lawler and J. Trujillo Ferreis, \newblock Random walk loop soup, \newblock {\em TAMS } 359  (2007), 565--588.

\bibitem{LW}
G. Lawler and W. Werner,  \newblock The Brownian loop soup, \newblock {\em PTRF} 44 (2004), 197--217.

\bibitem {LT} M. Ledoux  and 
M.  Talagrand, 
\newblock {\em Probability in {B}anach Spaces}.
\newblock New York: Springer-Verlag, 1991.


\bibitem{Le Jan1} Y. Le Jan, \newblock
{\em Markov paths, loops and fields. }   \'{E}cole d'\'{E}t\'{e} de Probabilit\'{e}s de Saint-Flour XXXVIII - 2008. Lecture Notes in Mathematics 2026. (2011), 
Springer-Verlag, Berlin-Heidelberg. 
 

\bibitem{LMR} Y. Le Jan, M. B.  Marcus and J.~Rosen,\newblock
{\em Permanental fields,  loop soups and continuous additive functionals.},  http://arxiv.org/pdf/1209.1804.pdf 

 \bibitem{MR96} M. B.  Marcus and J.~Rosen, Gaussian chaos and sample path properties of additive functionals of
  symmetric {M}arkov processes,
 {\it Ann. Probab.}\, {\bf 24}\, (1996), 1130--1177.

 
 \bibitem{MRmem} M. B.  Marcus and J.~Rosen,  Renormalized self-intersection local times and Wick power chaos processes, 
  \,\,{\it Memoirs of the A.M.S.},\, (1999), Number 675.
  
   \bibitem{MRcont} M. B.  Marcus and J.~Rosen,  An almost sure limit theorem for Wick powers of Gaussian differences quotients,  IMS Collections, High Dimensional Probability V: The Luminy Volume, Vol. 5, (2009), 258-272.


\bibitem{book} M. B.  Marcus and J.~Rosen, {\em Markov Processes,
Gaussian Processes and Local Times}, Cambridge University Press, New
York,  (2006).

 
  




\bibitem{VJ}  D.  Vere-Jones,
     {\em Alpha-permanents},
    New Zealand J. of Math.,  (1997), 26, 125--149.



\end{thebibliography}
\end{document}